\documentclass{article}

\usepackage{arxiv}

\usepackage[utf8]{inputenc} 
\usepackage[T1]{fontenc}    
\usepackage{hyperref}       
\usepackage{url}            
\usepackage{booktabs}       
\usepackage{amsmath}
\usepackage{amsfonts}       
\usepackage{amssymb}
\usepackage{amsthm}
\usepackage{algorithmic}
\usepackage{algorithm}
\usepackage{nicefrac}       
\usepackage{microtype}      
\usepackage{lipsum}		
\usepackage{graphicx,xcolor}
\usepackage[numbers]{natbib}
\usepackage{doi}
\usepackage{cleveref}
\usepackage{subcaption}
\usepackage{svg}
\usepackage{csquotes}
\usepackage{lineno}

\usepackage{localmacros}

\title{Curvature corrected tangent space-based approximation of manifold-valued data}

\author{ 
Willem Diepeveen
\\
	Faculty of Mathematics\\
	University of Cambridge\\
	Cambridge, UK \\
	\texttt{wd292@cam.ac.uk} \\
	\And
 Joyce Chew
 \\
	Department of Mathematics\\
	University of California, Los Angeles\\
	Los Angeles, CA 90095, USA \\
	\texttt{joycechew@math.ucla.edu} \\
    \And
 Deanna Needell
 \\
	Department of Mathematics\\
	University of California, Los Angeles\\
	Los Angeles, CA 90095, USA \\
	\texttt{deanna@math.ucla.edu} \\
}

\hypersetup{
pdfauthor={Willem Diepeveen, Joyce Chew, Deanna Needell},
pdfkeywords={First keyword, Second keyword, More},
}

\begin{document}
\maketitle

\begin{abstract}
	When generalizing schemes for real-valued data approximation or decomposition to data living in Riemannian manifolds, tangent space-based schemes are very attractive for the simple reason that these spaces are linear. An open challenge is to do this in such a way that the generalized scheme is applicable to general Riemannian manifolds, is global-geometry aware and is computationally feasible. Existing schemes have been unable to account for all three of these key factors at the same time.
    
    In this work, we take a systematic approach to developing a framework that is able to account for all three factors. First, we will restrict ourselves to the -- still general -- class of symmetric Riemannian manifolds and show how curvature affects general manifold-valued tensor approximation schemes. Next, we show how the latter observations can be used in a general strategy for developing approximation schemes that are also global-geometry aware. Finally, having general applicability and global-geometry awareness taken into account we restrict ourselves once more in a case study on low-rank approximation. Here we show how computational feasibility can be achieved 
    and propose the curvature-corrected truncated higher-order singular value decomposition (CC-tHOSVD), whose performance is subsequently tested in numerical experiments with both synthetic and real data living in symmetric Riemannian manifolds with both positive and negative curvature.
\end{abstract}




\blfootnote{Our code for the main algorithms is available at \href{https://github.com/wdiepeveen/manifold-valued-tensors}{\texttt{https://github.com/wdiepeveen/manifold-valued-tensors}} }


\section{Introduction.}
\label{sec:introduction}
Manifold-valued data are becoming ubiquitous. By now classical applications from the imaging sciences include non-linear color spaces \cite{chan2001total} such as the Chromaticity Brightness model ($\Sphere^2 \times \mathbb{R}$) and the Hue Saturation Value model ($\Sphere^1 \times \Real^2$), InSAR imaging \cite{massonnet1998radar} ($\Sphere^1$), Diffusion Tensor imaging {\cite{basser1994mr} with the manifold of positive definite symmetric $3\times 3$ matrices ($\mathcal{P}(3)$) and Electron Backscatter Diffraction imaging {\cite{adams1993orientation}} with the 3D rotation group ($\mathrm{SO}(3)$). Beyond the imaging sciences, the data sciences have popularized non-linear data embedding into Riemannian manifolds, as data in many domains have non-Euclidean features \cite{bronstein2017geometric,krioukov2010hyperbolic}. That is, besides the Euclidean embedding \cite{sammon1969nonlinear} ($\Real^\dimInd$) there has been a surge in work on spherical embeddings \cite{liu2017sphereface,wilson2014spherical} ($\Sphere^\dimInd$), hyperbolic embeddings \cite{keller2020hydra,nickel2017poincare,sonthalia2020tree,walter2004h} ($\Hyperbolic^\dimInd$), combinations of these through product manifold embedding \cite{gu2019learning,skopek2020mixed}, and embedding into other non-constant curvature spaces such as the Grassmann manifold \cite{huang2016building,cruceru2021computationally} ($\mathrm{Gr}(\dimInd,k)$) and the space of symmetric positive definite matrices \cite{tuzel2006region,cruceru2021computationally} ($\mathcal{P}(\dimInd)$).

All of the above spaces are examples of \emph{symmetric Riemannian manifolds} under their standard Riemannian metric. Arguably, there is a reason why symmetric spaces are very important to these communities. Indeed, the authors of \cite{lopez2021symmetric} argue that symmetric spaces in the setting of data embeddings are a very natural, versatile and rich class of spaces for capturing the geometry of real-world data. For the imaging sciences, symmetric spaces are particularly useful from a more practical point of view, i.e., for constructing numerically efficient models \cite{bacak2016second,bergmann2018bezier} and algorithms \cite{bergmann2016parallel,bergmann2021fenchel,diepeveen2021inexact} through inexpensive evaluation of differentials or covariant derivatives of important mappings. The latter is possible due to a key property of symmetric spaces. That is, knowing the curvature at only one point yields a large amount of information about the global geometry of the space of interest.

Surprisingly, in downstream processing and analysing of manifold-valued data -- e.g., dimension reduction for visualisation \cite{hamarneh2011perception} or for retrieving underlying structure from (principal) sub-manifolds \cite{fletcher2003statistics} -- the symmetric structure is hardly ever used. Typically, works seem to choose between methods that are either \emph{global, but manifold-specific} \cite{harandi2014manifold} or \emph{localized to a tangent space, but general} \cite{fletcher2004principal}. Focusing on the latter class of approaches, subsequently quantifying the localization error is hardly ever addressed. Even for the seminal work on principal geodesic analysis (PGA)\footnote{as a generalization for principal component analysis (PCA) to Riemannian manifolds} \cite{fletcher2003statistics,fletcher2004principal} -- which is in practice often implemented as a tangent space PCA -- an approximation error analysis 
due to curvature effects is -- to the best of our knowledge -- still open for general (symmetric) Riemannian manifolds. Existing work that attempts to take curvature into account is still too local in nature \cite{lazar2017scale}, only covers specific manifolds \cite{chakraborty2016efficient,said2007exact}, or does not yield practical algorithms \cite{sommer2014optimization}.

On the flip side, being able to use the full potential of symmetric Riemannian manifolds -- in a curvature-aware tangent space-based framework -- has two main advantages. First, it would enable one to open the doors for many more approximation methods beyond simple PCA, e.g., approaches based on classical decomposition schemes for matrices and tensors \cite{kolda2009tensor} and extensions to modern ones \cite{cai2021mode,iwen2021lower,iwen2021modewise}. Secondly, and more practically, popular modern software packages for optimization and data analysis on manifolds \cite{bergmann2022manopt,miolane2020geomstats} are designed for the general setting so that one does not have to reinvent every method for new Riemannian manifolds. 
So motivated by the above, in this work we focus on the general problem of local, but curvature-aware approximation of the multi-dimensional manifold-valued array $\Tensor \in \manifold^{\dimInd_1 \times \cdots \times \dimInd_{n}}$ from some point $\mPoint\in \manifold$ on a symmetric Riemannian manifold $\manifold$. That is, given a suitable class of approximators defined by a tangent space subset $\mathcal{S}_{\mPoint}(\manifold^{\dimInd_1\times \ldots\times \dimInd_n}) \subset \tangent_{\mPoint} \manifold^{\dimInd_1 \times \cdots \times \dimInd_{n}}$, we aim to construct a suitable approximation $\mTVector_{\mPoint} \in \mathcal{S}_{\mPoint}(\manifold^{\dimInd_1\times \ldots \times \dimInd_n})$ such that the \emph{global approximation error},
\begin{equation}
    \distance_{\manifold^{\dimInd_1 \times \cdots \times \dimInd_{n}}}(\Tensor, \exp_{\mPoint}(\mTVector_{\mPoint}))^2,
    \label{eq:manifold-approximation-error}
\end{equation}
is small, where $\exp_{\mPoint}:\tangent_{\mPoint} \manifold^{\dimInd_1 \times \cdots \times \dimInd_{n}} \to \manifold^{\dimInd_1 \times \cdots \times \dimInd_{n}}$ is the exponential mapping on the manifold $\manifold^{\dimInd_1 \times \cdots \times \dimInd_{n}}$ and $d_{\manifold^{\dimInd_1 \times \cdots \times \dimInd_{n}}}: \manifold^{\dimInd_1 \times \cdots \times \dimInd_{n}} \times \manifold^{\dimInd_1 \times \cdots \times \dimInd_{n}} \to \Real$ is the distance. In particular, we aim for a generally applicable, numerically feasible and theoretically motivated approach to the approximation of the minimizer of \cref{eq:manifold-approximation-error} that can be a framework for generalizing Euclidean schemes to Riemannian manifolds.

\subsection{Related work.} 
\label{sec:related-work}
Geometry of data has been considered in several forms. One worth mentioning is the setting of data in Euclidean space that lives on a lower-dimensional unknown manifold. Although this is strictly speaking a different setting than the problem stated in \cref{eq:manifold-approximation-error} -- where there is an explicit manifold and a specific notion of geometry --, algorithms developed for the implicit setting have to deal with similar challenges. In particular, there seems to be a trend that some notion of global geometry is necessary to build suitable algorithms. To give one example, but without going into detail, graph Laplacian-based approaches for low rank approximation \cite{zhang2004principal,zhang2012low}, non-negative matrix factorization \cite{cai2008non,guan2011manifold,huang2014robust} and low rank representation \cite{shu2017multiple} are superior to naive methods that only consider the ambient Euclidean geometry.

In the case of an explicit manifold, there is more information available than in the above setting. This allows one to work in an \emph{intrinsic} setting, where one can use several manifold mappings, rather than in an extrinsic one. It is less clear how the observations from the extrinsic case regarding necessity of global information should be translated to the intrinsic case. In the following we will briefly review developments in global and tangent space-based approaches to approximation of manifold-valued data. In particular, we will highlight several best practices and limitations and argue where there are missed opportunities. The following overview will mainly focus on the \emph{order 1 case}, i.e., when the manifold-valued data lives in $\manifold^{\dimInd_1}$ for some $\dimInd_1\in\Natural$. This is for the simple reason that to the best of our knowledge there is no existing work on tensor-decomposition approaches for manifold-valued data in the spirit of \cite{kolda2009tensor}.

\paragraph{Global approaches to approximation.}
In the intrinsic setting, it seems natural -- as Riemannian manifolds do not have a well-defined notion of an origin -- to design methods globally. For the decomposition of data into nested geodesic subspaces, numerically efficient schemes for several manifolds have been proposed in a global fashion. In particular, for nested geodesic subspaces of spheres and hyperbolic spaces \cite{tabaghi2023principal} or on polyspheres and Grassmann manifolds \cite{curry2019principal}, the manifold-specific Riemannian structure gives practical algorithms. If the constraint of having geodesic subspaces is relaxed, there is the famous nested spheres framework \cite{jung2012analysis} and the recent generalisation to Grassmann manifolds \cite{yang2021nested}.

Even though the above methods perform very well for their respective manifolds, it is hard to adapt these methods to a general setting. There have been attempts to generalize nested (non-)geodesic subspace frameworks, but these run into practical challenges. One issue is that even when a method is general, implementation relies on a manifold-specific mapping \cite{mardia2022principal}. More common is the issue that the proposed method has a hard time scaling up. In particular, it can be hard to scale up to high-dimensional manifolds due to non-convexity of the proposed subspace construction problem \cite{damon2014backwards,pennec2018barycentric}, or to scale up to multidimensional subspaces either due to construction restrictions \cite{hauberg2015principal,panaretos2014principal} or due to curse of dimensionality \cite{sommer2013horizontal,yao2016principal}.

Besides the approximation approaches outlined in the above, which are all similar in spirit to doing a low rank approximation, there have been methods trying to generalize other approximation schemes. All of these come with their own challenges in terms of generalizability, but nevertheless work well for their respective manifolds. In particular, low rank representation for the Grassmann manifold as proposed in \cite{wang2015low} and \cite{piao2020reweighted} relies on linearity even though the manifold is a quotient space, which renders the definition ambiguous and hard to generalize. Similarly, in \cite{wang2018low} low-rank representation is proposed for the space of symmetric positive definite matrices, where again the ambient linearity is used. Although using linearity is less of an issue for the latter approach, the definition depends on a specific metric, once more making generalization challenging. Another example of a numerically successful, but hard to generalize method is \cite{harandi2012sparse}, in which a kernel approach is used to generalize dictionary learning on the space of symmetric positive definite matrices. Here the design of the kernel -- reflecting the global geometry of the space -- is heavily dependent on the Riemannian manifold structure of this matrix manifold.

Finally, we would like to point out that there are exceptions to this recurring issue of global methods being unable to generalize. For example, basic multiscale methods such as \cite{rahman2005multiscale} do well in terms of general applicability and computational feasibility. However, adding more structure to such a scheme may arguably be challenging. Indeed, the authors of \cite{grohs2012definability} point out that adding extra structure such as bi-orthogonal wavelet structure comes with definability issues. 


\paragraph{Tangent space approaches to approximation.}
Perhaps unsurprisingly after the above discussion, tangent space-based methods are very popular due to their general applicability with respect to the underlying manifold.
Besides their generality from a geometric point of view, the tangent space setting is more often than not the only framework in which it is even possible to generalize approximation and decomposition schemes from the Euclidean case to the general Riemannian manifold setting. 
Indeed, even though low rank-like approaches make sense for both global versions -- as discussed above -- and tangent space versions -- e.g., the aforementioned PGA framework \cite{fletcher2003statistics,fletcher2004principal} --, other methods such as dictionary learning \cite{fu2015low,ho2013nonlinear,yin2015nonlinear} can only be defined in the tangent space. 
In other words, there is a case to make for the idea of generalizing Euclidean approaches to the tangent space setting. 

There are two main challenges that come with a tangent space-based approach: where to linearize and how to take the global structure of the manifold into account in a computationally feasible way. Regarding the former challenge, it should be mentioned that one has to be careful when choosing a point of linearization as argued by the authors of \cite{huckemann2006principal} in the case of spheres and later for shape space \cite{huckemann2010intrinsic}. However, these are issues one can deal with by being careful or passing to a framework where the point of linearization is also optimized over \cite{zhang2013probabilistic}.
Arguably, the key open problem for basically any tangent space-based method is the latter challenge, i.e., a principled way to capture global geometry without sacrificing computational feasibility. To be more concrete, even though it is occasionally possible to take global structure into account, computational feasibility tends to suffer from subsequent adaptations. Once more, a prime example is the seminal work \cite{fletcher2003statistics,fletcher2004principal} by Fletcher et al. on PGA. In particular, although the original method was proposed in an exact form that took global geometry into account, linearization 
allowed for actually computing it cheaply. Later work on solving the exact version in a general setting \cite{sommer2014optimization} or for specific manifolds \cite{chakraborty2016efficient,lazar2017scale,said2007exact} yielded slow algorithms without any convergence guarantees due to non-convexity of this geometry-aware subspace construction problem.

\subsection{Specific contributions.} 
The contributions in this article are the theory and algorithm development for computationally feasible and geometry-aware tangent space-based methods for approximating manifold-valued data living in symmetric Riemannian manifolds.

First, we would like to point out that this article is, to the best of our knowledge, the first work to address approximation of multidimensional manifold-valued arrays in a tensor decomposition setting. Furthermore, this article is, to the best of our knowledge, the first to give a provably complete account of the most important effects of the global geometry on naive tangent space-based approximation schemes on symmetric Riemannian manifolds and to propose a general way of correcting for these effects. 
Finally, we showcase how the developed theory can be used to construct a numerically feasible and geometry-aware low (multi-linear) rank approximation scheme for general symmetric Riemannian manifolds. 
More specifically, the contributions are the following.

\paragraph{Quantifying the role of curvature on symmetric Riemannian manifolds.} 
We analyze the global approximation error \cref{eq:manifold-approximation-error} on Riemannian manifolds and show that for symmetric Riemannian manifolds the leading term is fully determined by the tangent space approximation error and the curvature at the point of linearization, which is our main theoretical result. In particular, we find that, if not corrected for, negative curvature can amplify seemingly minor errors exponentially, whereas seemingly large errors on a positively curved space can be diminished.

\paragraph{Curvature corrected approximation of manifold-valued tensors.}
Based on the analysis, we propose a \emph{curvature corrected approximation error}. In addition, we show that under mild conditions, minimizing the proposed error will decrease the discrepancy to the original approximation error \cref{eq:manifold-approximation-error}. In other words, minimizing the curvature corrected approximation error will yield a low original approximation error. Moreover, we will show that the curvature corrected approximation error has several nice properties such as a lower bound on the approximation error and stability with respect to vanishing curvature effects.

\paragraph{Efficient and curvature corrected low multi-linear rank approximation.}
In a case study on low multi-linear rank approximation, we will use the developed theory and propose \emph{curvature corrected truncated higher-order singular value decomposition} (CC-tHOSVD), which boils down to a truncated higher-order singular value decomposition and a curvature correction step that consists of solving a linear system of equations. We test the proposed method on synthetic 1D and real 2D manifold-valued data and on positive and negative-curvature manifolds to showcase that this method indeed overcomes the long-standing open problems of being generally applicable to different manifolds, taking global geometry into account properly and being numerically feasible.

\subsection{Outline.}
This article is structured as follows. \Cref{sec:notation} covers several notation conventions that will help the reader to keep the various objects apart throughout the rest of the work, and introduces manifold-valued tensors, the main object of interest in this article, along with several related notions.
In \cref{sec:global-error-curvature-corrected,sec:curvature-corrected-approximation,sec:cc-lra} the core theoretical ideas for realizing generally applicable, global-geometry aware and computationally feasible manifold-valued tensor approximation are developed. In particular, in \cref{sec:global-error-curvature-corrected} we state and prove our main theoretical result (\cref{eq:general-error-symmetric spaces}), which states that on general symmetric Riemannian manifolds, the leading order term in the global approximation error \cref{eq:manifold-approximation-error} is fully determined by the tangent space approximation error and and the curvature at the point of linearization. In \cref{sec:curvature-corrected-approximation} we will show that the general result from \cref{sec:global-error-curvature-corrected} can be used to realize global-geometry aware manifold-valued tensor approximation. Next, in a case study on low multi-linear rank approximation in \cref{sec:cc-lra}, we argue that computational feasibility can also be attained, which will be tested thoroughly through several numerical experiments in \cref{sec:manifold-valued-tensors-numerics}. Finally, we summarize our findings in \cref{sec:manifold-valued-tensors-conclusions}.

\section{Notation.} 
\label{sec:notation}

\subsection{Differential and Riemannian geometry.}
\label{sec:riemannian-geometry}
Here we present some basics and notation from differential and Riemannian geometry, see \cite{boothby2003introduction,carmo1992riemannian,lee2013smooth,sakai1996riemannian} for more details. 

Let $\manifold$ be a smooth manifold. We write $C^\infty(\manifold)$ for the space of smooth functions over $\manifold$. The \emph{tangent space} at $\mPoint \in \manifold$ -- defined as the space of all \emph{derivations} at $\mPoint$ -- is denoted by $\tangent_\mPoint \manifold$ and for \emph{tangent vectors} we write $\mTVector_\mPoint$. For the \emph{tangent bundle} we write $\tangent\manifold$ and vector fields -- defined as \emph{sections} of the tangent bundle -- are written as $\vectorfield(\manifold) \subset \tangent\manifold$. Furthermore, additional structure and on manifolds and derived notions thereof are typically defined through tensor fields of the form $T~:~\vectorfield(\manifold)^n \to \vectorfield(\manifold)^m$ or $T~:~\vectorfield(\manifold)^n \to C^\infty(\manifold)$. A tensor field can -- by definition of being $C^\infty(\manifold)$-modules -- be restricted to a tangent space at each point $\mPoint \in \manifold$, for which we use the notation $T_{\mPoint}: \tangent_\mPoint \manifold^n \to \tangent_\mPoint \manifold^m$ or $T_{\mPoint}: \tangent_\mPoint \manifold^n \to \Real$.

A smooth manifold $\manifold$ becomes a \emph{Riemannian manifold} if it is equipped with a smoothly varying \emph{metric tensor field} $(\,\cdot\,, \,\cdot\,) \colon \vectorfield(\manifold) \times \vectorfield(\manifold) \to C^\infty(\manifold)$. This tensor field induces a \emph{(Riemannian) metric} $\distance_{\manifold} \colon \manifold\times\manifold\to\Real_{\geq 0}$. The metric tensor can also be used to construct a unique affine connection, the \emph{Levi-Civita connection} that is denoted by $\nabla_{(\,\cdot\,)}(\,\cdot\,) : \vectorfield(\manifold) \times \vectorfield(\manifold) \to \vectorfield(\manifold)$. 
This connection is in turn the cornerstone of a myriad of manifold mappings.
One is the notion of a \emph{geodesic}, which for two points $\mPoint,\mPointB \in \manifold$ is defined as a curve $\geodesic_{\mPoint,\mPointB} \colon [0,1] \to \manifold$ with minimal length that connects $\mPoint$ with $\mPointB$. This notion is well-defined if the manifold is \emph{geodesically connected}, i.e., any two points $\mPoint,\mPointB\in\manifold$ can be connected with a geodesic that is contained in $\manifold$.
Another closely related notion is the curve $t \mapsto \geodesic_{\mPoint,\mTVector_\mPoint}(t)$  for a geodesic starting from $\mPoint\in\manifold$ with velocity $\dot{\geodesic}_{\mPoint,\mTVector_\mPoint} (0) = \mTVector_\mPoint \in \tangent_\mPoint\manifold$. This can be used to define the \emph{exponential map} $\exp_\mPoint \colon \mathcal{G}_\mPoint \to \manifold$ as 
\[ 
\exp_\mPoint(\mTVector_\mPoint) := \geodesic_{\mPoint,\mTVector_\mPoint}(1)
\quad\text{where $\mathcal{G}_\mPoint \subset \tangent_\mPoint\manifold$ is the set on which $\geodesic_{\mPoint,\mTVector_\mPoint}(1)$ is defined.} 
\]
The manifold $\manifold$ is said to be \emph{complete} whenever $\mathcal{G}_p = \mathcal{T}_{\mPoint}\manifold$.
Furthermore, the \emph{logarithmic map} $\log_\mPoint \colon \exp(\mathcal{G}'_\mPoint ) \to \mathcal{G}'_\mPoint$ is defined as the inverse of $\exp_\mPoint$, so it is well-defined on  $\mathcal{G}'_{\mPoint} \subset \mathcal{G}_{\mPoint}$ where $\exp_\mPoint$ is a diffeomorphism. Finally, we write $\curvature (\,\cdot\,, \,\cdot\,)(\,\cdot\,) \colon \vectorfield(\manifold) \times \vectorfield(\manifold) \times \vectorfield(\manifold) \to \vectorfield(\manifold)$ for the \emph{curvature tensor field}, which can be used to define the sectional curvature $\sectionalcurvature_\mPoint: \tangent_\mPoint\manifold \times \tangent_\mPoint\manifold \to \Real$ at a point $\mPoint\in \manifold$ as
\[
    \sectionalcurvature_\mPoint( \mTVector_\mPoint,\mTVectorB_\mPoint) :=\left\{\begin{matrix}
    \frac{(\curvature_\mPoint(\mTVector_\mPoint, \mTVectorB_\mPoint) \mTVectorB_\mPoint, \mTVector_\mPoint)_\mPoint}{\|\mTVector_\mPoint\|_\mPoint^2 \|\mTVectorB_\mPoint\|_\mPoint^2 - ( \mTVector_\mPoint, \mTVectorB_\mPoint)_\mPoint^2} & \text{if  $\mTVector_\mPoint, \mTVectorB_\mPoint \in \tangent_\mPoint \manifold$ are linearly independent},  \\
    0 &  \text{otherwise}.\\
    \end{matrix}\right.
\]

\subsection{Manifold-valued tensors.} 
\label{sec:manifold-valued-tensors}

In the following we will be concerned with multi-dimensional manifold-valued arrays. To stay close to the vocabulary from the decomposition of multi-dimensional real-valued arrays, i.e., the field of tensor decompositions, we will use the following definition. 

\begin{definition}[manifold-valued tensors]
    The \emph{space of manifold-valued tensors of order $n$ and size $\dimInd_1 \times \cdots \times \dimInd_n$} is denoted as the Riemannian manifold $\manifold^{\dimInd_1 \times \cdots \times \dimInd_n}$ equipped with the standard product metric. The elements of a tensor $\Tensor \in \manifold^{\dimInd_1 \times \cdots \times \dimInd_n}$ are denoted by $\Tensor_{\sumIndA_1, \cdots, \sumIndA_n}$.
\end{definition}

The notion of a manifold-valued tensor or real-valued tensor is not to be confused with the notion of a tensor field in the setting of differential geometry. In this work, the only tensor fields will be the metric tensor field and the curvature tensor field, both of which have special names and were introduced with special notation in \cref{sec:riemannian-geometry}. To avoid any further confusion we will always talk about manifold-valued tensors and real-valued tensors and use an appropriate notation convention as well. In particular, manifold-valued tensors, real-valued tensors, real-valued matrices and vectors, points on a manifold, and scalars are denoted in different typeface for clarity. Throughout this paper, calligraphic and boldface capital letters are used for manifold-valued or real-valued tensors, boldface capital letters are used for real-valued matrices, lower boldface letters for real-valued vectors or points on a manifold, and regular letters for scalars. 

\paragraph{Distance and norm.}
Because $\Tensor\in \manifold^{\dimInd_1 \times \cdots \times \dimInd_n}$ does not live in a linear space, we do not have a manifold-valued tensor norm directly, but we do have a distance, which we inherit from the product manifold-structure 
\begin{equation}
    \distance_{\manifold^{\dimInd_1 \times \cdots \times \dimInd_n}}(\Tensor, \TensorB) := \sqrt{\sum_{\sumIndA_1, \ldots, \sumIndA_n=1}^{\dimInd_1, \cdots, \dimInd_n} \distance_{\manifold}(\Tensor_{\sumIndA_1,\ldots, \sumIndA_{n}}, \TensorB_{\sumIndA_1,\ldots, \sumIndA_{n}})^2}, \quad \Tensor, \TensorB\in \manifold^{\dimInd_1 \times \cdots \times \dimInd_n}.
\end{equation}
Furthermore, as mentioned in \cref{sec:introduction} we will develop all of the analysis in a chosen tangent space $(\tangent_\mPoint\manifold)^{\dimInd_1 \times \cdots \times \dimInd_n}$ for $\mPoint\in \manifold$. For notational convenience we will write $\tangent_{\mPoint}\manifold^{\dimInd_1 \times \cdots \times \dimInd_n}$ for the tangent space and for tangent vectors we write $\mTVector_{\mPoint} \in \tangent_{\mPoint}\manifold^{\dimInd_1 \times \cdots \times \dimInd_n}$. The tangent space $\tangent_{\mPoint} \manifold^{\dimInd_1 \times \cdots \times \dimInd_n}$ is normed through the product metric tensor, i.e., the norm of a tangent vector $\mTVector_{\mPoint} \in \tangent_{\mPoint} \manifold^{\dimInd_1 \times \cdots \times \dimInd_n}$ is given by 
\begin{equation}
    \|\mTVector_{\mPoint}\|_{\mPoint} := \Bigl(\sum_{i_1, \ldots, i_n}^{\dimInd_1, \ldots, \dimInd_n} \|(\mTVector_{\mPoint})_{i_1, \ldots, i_n}\|_{\mPoint}^2 \Bigr)^{\frac{1}{2}}
    \label{eq:tangent-vector-norm-tensor-case}
\end{equation} 

\paragraph{Pointwise manifold mappings.}
The norm in \cref{eq:tangent-vector-norm-tensor-case} is in some sense defined pointwise in $\mPoint$. In the following we will often use several manifold mappings in a similar fashion. That is, we write for geodesics
\begin{equation}
    \geodesic_{\mPoint, \Tensor}(t):= \bigl(\geodesic_{\mPoint,\Tensor_{\sumIndA_1, \ldots, \sumIndA_n}}(t)\bigr)_{\sumIndA_1, \ldots, \sumIndA_n=1}^{\dimInd_1, \cdots, \dimInd_n}.
\end{equation}
Similarly, for the exponential and logarithmic mappings we write
\begin{equation}
    \exp_{\mPoint}(\mTVector_{\mPoint} ) := \bigl(\exp_{\mPoint}((\mTVector_{\mPoint})_{\sumIndA_1, \ldots, \sumIndA_n}) \bigr)_{\sumIndA_1, \ldots, \sumIndA_n=1}^{\dimInd_1, \cdots, \dimInd_n}, \quad \text{and} \quad
    \log_{\mPoint}\Tensor := \bigl(\log_{\mPoint}\Tensor_{\sumIndA_1, \ldots, \sumIndA_n}\bigr)_{\sumIndA_1, \ldots, \sumIndA_n=1}^{\dimInd_1, \cdots, \dimInd_n}
\end{equation}
and for the distance we write
\begin{equation}
    \distance_{\manifold^{\dimInd_1 \times \cdots \times \dimInd_n}}(\mPoint, \Tensor) := \sqrt{\sum_{\sumIndA_1, \ldots, \sumIndA_n=1}^{\dimInd_1, \cdots, \dimInd_n} \distance_{\manifold}(\mPoint, \Tensor_{\sumIndA_1,\ldots, \sumIndA_{n}})^2}.
\end{equation}

\paragraph{Mode-$k$ unfolding.}
An $n$ mode manifold-valued tensor $\Tensor$ can be reshaped into a manifold-valued tensor array -- i.e., a manifold-valued tensor of order $1$ consisting of manifold-valued tensors of order $n-1$ --, in $n$ ways by unfolding it along each of the $n$ modes. The \emph{mode-$k$ unfolding} of manifold-valued tensor
$\Tensor \in \manifold^{\dimInd_1 \times \cdots \times \dimInd_n}$ is the manifold-valued array denoted by $\Tensor_{(k)} \in (\manifold^{d_1 \times \dimInd_{k-1} \times \dimInd_{k+1}\cdots \times d_n})^{d_k}$ whose entries are composed of a manifold-valued tensor obtained from $\Tensor$ by fixing all indices except for the $k$-th dimension. The mapping $\Tensor \mapsto \Tensor_{(k)}$ is called the mode-$k$ unfolding operator.

\paragraph{Mode-$k$ product.}
The \emph{mode-$k$ product} at ${\mPoint}\in \manifold$ between a tangent vector with a real-valued matrix $\times_{k}:\tangent_{\mPoint} \manifold^{\dimInd_1 \times \cdots \times \dimInd_n} \times \Real^{\dimInd_k'\times \dimInd_k}\to \tangent_{\mPoint} \manifold^{d_1 \times \dimInd_{k-1} \times \dimInd_{k}' \times \dimInd_{k+1}\cdots \times d_n}$ for $k\in\Natural$ with $k \leq n$ is defined as 
\begin{equation}
    (\mTVector_{\mPoint} \times_{k} \mathbf{A})_{i_1, \ldots, i_{k-1}, j, i_{k+1}, \ldots,  i_n} := \sum_{s=1}^{\dimInd_k} (\mTVector_{\mPoint})_{i_1, \ldots, i_{k-1}, s, i_{k+1}, \ldots,  i_n} \mathbf{A}_{j,s},
    \label{eq:mode-k-product}
\end{equation}
where $\mPoint\in \manifold$, $\mTVector_{\mPoint}\in \tangent_{\mPoint} \manifold^{\dimInd_1 \times \cdots \times \dimInd_n}$ and $\mathbf{A}\in \Real^{\dimInd_k'\times \dimInd_k}$. If we have a multiple mode product from different modes, we write $\mTVector_{\mPoint}\times_{k=t}^s \mathbf{A}_k:= \mTVector_{\mPoint} \times_t \mathbf{A}_t \times_{t+1} \ldots \times_s \mathbf{A}_s$. 

\begin{remark}
    Note that \cref{eq:mode-k-product} is well-defined as all components of $\mTVector_{\mPoint}$ live in the same tangent space $\tangent_\mPoint \manifold$.
\end{remark}


\section{Global error of tangent space-based manifold-valued tensors approximation.} 
\label{sec:global-error-curvature-corrected}

Recall the global approximation error
\begin{equation}
    \distance_{\manifold^{\dimInd_1 \times \cdots \times \dimInd_{n}}}(\Tensor, \exp_{\mPoint}(\mTVector_{\mPoint}))^2,
    \label{eq:global-approximation-error-2}
\end{equation}
from \cref{sec:introduction}. In this section we will analyze the above approximation error for a given manifold-valued tensor $\Tensor \in \manifold^{\dimInd_1 \times \cdots \times \dimInd_{n}}$, point of linearization $\mPoint\in\manifold$ and tangent vector approximation $\mTVector_{\mPoint} \in \tangent_{\mPoint}\manifold^{\dimInd_1 \times \cdots \times \dimInd_{n}}$ and make two key observations. First, we will show in \cref{sec:global-error-riemannian-manifolds} that for $\mTVector_{\mPoint}\approx \log_\mPoint \Tensor$, the dominant behavior of \cref{eq:global-approximation-error-2} can be understood through the solution of a \emph{Jacobi equation}. Subsequently, in \cref{sec:global-error-symmetric-riemannian-manifolds} we restrict ourselves to \emph{symmetric Riemannian manifolds}, on which the Jacobi equation of interest has a closed form solution. We will conclude this section with a discussion on practical consequences that naive approximation schemes will suffer from when not correcting for global geometry.


\subsection{Global error on Riemannian manifolds.} 
\label{sec:global-error-riemannian-manifolds}

As suggested above, \emph{Jacobi fields} will be the key tool for the analysis.
A vector field $J\in \mathscr{X}(\gamma)$ along a geodesic $\gamma:[0,1]\to\manifold$  is said to be a Jacobi field if it satisfies the \emph{Jacobi equation}
\begin{equation}
    \frac{\mathrm{D}^{2}}{\mathrm{d}t^{2}}J(t)+R(J(t),\dot{\gamma}(t))\dot{\gamma} (t)=0, \quad \text{for all $t\in [0,1]$,}
    \label{eq:jacobi-equation}
\end{equation}
where $\frac{\mathrm{D}}{\mathrm{d}t}$ is the covariant derivative with respect to the Levi-Civita connection and $\dot{\gamma}(t)\in \tangent_{\gamma(t)}\manifold$ the velocity tangent vector field of $\gamma$.

The (length-minimizing) geodesic, on which the Jacobi field will be defined, is $\geodesic_{\mPoint,\Tensor}$. To avoid ambiguity, i.e., for $\log_\mPoint \Tensor$ to be well-defined, these geodesics have to be the unique. This can be done by restricting the entries of the manifold-valued tensor $\Tensor_{\sumIndA_1, \ldots, \sumIndA_n} \in \exp_{\mPoint}(\mathcal{G}'_{\mPoint})$, where $\mathcal{G}'_{\mPoint}\subset \tangent_\mPoint\manifold$ as introduced in \cref{sec:riemannian-geometry}.

We are now ready to prove the following statement, which tells us that for $\mTVector_{\mPoint}\approx \log_\mPoint \Tensor$, the dominant behavior of the global approximation error can be understood through the solution of a Jacobi equation. 

\begin{lemma}[approximation error on the manifold]
\label{lem:approximation-error-Riemannian-manifold}
Let $(\manifold, (\cdot, \cdot))$ be a Riemannian manifold and $\mPoint\in \manifold$ be a point on the manifold. Furthermore, let $\Tensor\in \manifold^{d_1 \times \cdots \times d_n}$ be a manifold-valued tensor with entries satisfying $\Tensor_{\sumIndA_1, \ldots, \sumIndA_n} \in \exp_{\mPoint}(\mathcal{G}'_{\mPoint})$. Finally, let $\mTVector_{\mPoint}\in \tangent_{\mPoint} \manifold^{d_1 \times \cdots \times d_n}$ be any approximation of $\log_{\mPoint}\Tensor$ and $J\in \mathscr{X}(\gamma_{\mPoint,\Tensor})$ be the solution to the Jacobi equation along the length minimizing geodesic $\gamma_{\mPoint,\Tensor}:[0,1]\to \manifold^{d_1 \times \cdots \times d_n}$ with initial condition $J(0)= 0 \in \tangent_{\mPoint}\manifold^{d_1 \times \cdots \times d_n}$ and initial velocity $\frac{\mathrm{D}}{\mathrm{d}t}J(t)|_{t=0} = \frac{\mTVector_{\mPoint} - \log_{\mPoint}\Tensor}{\|\mTVector_{\mPoint} - \log_{\mPoint}\Tensor\|_{\mPoint}}\in \tangent_{\mPoint}\manifold^{d_1 \times \cdots \times d_n}$.

Then, the error \cref{eq:global-approximation-error-2} satisfies
\begin{equation}
    \distance_{\manifold^{d_1 \times \cdots \times d_{n}}}(\Tensor, \exp_{\mPoint}(\mTVector_{\mPoint}))^2 =  \epsilon^2 \|J(1)\|^2_{\Tensor} + \mathcal{O}(\epsilon^3), \quad \epsilon:= \|\mTVector_{\mPoint} - \log_{\mPoint}\Tensor\|_{\mPoint}.
    \label{eq:jacobi-bound-low-rank-mfld-tensor}
\end{equation}
\end{lemma}

\begin{proof}
Define the geodesic variation $\Gamma:[0,1]\times[-\epsilon, \epsilon]$ with $\epsilon := \|\mTVector_{\mPoint} - \log_{\mPoint}\Tensor\|_{\mPoint}$ as in the statement above by
\begin{equation}
    \Gamma(t,s):= \exp_{\mPoint}\Bigl(t (\log_{\mPoint}\Tensor + \frac{s}{\epsilon}(\mTVector_{\mPoint} - \log_{\mPoint}\Tensor))\Bigr).
\end{equation}
Note that $\Gamma(1,0) = \Tensor$ and $\Gamma(1,\epsilon) = \exp_{\mPoint}(\mTVector_{\mPoint})$. The corresponding Jacobi field is given by 
\begin{equation}
    J(t):= \frac{\partial}{\partial s}\Gamma(t,s)\mid_{s=0}
    \label{eq:corresponding-jacobi-field-geo-variation-thm}
\end{equation}
and it satisfies the initial condition
\begin{equation}
    J(0)= \frac{\partial}{\partial s}\Gamma(t,s)\mid_{s=0}\mid_{t=0} = \frac{\partial}{\partial s}\Gamma(0,s)\mid_{s=0} = \frac{\partial}{\partial s} \exp_{\mPoint}(0 ) \mid_{s=0} = 0
\end{equation}
and
\begin{multline}
    \frac{\mathrm{D}}{\mathrm{d}t}J(t)\mid_{t=0} = \frac{\mathrm{D}}{\mathrm{d}t} \frac{\partial}{\partial s}\Gamma(t,s)\mid_{s=0}\mid_{t=0} \overset{\text{symmetry}}{=} \frac{\mathrm{D}}{\mathrm{d}s} \frac{\partial}{\partial t}\Gamma(t,s)\mid_{t=0}\mid_{s=0} \\
    = \frac{\mathrm{D}}{\mathrm{d}s} \frac{\partial}{\partial t} \exp_{\mPoint}\Bigl(t (\log_{\mPoint}\Tensor + \frac{s}{\epsilon}(\mTVector_{\mPoint} - \log_{\mPoint}\Tensor))\Bigr)\mid_{t=0}\mid_{s=0} \\
    = \frac{\mathrm{D}}{\mathrm{d}s} \log_{\mPoint}\Tensor + \frac{s}{\epsilon}(\mTVector_{\mPoint} - \log_{\mPoint}\Tensor) \mid_{s=0} = \frac{1}{\epsilon}(\mTVector_{\mPoint} - \log_{\mPoint}\Tensor),
\end{multline}
where we used the fact that $\frac{\mathrm{d}}{\mathrm{d}t} \exp_x(tX)\mid_{t=0} = X$ in the second step.

We will show that \cref{eq:jacobi-bound-low-rank-mfld-tensor} holds through Taylor expansion of $d_{\manifold^{d_1 \times \cdots \times d_{n}}}(\Gamma(1,0), \Gamma(1,s))^2$ around $s=0$. For notational convenience we will write $\distance(\cdot, \cdot)$ for the distance rather than $d_{\manifold^{d_1 \times \cdots \times d_{n}}}(\cdot, \cdot)$. Expanding up to second order gives
\begin{multline}
    \distance(\Tensor, \exp_{\mPoint}(\mTVector_{\mPoint}))^2 = \distance(\Gamma(1,0), \Gamma(1,\epsilon))^2 
    = \distance(\Gamma(1,0), \Gamma(1,0))^2 +  \epsilon \frac{\mathrm{d}}{\mathrm{d} s} \distance(\Gamma(1,0), \Gamma(1,s))^2\mid_{s=0} \\
    + \frac{\epsilon^2}{2} \frac{\mathrm{d}^2}{\mathrm{d} s^2} \distance(\Gamma(1,0), \Gamma(1,s))^2\mid_{s=0} + \frac{1}{2}\int_0^\epsilon \frac{\mathrm{d}^3}{\mathrm{d} s'^3} \distance(\Gamma(1,0), \Gamma(1,s'))^2 (s - s')^2\mathrm{d}s'.
\end{multline}
For proving the result \cref{eq:jacobi-bound-low-rank-mfld-tensor}, we must show that the zeroth and first order terms are zero, i.e.,
\begin{equation}
    \distance(\Gamma(1,0), \Gamma(1,0))^2=0
    \label{eq:jacobi-bound-low-rank-mfld-tensor-thm-taylor0}
\end{equation}
and
\begin{equation}
    \frac{\mathrm{d}}{\mathrm{d} s} \distance(\Gamma(1,0), \Gamma(1,s))^2\mid_{s=0} = 0,
    \label{eq:jacobi-bound-low-rank-mfld-tensor-thm-taylor1}
\end{equation}
and that the second order term satisfies
\begin{equation}
    \frac{1}{2} \frac{\mathrm{d}^2}{\mathrm{d} s^2} \distance(\Gamma(1,0), \Gamma(1,s))^2\mid_{s=0} = \|J(1)\|^2_{\Tensor}.
    \label{eq:jacobi-bound-low-rank-mfld-tensor-thm-taylor2}
\end{equation}
It is clear that the remainder term is $\mathcal{O}(\epsilon^3)$.

Trivially, \cref{eq:jacobi-bound-low-rank-mfld-tensor-thm-taylor0} holds, because $\distance$ is a distance. For showing \cref{eq:jacobi-bound-low-rank-mfld-tensor-thm-taylor1} notice that
\begin{multline}
    \frac{\mathrm{d}}{\mathrm{d} s} \distance(\Gamma(1,0), \Gamma(1,s))^2 = D_{\Gamma(1,s)} \distance(\Gamma(1,0), \cdot)^2 \Bigl[ \frac{\partial}{\partial s} \Gamma(1,s) \Bigr] = \Bigl( \operatorname{grad} \distance(\Gamma(1,0), \cdot)^2 \mid_{\Gamma(1,s)}, \frac{\partial}{\partial s} \Gamma(1,s)\Bigr)_{\Gamma(1,s)} \\
    = \Bigl(- 2 \log_{\Gamma(1,s)} \Gamma(1,0), \frac{\partial}{\partial s} \Gamma(1,s)\Bigr)_{\Gamma(1,s)},
    \label{eq:towards-jacobi-bound-low-rank-mfld-tensor-thm-taylor1}
\end{multline}
where we used the chain rule in the first equality, the definition of the Riemannian gradient in the second equality and the fact that $\operatorname{grad} \distance(\mPointB,\cdot)^2\mid_{\mPointC}= -2 \log_\mPointC \mPointB$ on any Riemannian manifold. Then, evaluating \cref{eq:towards-jacobi-bound-low-rank-mfld-tensor-thm-taylor1} at $s=0$ yields \cref{eq:jacobi-bound-low-rank-mfld-tensor-thm-taylor1} as $\log_{\Gamma(1,0)} \Gamma(1,0) = 0\in \tangent_{\mPoint}\manifold^{d_1 \times \cdots \times d_{n}}$.

Finally, it remains to show that \cref{eq:jacobi-bound-low-rank-mfld-tensor-thm-taylor2} holds. From \cref{eq:towards-jacobi-bound-low-rank-mfld-tensor-thm-taylor1} we find that
\begin{multline}
    \frac{\mathrm{d}^2}{\mathrm{d} s^2} \distance(\Gamma(1,0), \Gamma(1,s))^2 = \frac{\mathrm{d} }{\mathrm{d} s} \Bigl(- 2 \log_{\Gamma(1,s)} \Gamma(1,0), \frac{\partial }{\partial s} \Gamma(1,s)\Bigr)_{\Gamma(1,s)} \\
    \frac{\mathrm{d} }{\mathrm{d} s} \Bigl(- 2 \log_{(\cdot)} \Gamma(1,0), \frac{\partial }{\partial s'} \Gamma(1,s') \mid_{s' = (\Gamma(1,\star))^{-1}(\cdot) }\Bigr)_{\Gamma(1,s)} \\
    = \frac{\partial }{\partial s} \Gamma(1,s) \Bigl(- 2 \log_{(\cdot)} \Gamma(1,0), \frac{\partial }{\partial s'} \Gamma(1,s') \mid_{s' = (\Gamma(1,\star))^{-1}(\cdot) } \Bigr)_{(\cdot)} \\
    \overset{\text{metric compatibility}}{=} \Bigl(- 2 \nabla_{\frac{\partial }{\partial s} \Gamma(1,s)}\log_{(\cdot)} \Gamma(1,0), \frac{\partial }{\partial s} \Gamma(1,s) \Bigr)_{\Gamma(1,s)} \\ + \Bigl(- 2 \log_{\Gamma(1,s)} \Gamma(1,0), \nabla_{\frac{\partial }{\partial s} \Gamma(1,s)} \frac{\partial }{\partial s'} \Gamma(1,s') \mid_{s' = (\Gamma(1,\star))^{-1}(\cdot) }\Bigr)_{\Gamma(1,s)}.
    \label{eq:towards-jacobi-bound-low-rank-mfld-tensor-thm-taylor2}
\end{multline}

Evaluating \cref{eq:towards-jacobi-bound-low-rank-mfld-tensor-thm-taylor2} at $s=0$ gives once again zero for the second term because $\log_{\Gamma(1,0)} \Gamma(1,0) = 0\in \tangent_{\mPoint}\manifold^{d_1 \times \cdots \times d_{n}}$. However, using that $(-\nabla_{X_{\mPointB}} \log_{(\cdot)}\mPointB, Y_{\mPointB})_{\mPointB} = (X,Y)_{\mPointB}$ on any Riemannian manifold, the first term in the final line of \cref{eq:towards-jacobi-bound-low-rank-mfld-tensor-thm-taylor2} reduces to
\begin{equation}
    \Bigl(- 2 \nabla_{\frac{\partial }{\partial s} \Gamma(1,s)}\log_{(\cdot)} \Gamma(1,0), \frac{\partial }{\partial s} \Gamma(1,s) \Bigr)_{\Gamma(1,s)}\mid_{s=0} = 2 \Bigl(\frac{\partial }{\partial s} \Gamma(1,s)\mid_{s=0} , \frac{\partial }{\partial s} \Gamma(1,s) \mid_{s=0}\Bigr)_{\Gamma(1,0)} \overset{\cref{eq:corresponding-jacobi-field-geo-variation-thm}}{=} 2 \|J(1)\|^2_{\Tensor}.
\end{equation}
When the above is substituted back into \cref{eq:towards-jacobi-bound-low-rank-mfld-tensor-thm-taylor2}, \cref{eq:jacobi-bound-low-rank-mfld-tensor-thm-taylor2} holds, which proves the claim.

\end{proof}


\subsection{Global error on symmetric Riemannian manifolds.} 
\label{sec:global-error-symmetric-riemannian-manifolds}

\Cref{lem:approximation-error-Riemannian-manifold} tells 
us that for $\mTVector_{\mPoint}\approx \log_\mPoint \Tensor$ the global approximation error can be understood in terms of the solution to the Jacobi equation along the length minimizing geodesic $\gamma_{\mPoint,\Tensor}:[0,1]\to \manifold^{d_1 \times \cdots \times d_n}$ with initial condition $J(0)= 0 \in \tangent_{\mPoint}\manifold^{d_1 \times \cdots \times d_n}$ and initial velocity $\frac{\mathrm{D}}{\mathrm{d}t}J(t)|_{t=0} = \frac{\mTVector_{\mPoint} - \log_{\mPoint}\Tensor}{\|\mTVector_{\mPoint} - \log_{\mPoint}\Tensor\|_{\mPoint}}\in \tangent_{\mPoint}\manifold^{d_1 \times \cdots \times d_n}$, assuming that $\|J(1)\|^2_{\Tensor} \in \mathcal{O}(1)$. In this part we will restrict ourselves to locally symmetric Riemannian manifolds, on which we can solve the Jacobi equation in closed-form. This will allow us to study the global approximation error more closely.

In the following, we use the definition of symmetric Riemannian manifolds through geodesic reflections. 

\begin{definition}[geodesic reflection]
    Let $(\manifold, (\cdot, \cdot))$ be a Riemannian manifold and $\mPoint\in \manifold$ be a point on the manifold. Furthermore, let $\mathcal{N}(\mPoint)\subset \manifold$ be a neighborhood of $\mPoint$. A mapping $\reflection_{\mPoint}: \mathcal{N}(\mPoint) \to \manifold$ is called a \emph{geodesic reflection} at $\mPoint \in \mathcal{N}(\mPoint)$ if
    \begin{equation}
        \reflection_{\mPoint}(\mPoint)=\mPoint \text { and } D_{\mPoint}\reflection_{\mPoint} (\cdot) =-\operatorname{id}_{\mPoint},
        \label{eq:geodesic-reflection-def-properties}
    \end{equation}
    where $D_{\mPoint}\reflection_{\mPoint} (\cdot):\tangent_{\mPoint}\manifold \to \tangent_{\mPoint}\manifold$ is the differential of $\reflection_{\mPoint}$.
\end{definition}

\begin{definition}[symmetric Riemannian manifold]
    A Riemannian manifold $(\manifold, (\cdot, \cdot))$ is called \emph{locally symmetric} if a geodesic reflection at any point $\mPoint \in \manifold$ is an isometry on a local neighborhood of $\mPoint$, i.e., for all $\mPointB, \mPointC \in \mathcal{N}(\mPoint)\subset \manifold$ we have
    \begin{equation}
        \distance_{\manifold}(\reflection_{\mPoint}(\mPointB), \reflection_{\mPoint}(\mPointC))=\distance_{\manifold}(\mPointB, \mPointC) .
        \label{eq:symmetry-properties}
    \end{equation}
    A locally symmetric space is said to be a \emph{(globally) symmetric} space if in addition its geodesic symmetries can be extended to isometries on all of $\manifold$. 
\end{definition}

A Riemannian manifold being locally symmetric is equivalent to a vanishing covariant derivative of the curvature tensor, i.e., $\nabla \curvature = 0$. The latter is key to proving our main result, which tells us that for $\mTVector_{\mPoint}\approx \log_\mPoint \Tensor$ -- and under an additional assumption of local symmetry --, the dominant behavior of the global approximation error can be understood through the curvature tensor at the point of linearization. 

\begin{theorem}[approximation error on locally symmetric manifolds]
\label{eq:general-error-symmetric spaces}
Let $(\manifold, (\cdot, \cdot))$ be a $\dimInd$-dimensional symmetric Riemannian manifold and $\mPoint\in \manifold$ be a point on the manifold. Furthermore, let $\Tensor\in \manifold^{d_1 \times \cdots \times d_n}$ be manifold-valued tensor with entries satisfying $\Tensor_{\sumIndA_1, \ldots, \sumIndA_n} \in \exp_{\mPoint}(\mathcal{G}'_{\mPoint})$. Finally, let $\{\Theta^{\sumIndB_1, \ldots, \sumIndB_n, j}_{\mPoint}\}_{\sumIndB_1, \ldots, \sumIndB_n, \sumIndB=1}^{d_1, \cdots, d_n,d}\subset \tangent_{\mPoint} \manifold^{d_1 \times \cdots \times d_n}$ be an orthonormal frame that diagonalizes the operator 
\begin{equation}
    \Theta_{\mPoint} \mapsto \curvature_{\mPoint}(\Theta_\mPoint, \log_{\mPoint}\Tensor) \log_{\mPoint}\Tensor
\end{equation}
with respective eigenvalues $\kappa_{\sumIndB_1, \ldots, \sumIndB_n, \sumIndB}$ and define $\beta:\Real\to\Real$ as
\begin{equation}
    \beta(\kappa) := \left\{\begin{matrix}
\frac{\sinh(\sqrt{-\kappa})}{\sqrt{-\kappa}}, & \kappa <0, \\
1, & \kappa = 0, \\
\frac{\sin(\sqrt{\kappa})}{\sqrt{\kappa}}, & \kappa >0.
\label{eq:beta}
\end{matrix}\right.
\end{equation}

Then, \cref{eq:jacobi-bound-low-rank-mfld-tensor} reduces to
\begin{equation}
    d_{\manifold^{d_1 \times \cdots \times d_{n}}}(\Tensor, \exp_{\mPoint}(\mTVector_{\mPoint}))^2 = \sum_{\sumIndB_1, \ldots, \sumIndB_n, \sumIndB=1}^{\dimInd_1, \ldots, \dimInd_n, \dimInd} \beta(\kappa_{\sumIndB_1, \ldots, \sumIndB_n, \sumIndB})^2 \Bigl( \mTVector_{\mPoint} - \log_{\mPoint}\Tensor,  \Theta^{\sumIndB_1, \ldots, \sumIndB_n, j}_{\mPoint} \Bigr)_{\mPoint}^2 + \mathcal{O}(\epsilon^3)
    \label{eq:jacobi-bound-low-rank-mfld-tensor-symmetric-space}
\end{equation}
where $\epsilon:= \|\mTVector_{\mPoint} -\log_{\mPoint}\Tensor\|_{\mPoint}$.
\end{theorem}

\begin{proof}
    Let $\{\Theta^{\sumIndB_1, \ldots, \sumIndB_n, \sumIndB}(t)\}_{\sumIndB_1, \ldots, \sumIndB_n, \sumIndB=1}^{d_1, \cdots, d_n,d}\subset \tangent_{\gamma_{\mPoint,\Tensor}(t)} \manifold^{d_1 \times \cdots \times d_n}$ be the parallel transported orthonormal basis $\{\Theta^{\sumIndB_1, \ldots, \sumIndB_n, j}_{\mPoint}\}_{\sumIndB_1, \ldots, \sumIndB_n, \sumIndB=1}^{d_1, \cdots, d_n,d}\subset \tangent_{\mPoint} \manifold^{d_1 \times \cdots \times d_n}$, and 
    \begin{equation}
        J(t) := \sum_{\sumIndB_1, \ldots, \sumIndB_n, \sumIndB=1}^{\dimInd_1, \ldots, \dimInd_n, \dimInd} \frac{1}{\epsilon}(\mTVector_{\mPoint} -\log_{\mPoint}\Tensor , \Theta^{\sumIndB_1, \ldots, \sumIndB_n, j}_{\mPoint})_{\mPoint} a_{\sumIndB_1, \ldots, \sumIndB_n, \sumIndB}(t) \Theta^{\sumIndB_1, \ldots, \sumIndB_n, \sumIndB}(t)
        \label{thm:jacobi-solution-template}
    \end{equation}
    be a Jacobi field along the length minimizing geodesic $\gamma_{\mPoint,\Tensor}:[0,1]\to \manifold^{d_1 \times \cdots \times d_n}$ with initial condition $J(0)= 0 \in \tangent_{\mPoint}\manifold^{d_1 \times \cdots \times d_n}$ and initial velocity $\frac{\mathrm{D}}{\mathrm{d}t}J(t)|_{t=0} = \frac{\mTVector_{\mPoint} - \log_{\mPoint}\Tensor}{\|\mTVector_{\mPoint} - \log_{\mPoint}\Tensor\|_{\mPoint}}\in \tangent_{\mPoint}\manifold^{d_1 \times \cdots \times d_n}$.

    Substituting \cref{thm:jacobi-solution-template} into the Jacobi equation \cref{eq:jacobi-equation} decomposes the equation into ODEs in the coefficients $ a_{\sumIndB_1, \ldots, \sumIndB_n, \sumIndB}(t)$ by \cite[Prop. 3.5]{bacak2016second}. That is,
    \begin{equation}
        a_{\sumIndB_1, \ldots, \sumIndB_n, \sumIndB}''(t) + \kappa_{\sumIndB_1, \ldots, \sumIndB_n, \sumIndB} a_{\sumIndB_1, \ldots, \sumIndB_n, \sumIndB}(t) = 0, \quad \text{for $\sumIndB_1 = 1,\ldots, \dimInd_1$, \ldots, $\sumIndB_n = 1,\ldots, \dimInd_n$ and $\sumIndB = 1,\ldots, \dimInd$},
        \label{eq:jacob-decom-odes}
    \end{equation}
    with initial conditions $a_{\sumIndB_1, \ldots, \sumIndB_n, \sumIndB}(0) = 0$ and $a_{\sumIndB_1, \ldots, \sumIndB_n, \sumIndB}'(0) = 1$. Hence, solving the ODE \cref{eq:jacob-decom-odes} with the given initial conditions yields $a_{\sumIndB_1, \ldots, \sumIndB_n, \sumIndB}(1)= \beta(\kappa_{\sumIndB_1, \ldots, \sumIndB_n, \sumIndB})$, i.e., we have
    \begin{equation}
        J(1)= \sum_{\sumIndB_1, \ldots, \sumIndB_n, \sumIndB=1}^{\dimInd_1, \ldots, \dimInd_n, \dimInd} \frac{\beta(\kappa_{\sumIndB_1, \ldots, \sumIndB_n, \sumIndB})}{\epsilon}(\mTVector_{\mPoint} -\log_{\mPoint}\Tensor , \Theta^{\sumIndB_1, \ldots, \sumIndB_n, j}_{\mPoint})_{\mPoint}  \Theta^{\sumIndB_1, \ldots, \sumIndB_n, \sumIndB}(1).
    \end{equation}

    Finally, using the orthonormality of the frame $\{\Theta^{\sumIndB_1, \ldots, \sumIndB_n, \sumIndB}(t)\}_{\sumIndB_1, \ldots, \sumIndB_n, \sumIndB=1}^{d_1, \cdots, d_n,d}$, i.e., 
    \begin{equation}
        (\Theta^{\sumIndB_1, \ldots, \sumIndB_n, \sumIndB}(1), \Theta^{\sumIndB_1', \ldots, \sumIndB_n', \sumIndB'}(1))_\Tensor = \delta_{(\sumIndB_1, \ldots, \sumIndB_n, \sumIndB), (\sumIndB_1', \ldots, \sumIndB_n', \sumIndB')},
    \end{equation}
    where 
    \begin{equation}
        \delta_{(\sumIndB_1, \ldots, \sumIndB_n, \sumIndB), (\sumIndB_1', \ldots, \sumIndB_n', \sumIndB')} := \begin{cases}
         1, & \text{ if } (\sumIndB_1, \ldots, \sumIndB_n, \sumIndB)= (\sumIndB_1', \ldots, \sumIndB_n', \sumIndB'),\\
         0, & \text{ else, }
        \end{cases}
    \end{equation}
    we find
    \begin{multline}
        \epsilon^2 \|J(1)\|_\Tensor^2 =
        \epsilon^2 \sum_{\sumIndB_1, \ldots, \sumIndB_n, \sumIndB=1}^{\dimInd_1, \ldots, \dimInd_n, \dimInd} \frac{\beta(\kappa_{\sumIndB_1, \ldots, \sumIndB_n, \sumIndB})^2 }{\epsilon^2}(\mTVector_{\mPoint} -\log_{\mPoint}\Tensor , \Theta^{\sumIndB_1, \ldots, \sumIndB_n, j}_{\mPoint})_{\mPoint}^2 (\Theta^{\sumIndB_1, \ldots, \sumIndB_n, \sumIndB}(1), \Theta^{\sumIndB_1, \ldots, \sumIndB_n, \sumIndB}(1))_\Tensor\\
        = \sum_{\sumIndB_1, \ldots, \sumIndB_n, \sumIndB=1}^{\dimInd_1, \ldots, \dimInd_n, \dimInd} \beta(\kappa_{\sumIndB_1, \ldots, \sumIndB_n, \sumIndB})^2  (\mTVector_{\mPoint} -\log_{\mPoint}\Tensor , \Theta^{\sumIndB_1, \ldots, \sumIndB_n, j}_{\mPoint})_{\mPoint}^2,
        \label{eq:jacobi-bound-low-rank-mfld-tensor-symmetric-space-final-step}
    \end{multline}
    which yields the desired result \cref{eq:jacobi-bound-low-rank-mfld-tensor-symmetric-space} through substituting \cref{eq:jacobi-bound-low-rank-mfld-tensor-symmetric-space-final-step} into \cref{eq:jacobi-bound-low-rank-mfld-tensor}.
\end{proof}

\begin{remark}
\label{rem:curvature-effects}
    From \cref{eq:jacobi-bound-low-rank-mfld-tensor-symmetric-space} in \cref{eq:general-error-symmetric spaces} we get an idea of what can go wrong when global geometry is not taken into account in a tangent space-based approach to approximation of manifold-valued tensors, i.e., when just picking $\mTVector_{\mPoint} \approx \log_{\mPoint}\Tensor$ in the sense that $\|\mTVector_{\mPoint} -\log_{\mPoint}\Tensor\|_{\mPoint}$ is small. In particular, on Riemannian manifolds with negative curvature small errors in the tangent space can be amplified and on Riemannian manifolds with positive curvature large errors can be diminished. 
    In other words, naive approximation in the tangent space might give an overly optimistic approximation error for spaces with negative curvature and an overly pessimistic error for spaces with positive curvature.
\end{remark}


\section{Curvature corrected tangent space-based approximation of manifold-valued tensors.} 
\label{sec:curvature-corrected-approximation}

The observation in \cref{rem:curvature-effects} 
begs the question whether curvature effects can be taken into account when approximating manifold-valued tensors. As mentioned in \cref{sec:related-work}, attempts at taking global geometry into account, e.g., for principal geodesic analysis \cite{chakraborty2016efficient,lazar2017scale,said2007exact,sommer2014optimization}, have yielded very slow algorithms without any convergence guarantees. So in general, directly trying to minimize the non-linear and non-convex global approximation error \cref{eq:manifold-approximation-error} over some possibly non-convex set $\mathcal{S}_{\mPoint}(\manifold^{\dimInd_1\times \ldots \times \dimInd_n})\subset \tangent_{\mPoint} \manifold^{\dimInd_1 \times \cdots \times \dimInd_{n}}$ of suitable approximators, i.e., minimizing
\begin{equation}
    \inf_{\mTVector_{\mPoint}\in \mathcal{S}_{\mPoint}(\manifold^{\dimInd_1\times \ldots \times \dimInd_n})} \Bigl\{ \distance_{\manifold^{\dimInd_1 \times \cdots \times \dimInd_{n}}}(\Tensor, \exp_{\mPoint}(\mTVector_{\mPoint}))^2 \Bigr\},
    \label{eq:global-approximation-error-minimization}
\end{equation}
is expected to suffer from similar problems. In the case of symmetric Riemannian manifolds, the leading term in the right hand side of \cref{eq:jacobi-bound-low-rank-mfld-tensor-symmetric-space} is convex and differentiable with Lipschitz gradients. Therefore, 
even though the set of suitable approximations $\mathcal{S}_{\mPoint}(\manifold^{\dimInd_1\times \ldots \times \dimInd_n})$ might still be non-convex, we find ourselves in a better position.

Hence, in this section we will investigate whether minimizing the \emph{curvature corrected approximation error}
\begin{equation}
    \inf_{\mTVector_{\mPoint}\in \mathcal{S}_{\mPoint}(\manifold^{\dimInd_1\times \ldots \times \dimInd_n})} \Bigl\{ \sum_{\sumIndB_1, \ldots, \sumIndB_n, \sumIndB=1}^{\dimInd_1, \ldots, \dimInd_n, \dimInd} \beta(\kappa_{\sumIndB_1, \ldots, \sumIndB_n, \sumIndB})^2 \Bigl( \mTVector_{\mPoint} - \log_{\mPoint}\Tensor,  \Theta^{\sumIndB_1, \ldots, \sumIndB_n, \sumIndB}_{\mPoint} \Bigr)_{\mPoint}^2 \Bigr\},
    \label{eq:inf-S-linearized-tangent-approximation}
\end{equation}
yields a more suitable framework for geometry-aware approximation of manifold-valued data. In particular, in \cref{sec:diminishing-discrepancy} we will argue that under mild conditions, minimizing \cref{eq:inf-S-linearized-tangent-approximation} is a good proxy for the original problem \cref{eq:global-approximation-error-minimization}. Subsequently, in \cref{sec:stability-curvature} we will argue that under mild conditions the solution to the naive problem 
\begin{equation}
    \inf_{\mTVector_{\mPoint}\in \mathcal{S}_{\mPoint}(\manifold^{\dimInd_1\times \ldots \times \dimInd_n})} \Bigl\{ \| \mTVector_{\mPoint} - \log_{\mPoint}\Tensor \|_{\mPoint}^2 \Bigr\},
    \label{eq:inf-S-naive-tangent-approximation}
\end{equation}
will be close to the curvature corrected solution of \cref{eq:inf-S-linearized-tangent-approximation}. In particular, we will see that a minimizer of \cref{eq:inf-S-linearized-tangent-approximation} converges to the naive solution under vanishing curvature. In other words, the naive solution can be useful as a starting point for a (possibly iterative) curvature correction scheme based on a nicer minimization problem \cref{eq:inf-S-linearized-tangent-approximation} than the one based on minimizing the global approximation error \cref{eq:global-approximation-error-minimization}.

In the following it will be useful to define the the family of functions $F(\; \cdot\; ;\kappa): \tangent_{\mPoint} \manifold^{d_1 \times \cdots \times d_n} \to \Real$ parameterized by $\kappa \in \Real^{\dimInd_1\times \ldots \times \dimInd_n \times \dimInd}$ and given by
\begin{equation}
    F(\mTVector_{\mPoint}; \kappa):= \sum_{\sumIndB_1, \ldots, \sumIndB_n, \sumIndB=1}^{\dimInd_1, \ldots, \dimInd_n, \dimInd} \beta(\kappa_{\sumIndB_1, \ldots, \sumIndB_n, \sumIndB})^2 \Bigl( \mTVector_{\mPoint} - \log_{\mPoint}\Tensor,  \Theta^{\sumIndB_1, \ldots, \sumIndB_n, \sumIndB}_{\mPoint} \Bigr)_{\mPoint}^2,
    \label{eq:ccl-short-hand}
\end{equation}
where $\{\Theta^{\sumIndB_1, \ldots, \sumIndB_n, j}_{\mPoint}\}_{\sumIndB_1, \ldots, \sumIndB_n, \sumIndB=1}^{d_1, \cdots, d_n,d}\subset \tangent_{\mPoint} \manifold^{d_1 \times \cdots \times d_n}$ is a fixed orthonormal basis. Note that 
\begin{equation}
    F(\mTVector_{\mPoint}; 0)=\| \mTVector_{\mPoint} - \log_{\mPoint}\Tensor\|_{\mPoint}^2.
    \label{eq:F-for-kappa0}
\end{equation}


\subsection{Diminishing discrepancy.} 
\label{sec:diminishing-discrepancy}
As the discrepancy between the global approximation error and the curvature corrected approximation error scales with the tangent space error as $\|\mTVector_{\mPoint} - \log_{\mPoint}\Tensor\|_{\mPoint}^3$, the discrepancy will diminish if the curvature corrected approximation error is an upper bound to the tangent space error. We can formalize this in the following lemma.

\begin{lemma}
\label{lem:lower-bound}
Let $(\manifold, (\cdot, \cdot))$ be a $\dimInd$-dimensional Riemannian manifold, let $\mPoint\in \manifold$ be a point on the manifold. Furthermore, consider any real-valued tensor $\kappa \in \Real^{d_1 \times \cdots \times d_n\times \dimInd}$ such that $\kappa_{\sumIndB_1, \ldots, \sumIndB_n, \sumIndB} \leq \pi^2$ and define the function $F(\; \cdot \; ; \kappa):\tangent_{\mPoint} \manifold^{d_1 \times \cdots \times d_n} \to \Real$ as in \cref{eq:ccl-short-hand}.

Then, for any $\mTVector_{\mPoint} \in \tangent_{\mPoint} \manifold^{d_1 \times \cdots \times d_n}$
\begin{equation}
     F(\mTVector_{\mPoint}; \kappa)  \geq \beta(\kappa_{\max})^2  \; F(\mTVector_{\mPoint}; 0),
     \label{eq:lem-lower-bound-main}
 \end{equation}
 where $\kappa_{\max} : = \max_{\sumIndB_1, \ldots, \sumIndB_n, \sumIndB} \kappa_{\sumIndB_1, \ldots, \sumIndB_n, \sumIndB}$.
\end{lemma}

\begin{proof}
The result follows directly from the fact that the function $\beta$ from \cref{eq:beta} is decreasing on $(-\infty, \pi^2]$. Indeed,
\begin{multline}
    F(\mTVector_{\mPoint}; \kappa) = \sum_{\sumIndB_1, \ldots, \sumIndB_n, \sumIndB=1}^{\dimInd_1, \ldots, \dimInd_n, \dimInd} \beta(\kappa_{\sumIndB_1, \ldots, \sumIndB_n, \sumIndB})^2 \Bigl( \mTVector_{\mPoint} - \log_{\mPoint}\Tensor,  \Theta^{\sumIndB_1, \ldots, \sumIndB_n, \sumIndB}_{\mPoint} \Bigr)_{\mPoint}^2 \\
    \geq \sum_{\sumIndB_1, \ldots, \sumIndB_n, \sumIndB=1}^{\dimInd_1, \ldots, \dimInd_n, \dimInd} \beta(\kappa_{\max})^2 \Bigl( \mTVector_{\mPoint} - \log_{\mPoint}\Tensor,  \Theta^{\sumIndB_1, \ldots, \sumIndB_n, \sumIndB}_{\mPoint} \Bigr)_{\mPoint}^2 = \beta(\kappa_{\max})^2 \sum_{\sumIndB_1, \ldots, \sumIndB_n, \sumIndB=1}^{\dimInd_1, \ldots, \dimInd_n, \dimInd} \Bigl( \mTVector_{\mPoint} - \log_{\mPoint}\Tensor,  \Theta^{\sumIndB_1, \ldots, \sumIndB_n, \sumIndB}_{\mPoint} \Bigr)_{\mPoint}^2\\
    \overset{\cref{eq:F-for-kappa0}}{=} \beta(\kappa_{\max})^2  \; F(\mTVector_{\mPoint}; 0).
\end{multline}
\end{proof}

From the lemma above we find that if $\kappa_{\sumIndB_1, \ldots, \sumIndB_n, \sumIndB} < \pi^2$, minimizing the curvature corrected approximation error will automatically diminish the discrepancy with the global approximation error, because $\beta(\kappa_{\sumIndB_1, \ldots, \sumIndB_n, \sumIndB})^2>0$. Then, in the case of non-positive curvature, this condition is automatically satisfied, as $\kappa_{\sumIndB_1, \ldots, \sumIndB_n, \sumIndB} \leq 0$ and subsequently $\beta(\kappa_{\max}) \geq 1$. In the case of positive curvature, the restriction is in general not a problem either. To see this, first note that without loss of generality we may assume that 
\begin{equation}
    (\Theta^{\sumIndB_1, \ldots, \sumIndB_n, \sumIndB}_{\mPoint})_{\sumIndA_1, \ldots, \sumIndA_n} = 0 \in \tangent_\mPoint \manifold, \quad \text{unless } \sumIndA_1 = \sumIndB_1, \ldots, \sumIndA_n = \sumIndB_n,
\end{equation}
due to the assumed product manifold structure, and additionally we may assume that
\begin{equation}
    (\Theta^{\sumIndB_1, \ldots, \sumIndB_n, 1}_{\mPoint})_{\sumIndA_1, \ldots, \sumIndA_n} = \left\{\begin{matrix}
    \frac{(\log_{\mPoint}\Tensor)_{\sumIndB_1, \ldots, \sumIndB_n}}{\|(\log_{\mPoint}\Tensor)_{\sumIndB_1, \ldots, \sumIndB_n}\|_{\mPoint}} & \text{if } \sumIndA_1 = \sumIndB_1, \ldots, \sumIndA_n = \sumIndB_n,  \\
    0 &  \text{otherwise},\\
    \end{matrix}\right.
    \label{eq:first-eigenvector-curvature-operator-1}
\end{equation}
because the right-hand-side of \cref{eq:first-eigenvector-curvature-operator-1} is an eigenvector of $R(\cdot, \log_{\mPoint}\Tensor) \log_{\mPoint}\Tensor$ with eigenvalue $\kappa_{\sumIndB_1, \ldots, \sumIndB_n, 1} = 0$. In other words, for $\sumIndB = 1$ we find
\begin{equation}
    \kappa_{\sumIndB_1, \ldots, \sumIndB_n, 1} = 0 \leq \pi^2,
\end{equation}

For the remaining cases $j \geq 2$, the above discussion and some rewriting tells us that we have the equality
\begin{equation}
    \kappa_{\sumIndB_1, \ldots, \sumIndB_n, \sumIndB} = \sectionalcurvature_\mPoint((\Theta^{\sumIndB_1, \ldots, \sumIndB_n, \sumIndB}_{\mPoint})_{\sumIndB_1, \ldots, \sumIndB_n}, \log_{\mPoint}\Tensor_{\sumIndB_1, \ldots, \sumIndB_n}) \distance_\manifold(\mPoint, \Tensor_{\sumIndB_1, \ldots, \sumIndB_n})^2,
    \label{eq:kappa-sectional-curvature-distance}
\end{equation}
where $\sectionalcurvature_\mPoint:\tangent_\mPoint \manifold \times \tangent_\mPoint \manifold \to \Real$ is the sectional curvature. In the special case of constant sectional curvature $\bar{\kappa}$ we then find
\begin{multline}
    \kappa_{\sumIndB_1, \ldots, \sumIndB_n, \sumIndB} = \sectionalcurvature_\mPoint((\Theta^{\sumIndB_1, \ldots, \sumIndB_n, \sumIndB}_{\mPoint})_{\sumIndB_1, \ldots, \sumIndB_n}, \log_{\mPoint}\Tensor_{\sumIndB_1, \ldots, \sumIndB_n}) \distance_\manifold(\mPoint, \Tensor_{\sumIndB_1, \ldots, \sumIndB_n})^2 \\
    \overset{\text{constant curvature}}{=} \bar{\kappa} \distance_\manifold(\mPoint, \Tensor_{\sumIndB_1, \ldots, \sumIndB_n})^2 \overset{\text{Bonnet–Myers}}{\leq} \bar{\kappa} \frac{\pi^2}{\bar{\kappa}} = \pi^2.
\end{multline}
In other words, in the case of a sphere we only get $\kappa_{\sumIndB_1, \ldots, \sumIndB_n, \sumIndB} = \pi^2$ if one or more of the entries of the manifold-valued tensor $\Tensor$ lies on the cut locus of $\mPoint$, which is on the opposite pole of the sphere. In practice, for the case of the sphere one has to be careful with picking the point $\mPoint$ given the data. For general symmetric Riemannian manifolds with positive curvature we expect that similar care will be sufficient for avoiding $\kappa_{\sumIndB_1, \ldots, \sumIndB_n, \sumIndB} = \pi^2$.

We conclude this part with noting that \cref{lem:lower-bound} gives us a practical lower bound for the approximation error.
\begin{corollary}[Approximation error lower bounds]
\label{cor:approximation-error-lower-bound}
Under the assumptions of \cref{lem:lower-bound} and for a subset  $\mathcal{S}_{\mPoint}(\manifold^{\dimInd_1\times \ldots \times \dimInd_n})\subset \tangent_{\mPoint} \manifold^{\dimInd_1 \times \cdots \times \dimInd_{n}}$
\begin{equation}
    \hspace{-0.8cm}
     \inf_{\mTVector_{\mPoint}\in \mathcal{S}_{\mPoint}(\manifold^{\dimInd_1\times \ldots \times \dimInd_n})} \Bigl\{ \sum_{\sumIndB_1, \ldots, \sumIndB_n, \sumIndB=1}^{\dimInd_1, \ldots, \dimInd_n, \dimInd} \beta(\kappa_{\sumIndB_1, \ldots, \sumIndB_n, \sumIndB})^2 \Bigl( \mTVector_{\mPoint} - \log_{\mPoint}\Tensor,  \Theta^{\sumIndB_1, \ldots, \sumIndB_n, \sumIndB}_{\mPoint} \Bigr)_{\mPoint}^2 \Bigr\} \geq \inf_{\mTVector_{\mPoint}\in \mathcal{S}_{\mPoint}(\manifold^{\dimInd_1\times \ldots \times \dimInd_n})} \Bigl\{ \beta(\kappa_{\max})^2 \| \mTVector_{\mPoint} - \log_{\mPoint}\Tensor\|_{\mPoint}^2 \Bigr\} .
     \label{eq:cor-lower-bound}
\end{equation}
\end{corollary}



\subsection{Stability with respect to curvature effects.} 
\label{sec:stability-curvature}
Because sectional curvature is by definition scale invariant, the identity \cref{eq:kappa-sectional-curvature-distance} tells us that not only the intrinsic curvature of the Riemannian manifold will be a determining factor for the severity of the curvature effects, but even so will the distance of the entries of the manifold-valued tensor $\Tensor \in \manifold^{\dimInd_1\times \ldots \times \dimInd_n}$ to the point $\mPoint\in \manifold$.

The above observation motivates one to consider the relation between minimizers of the curvature corrected approximation problem \cref{eq:inf-S-linearized-tangent-approximation} and the naive problem \cref{eq:inf-S-naive-tangent-approximation}. A natural question to ask is whether the minimizers of \cref{eq:inf-S-naive-tangent-approximation} are stable under curvature effects. That is, will minimizers of \cref{eq:inf-S-linearized-tangent-approximation} converge to minimizers of \cref{eq:inf-S-naive-tangent-approximation} under vanishing $\kappa\in \Real^{\dimInd_1\times \ldots \times \dimInd_n}$? Note that this is a non-trivial question since the admissible set $\mathcal{S}_{\mPoint}(\manifold^{\dimInd_1\times \ldots \times \dimInd_n})$ is non-convex and we have not assumed any additional structure yet.

The main tool for answering this question will be the notion of $\Gamma$-convergence and equi-coercivity from calculus of variations.

\begin{definition}[$\Gamma$-convergence]
    Let $(F_m)_{m=1}^\infty$ be a sequence of functions $F_m: X \rightarrow \Real$ on some topological space $X$. Then we say that $F_m$ $\Gamma$-converges to a limit $F$ and write $F_m \stackrel{\Gamma}{\rightarrow} F$ if the following conditions hold:
    \begin{itemize}
        \item For any $x \in X$ and any sequence $(x_m)_{m=1}^\infty$ such that $x_m \rightarrow x$, we have
        \begin{equation}
            F(x) \leq \liminf _{m \rightarrow \infty} F_m (x_m).
        \end{equation}
        \item For any $x \in X$, we can find a sequence $x_i \rightarrow x$ such that
        \begin{equation}
            F(x) \geq \limsup _{m \rightarrow \infty} F_m (x_m).
        \end{equation}
    \end{itemize}
\end{definition}

\begin{definition}[equi-coercivity]
    A family of functions $F_m: X \rightarrow \Real$ for $m\in \Natural$ is equi-coercive if for all $\alpha\in \Real$, there exists a compact set $K_\alpha \subset X$ such that $\left\{x \in X: F_m(x) \leq \alpha\right\} \subset K_\alpha$ for all $m$.
\end{definition}

Combining $\Gamma$-convergence and equi-coercivity guarantees convergence of minimizers through the following well-known result.
\begin{theorem}[{\cite[Fundamental Theorem of $\Gamma$-Convergence]{braides2006handbook}}]
\label{thm:fundamental-theorem-gamma-convergence}
    If $F_m \stackrel{\Gamma}{\rightarrow} F$ and the family $(F_m)_{m=1}^\infty$ is equi-coercive, then the sequence of minimizers $x_m \in \arg \min _{x \in X} F_m(x)$ converges to some $x \in \arg \min _{x \in X} F(x)$.
\end{theorem}

Before we are ready to answer the question above, we will prove the following lemma.

\begin{lemma}
\label{lem:upper-bound}
Under the assumptions of \cref{lem:lower-bound}
    \begin{equation}
     |F(\mTVector_{\mPoint}; \kappa) - F(\mTVector_{\mPoint}; 0)| \leq  \max\{( \beta(\kappa_{\min})^2 - 1 ), (1 - \beta(\kappa_{\max})^2) \} \|\mTVector_{\mPoint}\|_{\mPoint}^2,
     \label{eq:lem-upper-bound-main}
 \end{equation}
 where $\kappa_{\max} : = \max_{\sumIndB_1, \ldots, \sumIndB_n, \sumIndB} \kappa_{\sumIndB_1, \ldots, \sumIndB_n, \sumIndB}$ and $\kappa_{\min} : = \min_{\sumIndB_1, \ldots, \sumIndB_n, \sumIndB} \kappa_{\sumIndB_1, \ldots, \sumIndB_n, \sumIndB}$.
\end{lemma}

\begin{proof}
We will first show that the equality
    \begin{equation}
        F(\mTVector_{\mPoint}; \kappa)
        = \sum_{\sumIndB_1, \ldots, \sumIndB_n, j=1}^{\dimInd_1, \ldots, \dimInd_n, \dimInd} \biggl(\beta(\kappa_{\sumIndB_1, \ldots, \sumIndB_n, \sumIndB})^2 \Bigl( \mTVector_{\mPoint},  \Theta^{\sumIndB_1, \ldots, \sumIndB_n, j}_{\mPoint} \Bigr)_{\mPoint}^2 \biggr) - 2\Bigl( \mTVector_{\mPoint} , \log_{\mPoint}\Tensor \Bigr)_{\mPoint} + \| \log_{\mPoint}\Tensor\|_{\mPoint}^2
        \label{eq:lem-upper-bound-betas-rewrite-ccloss}
    \end{equation}
    holds. For that, first note (once more) that without loss of generality we may assume that 
    \begin{equation}
        (\Theta^{\sumIndB_1, \ldots, \sumIndB_n, 1}_{\mPoint})_{\sumIndB_1', \ldots, \sumIndB_n'} = \left\{\begin{matrix}
        \frac{(\log_{\mPoint}\Tensor)_{\sumIndB_1, \ldots, \sumIndB_n}}{\|(\log_{\mPoint}\Tensor)_{\sumIndB_1, \ldots, \sumIndB_n}\|_{\mPoint}} & \text{if } \sumIndB_1' = \sumIndB_1, \ldots, \sumIndB_n' = \sumIndB_n,  \\
        0 &  \text{otherwise}.\\
        \end{matrix}\right.
        \label{eq:first-eigenvector-curvature-operator}
    \end{equation}
    Indeed, the right-hand-side of \cref{eq:first-eigenvector-curvature-operator} is an eigenvector of $R(\cdot, \log_{\mPoint}\Tensor) \log_{\mPoint}\Tensor$ with eigenvalue $\kappa_{\sumIndB_1, \ldots, \sumIndB_n, 1} = 0$. 
    
    Then, this also means that $\log_{\mPoint}\Tensor \perp \Theta^{\sumIndB_1, \ldots, \sumIndB_n, \sumIndB}_{\mPoint}$ for all $\sumIndB_1 = 1,\ldots, \dimInd_1$ to $\sumIndB_n = 1, \ldots, \dimInd_n$ and for $\sumIndB = 2, \ldots, \dimInd$ and that 
    \begin{equation}
        \sum_{\sumIndB_1, \ldots, \sumIndB_n=1}^{\dimInd_1, \ldots, \dimInd_n} \|(\log_{\mPoint}\Tensor)_{\sumIndB_1, \ldots, \sumIndB_n}\|_{\mPoint}  \Theta^{\sumIndB_1, \ldots, \sumIndB_n, 1}_{\mPoint} = \log_{\mPoint}\Tensor
        \label{eq:sum-identity-basis-vectors-curvature-operator}
    \end{equation}
    Furthermore, for $\Phi_{\mPoint} \in \tangent_{\mPoint} \manifold^{\dimInd_1 \times \cdots \times \dimInd_n}$
    \begin{multline}
        \sum_{\sumIndB_1, \ldots, \sumIndB_n, j=1}^{\dimInd_1, \ldots, \dimInd_n, \dimInd} \beta(\kappa_{\sumIndB_1, \ldots, \sumIndB_n, \sumIndB})^2 \Bigl( \Phi_{\mPoint} ,  \Theta^{\sumIndB_1, \ldots, \sumIndB_n, \sumIndB}_{\mPoint} \Bigr)_{\mPoint}\Bigl(\log_{\mPoint}\Tensor,  \Theta^{\sumIndB_1, \ldots, \sumIndB_n, \sumIndB}_{\mPoint} \Bigr)_{\mPoint} \\
        \overset{\text{orthogonality}}{=} \sum_{\sumIndB_1, \ldots, \sumIndB_n=1}^{\dimInd_1, \ldots, \dimInd_n} \beta(0)^2 \Bigl( \Phi_{\mPoint} ,  \Theta^{\sumIndB_1, \ldots, \sumIndB_n, 1}_{\mPoint} \Bigr)_{\mPoint}\Bigl(\log_{\mPoint}\Tensor,  \Theta^{\sumIndB_1, \ldots, \sumIndB_n, 1}_{\mPoint} \Bigr)_{\mPoint} \\ 
        \overset{\cref{eq:first-eigenvector-curvature-operator}}{=} \sum_{\sumIndB_1, \ldots, \sumIndB_n=1}^{\dimInd_1, \ldots, \dimInd_n} \Bigl( \Phi_{\mPoint} ,  \Theta^{\sumIndB_1, \ldots, \sumIndB_n, 1}_{\mPoint} \Bigr)_{\mPoint} \|(\log_{\mPoint}\Tensor)_{\sumIndB_1, \ldots, \sumIndB_n}\|_{\mPoint} \\
        = \Bigl( \Phi_{\mPoint} , \sum_{\sumIndB_1, \ldots, \sumIndB_n=1}^{\dimInd_1, \ldots, \dimInd_n} \|(\log_{\mPoint}\Tensor)_{\sumIndB_1, \ldots, \sumIndB_n}\|_{\mPoint}  \Theta^{\sumIndB_1, \ldots, \sumIndB_n, 1}_{\mPoint} \Bigr)_{\mPoint} 
        \overset{\cref{eq:sum-identity-basis-vectors-curvature-operator}}{=} \Bigl( \Phi_{\mPoint} , \log_{\mPoint}\Tensor \Bigr)_{\mPoint}.
        \label{eq:phi-q-inner-sum-rewrite-orthogonality}
    \end{multline}

    Finally, we can show \cref{eq:lem-upper-bound-betas-rewrite-ccloss}. We have
    \begin{multline}
        F(\mTVector_{\mPoint}; \kappa) = \sum_{\sumIndB_1, \ldots, \sumIndB_n, j=1}^{\dimInd_1, \ldots, \dimInd_n, \dimInd} \beta(\kappa_{\sumIndB_1, \ldots, \sumIndB_n, j})^2 \Bigl( \mTVector_{\mPoint} - \log_{\mPoint}\Tensor,  \Theta^{\sumIndB_1, \ldots, \sumIndB_n, j}_{\mPoint} \Bigr)_{\mPoint}^2 \\
        = \sum_{\sumIndB_1, \ldots, \sumIndB_n, j=1}^{\dimInd_1, \ldots, \dimInd_n, \dimInd} \beta(\kappa_{\sumIndB_1, \ldots, \sumIndB_n, j})^2 \Bigl[ \Bigl( \mTVector_{\mPoint},  \Theta^{\sumIndB_1, \ldots, \sumIndB_n, j}_{\mPoint} \Bigr)_{\mPoint}^2 - 2\Bigl( \mTVector_{\mPoint} ,  \Theta^{\sumIndB_1, \ldots, \sumIndB_n, j}_{\mPoint} \Bigr)_{\mPoint}\Bigl(\log_{\mPoint}\Tensor,  \Theta^{\sumIndB_1, \ldots, \sumIndB_n, j}_{\mPoint} \Bigr)_{\mPoint} + \Bigl( \log_{\mPoint}\Tensor,  \Theta^{\sumIndB_1, \ldots, \sumIndB_n, j}_{\mPoint} \Bigr)_{\mPoint}^2 \Bigr]\\
        \overset{\cref{eq:phi-q-inner-sum-rewrite-orthogonality}}{=} \sum_{\sumIndB_1, \ldots, \sumIndB_n, j=1}^{\dimInd_1, \ldots, \dimInd_n, \dimInd} \biggl(\beta(\kappa_{\sumIndB_1, \ldots, \sumIndB_n, \sumIndB})^2 \Bigl( \mTVector_{\mPoint},  \Theta^{\sumIndB_1, \ldots, \sumIndB_n, j}_{\mPoint} \Bigr)_{\mPoint}^2 \biggr) - 2( \mTVector_{\mPoint} , \log_{\mPoint}\Tensor )_{\mPoint} + \| \log_{\mPoint}\Tensor\|_{\mPoint}^2,
        \label{eq:lem-upper-bound-rewrite-first-term}
    \end{multline}
    which gives us \cref{eq:lem-upper-bound-betas-rewrite-ccloss}.
    
    Next, as \cref{eq:first-eigenvector-curvature-operator} is an eigenvector with eigenvalue $0$, we have that $\kappa_{\min}\leq 0 \leq \kappa_{\max}$. It follows that 
    \begin{equation}
        \beta(\kappa_{\min})^2 - 1 \geq 0 \quad \text{and} \quad 1 - \beta(\kappa_{\max})^2 \geq 0.
    \end{equation}
     So for proving the claim \cref{eq:lem-upper-bound-main}, it is sufficient to show
    \begin{equation}
         F(\mTVector_{\mPoint}; \kappa) - F(\mTVector_{\mPoint}; 0) \leq (\beta(\kappa_{\min})^2 - 1) \| \mTVector_{\mPoint}\|_{\mPoint}^2
         \label{eq:lem-upper-bound-upper}
     \end{equation}
     and 
     \begin{equation}
         F(\mTVector_{\mPoint}; \kappa) - F(\mTVector_{\mPoint}; 0) \geq - (1 - \beta(\kappa_{\max})^2) \| \mTVector_{\mPoint}\|_{\mPoint}^2.
         \label{eq:lem-upper-bound-lower}
     \end{equation}
     Rewriting $\| \mTVector_{\mPoint} - \log_{\mPoint}\Tensor\|_{\mPoint}^2 = \| \mTVector_{\mPoint}\|_{\mPoint}^2 - 2 ( \mTVector_{\mPoint} , \log_{\mPoint}\Tensor )_{\mPoint} + \|\log_{\mPoint}\Tensor\|_{\mPoint}^2$ we find
     \begin{multline}
         F(\mTVector_{\mPoint}; \kappa) - F(\mTVector_{\mPoint}; 0) = \sum_{\sumIndB_1, \ldots, \sumIndB_n, \sumIndB=1}^{\dimInd_1, \ldots, \dimInd_n, \dimInd} \beta(\kappa_{\sumIndB_1, \ldots, \sumIndB_n, \sumIndB})^2 \Bigl( \mTVector_{\mPoint} - \log_{\mPoint}\Tensor,  \Theta^{\sumIndB_1, \ldots, \sumIndB_n, \sumIndB}_{\mPoint} \Bigr)_{\mPoint}^2  - \| \mTVector_{\mPoint} - \log_{\mPoint}\Tensor\|_{\mPoint}^2 \\
         \overset{\cref{eq:lem-upper-bound-rewrite-first-term}}{=} \sum_{\sumIndB_1, \ldots, \sumIndB_n, j=1}^{\dimInd_1, \ldots, \dimInd_n, \dimInd} \biggl(\beta(\kappa_{\sumIndB_1, \ldots, \sumIndB_n, \sumIndB})^2 \Bigl( \mTVector_{\mPoint},  \Theta^{\sumIndB_1, \ldots, \sumIndB_n, j}_{\mPoint} \Bigr)_{\mPoint}^2 \biggr) - \| \mTVector_{\mPoint}\|_{\mPoint}^2\\
         \overset{\text{$\beta$ decreasing}}{\leq} \sum_{\sumIndB_1, \ldots, \sumIndB_n, j=1}^{\dimInd_1, \ldots, \dimInd_n, \dimInd} \biggl(\beta(\kappa_{\min})^2 \Bigl( \mTVector_{\mPoint},  \Theta^{\sumIndB_1, \ldots, \sumIndB_n, j}_{\mPoint} \Bigr)_{\mPoint}^2 \biggr) - \| \mTVector_{\mPoint}\|_{\mPoint}^2\\
         = (\beta(\kappa_{\min})^2 - 1) \| \mTVector_{\mPoint}\|_{\mPoint}^2
     \end{multline}
     and
     \begin{multline}
         F(\mTVector_{\mPoint}; \kappa) - F(\mTVector_{\mPoint}; 0) = \sum_{\sumIndB_1, \ldots, \sumIndB_n, \sumIndB=1}^{\dimInd_1, \ldots, \dimInd_n, \dimInd} \beta(\kappa_{\sumIndB_1, \ldots, \sumIndB_n, \sumIndB})^2 \Bigl( \mTVector_{\mPoint} - \log_{\mPoint}\Tensor,  \Theta^{\sumIndB_1, \ldots, \sumIndB_n, \sumIndB}_{\mPoint} \Bigr)_{\mPoint}^2  - \| \mTVector_{\mPoint} - \log_{\mPoint}\Tensor\|_{\mPoint}^2 \\
         \overset{\cref{eq:lem-upper-bound-rewrite-first-term}}{=} \sum_{\sumIndB_1, \ldots, \sumIndB_n, j=1}^{\dimInd_1, \ldots, \dimInd_n, \dimInd} \biggl(\beta(\kappa_{\sumIndB_1, \ldots, \sumIndB_n, \sumIndB})^2 \Bigl( \mTVector_{\mPoint},  \Theta^{\sumIndB_1, \ldots, \sumIndB_n, j}_{\mPoint} \Bigr)_{\mPoint}^2 \biggr) - \| \mTVector_{\mPoint}\|_{\mPoint}^2\\
         \overset{\text{$\beta$ decreasing}}{\geq} \sum_{\sumIndB_1, \ldots, \sumIndB_n, j=1}^{\dimInd_1, \ldots, \dimInd_n, \dimInd} \biggl(\beta(\kappa_{\max})^2 \Bigl( \mTVector_{\mPoint},  \Theta^{\sumIndB_1, \ldots, \sumIndB_n, j}_{\mPoint} \Bigr)_{\mPoint}^2 \biggr) - \| \mTVector_{\mPoint}\|_{\mPoint}^2\\
         = - (1 - \beta(\kappa_{\max})^2) \| \mTVector_{\mPoint}\|_{\mPoint}^2,
     \end{multline}
     which gives us \cref{eq:lem-upper-bound-upper,eq:lem-upper-bound-lower} and proves the claim.
\end{proof}

We are now ready to state and prove the following result.

\begin{theorem}[Convergence of minimizers for vanishing curvature]
\label{thm:convergence-of-minimizers}
    Let $(\manifold, (\cdot, \cdot))$ be a Riemannian manifold, let $\mPoint\in \manifold$ be a point on the manifold and consider any subset $\mathcal{S}_{\mPoint}(\manifold^{\dimInd_1\times \ldots \times \dimInd_n})\subset \tangent_{\mPoint} \manifold^{d_1 \times \cdots \times d_n}$. Furthermore, consider any manifold-valued tensor $\Tensor\in \manifold^{\dimInd_1 \times \cdots \times \dimInd_n}$, any sequence of real-valued tensors $(\kappa^m)_{m=1}^\infty \subset \Real^{\dimInd_1 \times \cdots \times \dimInd_n \times \dimInd}$ such that each $\kappa^m_{\sumIndB_1, \ldots, \sumIndB_n, \sumIndB} \to 0$ as $m\to \infty$ and the corresponding family of functions $\{ F(\; \cdot \; ; \kappa^m)\}_{m=1}^\infty$ as in \cref{eq:ccl-short-hand}.

    If the set $\mathcal{S}_{\mPoint}(\manifold^{\dimInd_1\times \ldots \times \dimInd_n})\subset \tangent_{\mPoint} \manifold^{d_1 \times \cdots \times d_n}$ is closed under the product metric, i.e., under $\|\cdot\|_\mPoint$ on $\tangent_{\mPoint} \manifold^{d_1 \times \cdots \times d_n}$, then the sequence of minimizers $(\Xi^{*,m}_{\mPoint})_{m=1}^\infty \subset \mathcal{S}_{\mPoint}(\manifold^{\dimInd_1\times \ldots \times \dimInd_n})$ of $F(\mTVector_{\mPoint}; \kappa^m)$ over $ \mathcal{S}_{\mPoint}(\manifold^{\dimInd_1\times \ldots \times \dimInd_n})$ converges in the product metric on the tangent space $\tangent_{\mPoint}\manifold^{\dimInd_1\times \ldots \times \dimInd_n}$ to a minimizer of the zero-curvature case:
    \begin{equation}
        \Xi^{*,m}_{\mPoint} \to \Xi^{*}_{\mPoint} \in \argmin_{\Xi_\mPoint\in \mathcal{S}_{\mPoint}(\manifold^{\dimInd_1\times \ldots \times \dimInd_n})} \| \mTVector_{\mPoint} - \log_{\mPoint}\Tensor\|_{\mPoint}^2.
    \end{equation}
\end{theorem}

\begin{proof}
    For notational convenience we define $F_m : \mathcal{S}_{\mPoint}(\manifold^{\dimInd_1\times \ldots \times \dimInd_n}) \to \Real$ and $F: \mathcal{S}_{\mPoint}(\manifold^{\dimInd_1\times \ldots \times \dimInd_n})\to \Real$ as
    \begin{equation}
        F_m (\mTVector_{\mPoint}):= F(\mTVector_{\mPoint} ; \kappa^m) \quad \text{and} \quad F(\mTVector_{\mPoint}) := F(\mTVector_{\mPoint};0).
    \end{equation}
    By \cref{thm:fundamental-theorem-gamma-convergence} it is sufficient to show that $F_m \overset{\Gamma}{\to} F$ and that $(F_m)_{m=1}^\infty$ is equi-coercive on $\mathcal{S}_{\mPoint}(\manifold^{\dimInd_1\times \ldots \times \dimInd_n})$.

    [$\Gamma$-convergence] Before we move on to proving the claim, we point out that for any tangent vector $\mTVector_{\mPoint} \in \mathcal{S}_{\mPoint}(\manifold^{\dimInd_1\times \ldots \times \dimInd_n})$ and any sequence $(\Xi^m_{\mPoint})_{m=1}^\infty \subset \mathcal{S}_{\mPoint}(\manifold^{\dimInd_1\times \ldots \times \dimInd_n})$ with $\Xi^m_{\mPoint} \to \mTVector_{\mPoint}$  -- in the topology generated by the metric on $\tangent_{\mPoint} \manifold^{\dimInd_1 \times \ldots \times \dimInd_n}$ -- we have by assumptions on $\mathcal{S}_{\mPoint}(\manifold^{\dimInd_1\times \ldots \times \dimInd_n})$ being closed under the metric on $\tangent_{\mPoint} \manifold^{\dimInd_1 \times \ldots \times \dimInd_n}$ that 
    \begin{equation}
        \lim_{m\to\infty} \|\Xi^m_{\mPoint}\|_{\mPoint} = \|\mTVector_{\mPoint}\|_{\mPoint} \quad \text{and} \quad \lim_{m\to\infty} \| \Xi^m_{\mPoint} - \log_{\mPoint}\Tensor\|_{\mPoint}^2 = \| \mTVector_{\mPoint} - \log_{\mPoint}\Tensor\|_{\mPoint}^2.
        \label{eq:lemma-Sq-closed-metric-convergence}
    \end{equation}

     Now, choose any $\mTVector_{\mPoint}\in \mathcal{S}_{\mPoint}(\manifold^{\dimInd_1\times \ldots \times \dimInd_n})$ and any sequence $(\Xi^m_{\mPoint})_{m=1}^\infty \subset \mathcal{S}_{\mPoint}(\manifold^{\dimInd_1\times \ldots \times \dimInd_n})$ such that $\Xi^m_{\mPoint} \to \mTVector_{\mPoint}$. To verify $\Gamma$-convergence we have to establish the two defining properties, the limsup and the liminf inequality. For proving the claim,  we will prove the stronger statement $\lim_{m\to\infty} F_m (\Xi^m_{\mPoint}) = F(\mTVector_{\mPoint})$. In particular, using the triangle inequality we will show that
    \begin{equation}
        |F_m(\Xi^m_{\mPoint}) - F(\mTVector_{\mPoint})| \leq |F_m(\Xi^m_{\mPoint}) - F(\Xi^m_{\mPoint})|  + |F(\Xi^m_{\mPoint}) - F(\mTVector_{\mPoint})| \to 0.
    \end{equation}
    For showing that the the first term vanishes, we invoke \cref{lem:upper-bound} and see that
    \begin{equation}
        \lim_{m\to\infty} |F_m(\Xi^m_{\mPoint}) - F(\Xi^m_{\mPoint})| \overset{\cref{eq:lem-upper-bound-main}}{\leq} \lim_{m\to\infty}  \max\{( \beta(\kappa_{\min})^2 - 1 ), (1 - \beta(\kappa_{\max})^2) \} \|\Xi^m_{\mPoint}\|_{\mPoint}^2 \overset{\cref{eq:beta}}{=} 0,
    \end{equation}
    because $\beta$ is continuous, $\kappa^m_{\max} \to 0$, $\kappa^m_{\min} \to 0$ and because $\|\Xi^m_{\mPoint}\|_{\mPoint}^2$ is bounded as $\|\Xi^m_{\mPoint}\|_{\mPoint}^2 \to \|\mTVector_{\mPoint}\|_{\mPoint}^2 < \infty$ by the first equality in \cref{eq:lemma-Sq-closed-metric-convergence}. The second term vanishes because of the second equality in \cref{eq:lemma-Sq-closed-metric-convergence}.

    [Equi-coercivity] We need to show that $(F_m)_{m=1}^\infty$ is equi-coercive on $\mathcal{S}_{\mPoint}(\manifold^{\dimInd_1\times \ldots \times \dimInd_n})$. For that, choose $\alpha \in \Real$. Without loss of generality, we may assume that $\kappa^m_{\sumIndB_1, \ldots, \sumIndB_n, \sumIndB}\leq \pi^2 -\delta$ for some $\delta \in (0,\pi^2)$. Then, by \cref{lem:lower-bound} we find that
    \begin{equation}
        F_m (\mTVector_{\mPoint}) \overset{\cref{eq:lem-lower-bound-main}}{\geq} \beta(\kappa^m_{\max})^2  \; F(\mTVector_{\mPoint}) \overset{\text{$\beta$ decreasing}}{\geq}  \beta(\pi^2 -\delta)^2  \; F(\mTVector_{\mPoint}).
        \label{eq:th-equi-coerc-level-set}
    \end{equation}
    From \cref{eq:th-equi-coerc-level-set} we find that 
    \begin{equation}
        \{\mTVector_{\mPoint}\in\tangent_{\mPoint} \manifold^{\dimInd_1 \times \ldots \times \dimInd_n}\;\mid\;F_m (\mTVector_{\mPoint})\leq\alpha\} \subset \{\mTVector_{\mPoint}\in\tangent_{\mPoint} \manifold^{\dimInd_1 \times \ldots \times \dimInd_n}\;\mid\; \beta(\pi^2 -\delta)^2  \;  F (\mTVector_{\mPoint})\leq\alpha\}.
    \end{equation}
    Since every $F$ is strongly convex on the tangent space $\tangent_{\mPoint} \manifold^{\dimInd_1 \times \ldots \times \dimInd_n}$ -- which is isomorphic to $\Real^\dimInd$ -- and $\beta(\pi^2 -\delta)^2 >0$ the sub-level sets $\{\mTVector_{\mPoint}\in\tangent_{\mPoint} \manifold^{\dimInd_1 \times \ldots \times \dimInd_n}\;\mid\;\beta(\pi^2 -\delta)^2 F (\mTVector_{\mPoint})\leq\alpha\}$ are closed and bounded and hence compact. Since $\mathcal{S}_{\mPoint}(\manifold^{\dimInd_1\times \ldots \times \dimInd_n})$ is a closed subset of $\tangent_{\mPoint} \manifold^{\dimInd_1 \times \ldots \times \dimInd_n}$ by assumption, the restriction $K_\alpha := \{\mTVector_{\mPoint}\in\mathcal{S}_{\mPoint}(\manifold^{\dimInd_1\times \ldots \times \dimInd_n}) \;\mid\; \beta(\pi^2 -\delta)^2 F (\mTVector_{\mPoint})\leq\alpha\}$ is compact as well. Finally, we have for all $m$
    \begin{equation}
        \{ \mTVector_{\mPoint}\in\mathcal{S}_{\mPoint}(\manifold^{\dimInd_1\times \ldots \times \dimInd_n})\;\mid\;F_m (\mTVector_{\mPoint})\leq\alpha\} \subset K_\alpha, 
    \end{equation}
    which proves the claim.
\end{proof}


\section{Curvature corrected low-rank approximation.} 
\label{sec:cc-lra}

In \cref{sec:curvature-corrected-approximation} we have shown that minimizing the curvature corrected approximation error \cref{eq:inf-S-linearized-tangent-approximation} can be suitable for tangent space-based approximation of manifold-valued tensors. In this section we will use the observations from the above in a case study. In particular, we will be concerned with low manifold-valued tensor rank, which will be introduced in \cref{sec:manifold-valued-tensor-rank} along with a generalization of the \emph{higher-order singular value decomposition} (HOSVD) to the tangent space. Subsequently, in \cref{sec:curvature-correcteed-hosvd} we will consider curvature corrected low multi-linear rank approximation on general symmetric Riemannian manifolds and propose the \emph{curvature corrected truncated higher-order singular value decomposition} (CC-tHOSVD). We remind the reader again that works that have tried to construct a generalization of low rank approximation have had to choose between their method being generally applicable to arbitrary manifolds, numerically feasible or global-geometry aware, and that typically only two out of three of these features has been accounted for. In \cref{sec:curvature-correcteed-hosvd} we will discuss how CC-tHOSVD takes all three considerations into account. Finally, in \cref{sec:metric-corrected-hosvd} for the sake of completeness and for later benchmarking we will also propose the \emph{metric corrected truncated higher-order singular value decomposition} (MC-tHOSVD), which is an algorithm for minimizing the unaltered global approximation error \cref{eq:global-approximation-error-minimization} on general symmetric Riemannian manifolds.



\subsection{Rank-based decomposition of manifold-valued tensors.} 
\label{sec:manifold-valued-tensor-rank}
Unlike the well-defined and easy-to-compute notion of (real-valued) matrix rank, there are several rank-related notions for real-valued tensors. Pioneering work was done by Hitchcock \cite{hitchcock1927expression,hitchcock1928multiple} in 1927, who introduced a decomposition of a real-valued tensor as the sum of a finite number of rank-one real-valued tensors, which comes with a notion of rank. The decomposition was introduced again in 1944 by Cattell \cite{cattell1944parallel,cattell1952three}, but gained popularity after its third introduction by Carroll and Chang \cite{carroll1970analysis} and Harshman \cite{harshman1970foundations} in 1970 as CANDECOMP (canonical decomposition) and PARAFAC (parallel factors), respectively. The rank-one tensor decomposition was discovered independently once more by M\"ocks \cite{mocks1988topographic}, but since Kiers \cite{kiers2000towards} the decomposition is referred to as the CANDECOMP/PARAFAC decomposition or CP decomposition \cite{kolda2009tensor}. Subsequently, the corresponding notion of real-valued tensor rank is referred to as the CP rank of a tensor.

In the case of manifold-valued tensors, we can generalize the CP rank. That is, for a manifold-valued tensor $\Tensor\in \manifold^{\dimInd_1\times \ldots \times \dimInd_n}$ we define the \emph{CP rank at a point $\mPoint\in \manifold$} as the smallest integer $r$ for which we can write
\begin{equation}
    \log_{\mPoint} \Tensor = \sum_{\sumIndA=1}^r  (\mathbf{a}^{1,\sumIndA} \otimes \ldots \otimes \mathbf{a}^{n,\sumIndA})\eta_{\mPoint}^\sumIndA := \Bigl(\sum_{\sumIndA = 1}^r \bigl[(\mathbf{a}^{1,\sumIndA})_{\sumIndB_1}\cdot \ldots \cdot (\mathbf{a}^{n,\sumIndA})_{\sumIndB_n} \bigl]\eta_{\mPoint}^\sumIndA \Bigr)_{\sumIndB_1, \ldots, \sumIndB_n=1}^{\dimInd_1, \ldots, \dimInd_n}, 
\end{equation}
where $\eta_{\mPoint}^\sumIndA \in \tangent_{\mPoint} \manifold,\; \mathbf{a}^{\sumIndC, \sumIndA}\in \Real^{\dimInd_\sumIndC}$ for $\sumIndA= 1, \ldots, r$, $\sumIndC = 1, \ldots, n$.

However, already in the real-valued case actually computing the CP decomposition of a tensor (or the CP rank for that matter) is NP hard, although some approximation guarantees exist \cite{song2019relative}. An alternative has been proposed by Tucker \cite{tucker1963implications,tucker1966some} in 1963, in which tensor is decomposed into a core tensor multiplied by a matrix along each mode. Algorithms for computing specific instances of the Tucker decomposition, i.e., with core tensors of a specific size, go by many names besides the originally coined three-mode factor analysis \cite{tucker1966some}. Most notable are the three-mode or more generally N-mode PCA \cite{kroonenberg1980principal,kapteyn1986approach}, the higher-order SVD \cite{de2000multilinear} and the N-mode SVD \cite{vasilescu2002multilinear}. 

The multi-linear rank underlying the Tucker decomposition can also be generalized to manifold-valued tensors, but should not be confused with the idea of rank (i.e., the minimum number of rank-one components). That is, we say that a a manifold-valued tensor $\Tensor\in \manifold^{\dimInd_1\times \ldots \times \dimInd_n}$ has \emph{multi-linear rank $\mathbf{r}:=(r_1, \ldots, r_n)\in \Natural^n$ at a point $\mPoint\in \manifold$} if we can write it as
\begin{equation}
    \log_{\mPoint} \Tensor = \sum_{\sumIndA_1, \ldots, \sumIndA_n=1}^{r_1, \ldots, r_n}  (\mathbf{u}^{1,\sumIndA_1} \otimes \ldots \otimes \mathbf{u}^{n,\sumIndA_n})(\mathbfcal{R}_{\mPoint})_{\sumIndA_1, \ldots, \sumIndA_n} := \Bigl(\sum_{\sumIndA_1, \ldots, \sumIndA_n=1}^{r_1, \ldots, r_n}  \bigl[(\mathbf{u}^{1,\sumIndA_1})_{\sumIndB_1}\cdot \ldots \cdot (\mathbf{u}^{n,\sumIndA_n})_{\sumIndB_n} \bigl](\mathbfcal{R}_{\mPoint})_{\sumIndA_1, \ldots, \sumIndA_n} \Bigr)_{\sumIndB_1, \ldots, \sumIndB_n=1}^{\dimInd_1, \ldots, \dimInd_n}, 
\end{equation}
where $\mathbfcal{R}_{\mPoint} \in \tangent_{\mPoint} \manifold^{r_1 \times \ldots \times r_n},\; \mathbf{u}^{\sumIndC,\sumIndA_\sumIndC}\in \Real^{\dimInd_\sumIndC}$ and $\sumIndA_\sumIndC = 1, 
\ldots, r_\sumIndC$, $\sumIndC = 1, \ldots, n$. An equivalent and more convenient notation is
\begin{equation}
    \log_{\mPoint} \Tensor = \mathbfcal{R}_{\mPoint} \times_1 \mathbf{U}^1 \times_2 \ldots \times_n \mathbf{U}^n = \mathbfcal{R}_{\mPoint} \times_{\sumIndC=1}^n \mathbf{U}^\sumIndC,
    \label{eq:mfld-valued-multi-rank-decom}
\end{equation}
where $\mathbf{U}^\sumIndC \in \Real^{\dimInd_\sumIndC \times r_\sumIndC}$.



As to computing a decomposition of the form \cref{eq:mfld-valued-multi-rank-decom}, algorithms in \cite{de2000multilinear,kapteyn1986approach,kroonenberg1980principal,tucker1966some,vasilescu2002multilinear} for the real-valued case assume a specific multi-linear rank. In particular, let the matrix $\mathbfcal{C}_{(\sumIndC)}\in \Real^{\dimInd_\sumIndC \times \prod_{\ell \neq \sumIndC}^n \dimInd_\ell}$ be the mode-$\sumIndC$ unfolding of the real-valued tensor $\mathbfcal{C}\in \Real^{\dimInd_1 \times \ldots\times \dimInd_n}$. Then, these methods find a multi-linear rank-$(R_1, \ldots, R_n)$ decomposition, where $R_\sumIndC := \operatorname{rank}(\mathbfcal{C}_{(\sumIndC)})$. For the manifold-valued case something similar can be constructed. Now, let $\Tensor_{(\sumIndC)} \in (\manifold^{\dimInd_1 \times \ldots\times \dimInd_{\sumIndC-1}\times  \dimInd_{\sumIndC+1}\times \ldots \times \dimInd_n})^{\dimInd_\sumIndC}$ be the mode-$\sumIndC$ unfolding of manifold-valued tensor $\Tensor \in \manifold^{\dimInd_1 \times \cdots \times \dimInd_n}$. Next, let $D_\sumIndC:= \dimInd \cdot \prod_{\ell \neq \sumIndC}^n \dimInd_\ell$ where $\dimInd:= \dim(\manifold)$, consider any orthonormal basis $\{\Phi_{\mPoint}^{J_{\sumIndC}}\}_{J_{\sumIndC}=1}^{D_\sumIndC} \subset \tangent_\mPoint\manifold^{\dimInd_1 \times \ldots\times \dimInd_{\sumIndC-1}\times  \dimInd_{\sumIndC+1}\times \ldots \times \dimInd_n}$ and define $\mathbf{X}^{\sumIndC}\in \Real^{\dimInd_\sumIndC \times D_\sumIndC}$
\begin{equation}
    (\mathbf{X}^{\sumIndC})_{\sumIndB_\sumIndC, J_\sumIndC} = ((\log_\mPoint \Tensor_{(\sumIndC)})_{\sumIndB_\sumIndC}, \Phi_{\mPoint}^{J_\sumIndC})_\mPoint, \quad \text{where } \sumIndB_\sumIndC = 1, \ldots, \dimInd_k, \; J_\sumIndC = 1, \ldots , D_\sumIndC.
\end{equation}
Note that 
\begin{equation}
    (\log_\mPoint \Tensor_{(\sumIndC)})_{\sumIndB_\sumIndC} = \sum_{J_\sumIndC=1}^{D_\sumIndC} (\mathbf{X}^{\sumIndC})_{\sumIndB_\sumIndC, J_\sumIndC} \Phi_{\mPoint}^{J_\sumIndC}.
    \label{eq:Xk-in-coordinates}
\end{equation}
The matrix $\mathbf{X}^{\sumIndC}$ allows us to compute the tangent space SVD for the unfolded $\Tensor_{(\sumIndC)}$ at $\mPoint$. That is, we first compute the SVD of $\mathbf{X}^\sumIndC$, which gives
\begin{equation}
    \mathbf{X}^\sumIndC = \mathbf{U}^\sumIndC \Sigma^\sumIndC (\mathbf{W}^\sumIndC)^\top,
\end{equation}
where $\mathbf{W}^\sumIndC \in \Real^{D_\sumIndC \times R_\sumIndC}$, $\Sigma^\sumIndC = \operatorname{diag}(\sigma_1, \ldots, \sigma_{R_\sumIndC})\in \Real^{R_\sumIndC \times R_\sumIndC}$ with $\sigma_1\geq \ldots \geq \sigma_{R_\sumIndC}$, $\mathbf{U}^\sumIndC\in \Real^{\dimInd_\sumIndC \times R_\sumIndC}$ and where $R_\sumIndC := \operatorname{rank}(\mathbf{X}^\sumIndC)$. Next, we can decompose $\log_\mPoint \Tensor_{(\sumIndC)}$ into factors independent of the basis $\{ \Phi_{\mPoint}^{J_{\sumIndC}} \}_{J_{\sumIndC}=1}^{D_\sumIndC}$ as
\begin{equation}
    \log_\mPoint \Tensor_{(\sumIndC)} = \mathbfcal{W}_\mPoint^\sumIndC \times_1 \Sigma^\sumIndC \times_1 \mathbf{U}^\sumIndC
\end{equation}
where $\mathbfcal{W}^\sumIndC_\mPoint \in\tangent_\mPoint(\manifold^{\dimInd_1 \times \ldots\times \dimInd_{\sumIndC-1}\times  \dimInd_{\sumIndC+1}\times \ldots \times \dimInd_n})^{R_\sumIndC}$ is given by $\mathbfcal{W}^\sumIndC_\mPoint = \sum_{J_\sumIndC=1}^{D_\sumIndC} (\mathbf{W}^\sumIndC)_{J_\sumIndC,\sumIndB_\sumIndC} \Phi_{\mPoint}^{J_\sumIndC}$. 

\begin{remark}
Note that the choice of orthonormal basis $\{ \Phi_{\mPoint}^{J_{\sumIndC}}\}_{J_{\sumIndC}=1}^{D_\sumIndC}$ is truly irrelevant as $\Sigma^\sumIndC$ and $\mathbf{U}^\sumIndC$ in the SVD of $\mathbf{X}^\sumIndC$ can be constructed from the eigendecomposition of $\mathbf{X}^\sumIndC (\mathbf{X}^\sumIndC)^\top$, which is independent of the choice of basis
\begin{equation}
    (\mathbf{X}^\sumIndC (\mathbf{X}^\sumIndC)^\top )_{\sumIndB_\sumIndC, \sumIndB_\sumIndC'} = \sum_{J_\sumIndC=1}^{D_\sumIndC} (\log_\mPoint \Tensor_{\sumIndB_\sumIndC}, \Phi_{\mPoint}^{J_\sumIndC})_\mPoint (\log_\mPoint \Tensor_{\sumIndB_\sumIndC'}, \Phi_{\mPoint}^{J_\sumIndC})_\mPoint = (\log_\mPoint \Tensor_{\sumIndB_\sumIndC},\log_\mPoint \Tensor_{\sumIndB_\sumIndC'})_\mPoint.
\end{equation}
Subsequently, $\mathbfcal{W}^\sumIndC_\mPoint$ also satisfies $\mathbfcal{W}^\sumIndC_\mPoint = \log_\mPoint \Tensor_{(\sumIndC)} \times_1 (\mathbf{U}^\sumIndC)^\top \times_1 (\Sigma^\sumIndC)^{-1}$, which is also independent of basis.
\end{remark}

Repeating the process outlined above yields a generalization to the HOSVD, which we summarize in \cref{alg:tangent-space-HOSVD}.

\begin{algorithm}[h!]
\caption{Tangent space HOSVD}
\label{alg:tangent-space-HOSVD}
\begin{algorithmic}
\STATE{\textit{Initialisation}: $\Tensor\in \manifold^{\dimInd_1\times \ldots\times \dimInd_n}$, $\mPoint\in\manifold$}
\FOR{$\sumIndC = 1, \ldots, n$}
\STATE Unfold $\Tensor$ along mode $\sumIndC$ to get $\Tensor_{(\sumIndC)}$
\STATE Compute the tangent space SVD $\log_\mPoint \Tensor_{(\sumIndC)} = \mathbfcal{W}_\mPoint^\sumIndC \times_1 \Sigma^\sumIndC \times_1 \mathbf{U}^\sumIndC$
\ENDFOR
\STATE{$\mathbfcal{R}_{\mPoint}:= \log_{\mPoint} \Tensor \times_{1} \mathbf{U}_1^\top \times_{2} \ldots \times_{n} \mathbf{U}_n^\top$}
\RETURN $\mathbfcal{R}_{\mPoint} \in \tangent_{\mPoint} \manifold^{R_1 \times \ldots \times R_n}$ and $\{\mathbf{U}^\sumIndC\}_{\sumIndC=1}^n$ with $\mathbf{U}^\sumIndC \in \Real^{\dimInd_\sumIndC \times R_\sumIndC}$
\end{algorithmic}
\end{algorithm}

\subsection{Curvature corrected low multi-linear rank approximation.} 
\label{sec:curvature-correcteed-hosvd}
A low multi-linear rank-$(r_1, \ldots, r_n)$ approximation of a manifold-valued tensor $\Tensor\in \manifold^{\dimInd_1, \ldots, \dimInd_n}$ with $r_\sumIndC < R_\sumIndC = \operatorname{rank}(\mathbf{X}^\sumIndC)$ for $\mathbf{X}^\sumIndC$ from \cref{eq:Xk-in-coordinates} is in the real-valued case typically obtained from truncation of the HOSVD -- even though this is typically non-optimal due to the absence of an Eckart-Young-like theorem for the HOSVD. We intend to take a similar approach for the manifold-valued case, but note that after truncation the low multi-linear rank tensor is also suffering from curvature effects. In particular, in the following we will use the truncated tangent space HOSVD in a curvature correction scheme, which leads to the curvature corrected truncated HOSVD (CC-tHOSVD) (\cref{alg:curvature-corrected-HOSVD}). 

As mentioned in the section introduction, we are aiming for a scheme that is generally applicable to arbitrary manifolds, numerically feasible and global-geometry aware. Arguably, the curvature corrected problem 
\begin{equation}
    \inf_{\mTVector_{\mPoint}\in \mathcal{S}^{\mathbf{r}}_{\mPoint}(\manifold^{\dimInd_1\times \ldots \times \dimInd_n})} \Bigl\{ \sum_{\sumIndB_1, \ldots, \sumIndB_n, \sumIndB=1}^{\dimInd_1, \ldots, \dimInd_n, \dimInd} \beta(\kappa_{\sumIndB_1, \ldots, \sumIndB_n, \sumIndB})^2 \Bigl( \mTVector_{\mPoint} - \log_{\mPoint}\Tensor,  \Theta^{\sumIndB_1, \ldots, \sumIndB_n, \sumIndB}_{\mPoint} \Bigr)_{\mPoint}^2 \Bigr\},
    \label{eq:inf-Sr-linearized-tangent-approximation}
\end{equation}
where 
\begin{equation}
    \mathcal{S}^{\mathbf{r}}_{\mPoint}(\manifold^{\dimInd_1\times \ldots \times \dimInd_n}) := \{ \mTVector_\mPoint = \mathbfcal{R}_{\mPoint} \times_{\sumIndC=1}^n \mathbf{U}^\sumIndC \in  \tangent_{\mPoint} \manifold^{\dimInd_1\times \ldots\times \dimInd_n}\; \mid \;  \mathbfcal{R}_{\mPoint} \in \tangent_{\mPoint} \manifold^{r_1 \times \ldots \times r_n}, \mathbf{U}^\sumIndC \in \Real^{\dimInd_\sumIndC \times r_\sumIndC}\}
    \label{eq:Sr}
\end{equation}
is reasonably generally applicable as symmetric Riemannian manifolds include many practically used manifolds, and is global-geometry aware as we have seen in \cref{sec:diminishing-discrepancy} that the discrepancy with the global approximation error diminishes under the mild condition that all $\kappa_{\sumIndB_1, \ldots, \sumIndB_n, \sumIndB}< \pi^2$. So it remains to solve \cref{eq:inf-Sr-linearized-tangent-approximation} or approximate a solution efficiently.

The truncated tangent space HOSVD will be a good point to start from. Indeed, since the set \cref{eq:Sr} is closed in the tangent space topology, \cref{thm:convergence-of-minimizers} in \cref{sec:stability-curvature} tells us that the minimizer of \cref{eq:inf-Sr-linearized-tangent-approximation} will be close to the truncated tangent space HOSVD. Nevertheless, because of the non-convexity and high-dimensionality of $\mathcal{S}^{\mathbf{r}}_{\mPoint}(\manifold^{\dimInd_1\times \ldots \times \dimInd_n})$ minimizing over the full space will still be numerically expensive. To alleviate this issue, we will make the additional assumption that the optimal $\hat{\mathbf{U}}^\sumIndC$ are close to the $\mathbf{U}^\sumIndC$ obtained from the truncated HOSVD and that only the core tensor $\mathbfcal{R}_{\mPoint}$ needs correction. In other words, we will solve the problem
\begin{equation}
    \inf_{\mathbfcal{R}_{\mPoint}\in \tangent_{\mPoint}\manifold^{r_1\times \ldots \times r_n}} \Bigl\{ \sum_{\sumIndB_1, \ldots, \sumIndB_n, \sumIndB=1}^{\dimInd_1, \ldots, \dimInd_n, \dimInd} \beta(\kappa_{\sumIndB_1, \ldots, \sumIndB_n, \sumIndB})^2 \Bigl( \mathbfcal{R}_{\mPoint} \times_{\sumIndC=1}^n \mathbf{U}^\sumIndC - \log_{\mPoint}\Tensor,  \Theta^{\sumIndB_1, \ldots, \sumIndB_n, \sumIndB}_{\mPoint} \Bigr)_{\mPoint}^2 \Bigr\},
    \label{eq:inf-Sr-linearized-tangent-approximation-only-R}
\end{equation}
which is a quadratic minimization problem.

Instead of minimizing over $\mathbfcal{R}_{\mPoint}$ directly, we will expand the core tensor into an orthonormal basis, i.e., 
\begin{equation}
    \mathbfcal{R}_{\mPoint} = \mathbfcal{V} \times_{n+1} \phi_{\mPoint} :=\Bigl(\sum_{\sumIndA=1}^\dimInd\mathbfcal{V}_{\sumIndA_1, \ldots, \sumIndA_n, \sumIndA} (\phi_{\mPoint})_{\sumIndA}\Bigr)_{\sumIndA_1, \ldots, \sumIndA_n =1}^{r_1, \ldots, r_n}, 
    \label{eq:reparam-Rp}
\end{equation}
where $\mathbfcal{V}\in \Real^{r_1\times \ldots\times r_n\times \dimInd}$ is a real-valued tensor for $\dimInd:= \dim(\manifold)$ and $\phi_{\mPoint} := (\phi_{\mPoint}^{1}, \dots, \phi_{\mPoint}^{n}) \subset \tangent_\mPoint\manifold^{\dimInd}$ where $\{\phi_{\mPoint}^{\sumIndA}\}_{\sumIndA=1}^\dimInd\subset \tangent_{\mPoint}\manifold$ is an orthonormal basis, so that $(\phi_{\mPoint})_{\sumIndA} = \phi_{\mPoint}^{\sumIndA}$.

Rewriting \cref{eq:inf-Sr-linearized-tangent-approximation-only-R} through insertion of \cref{eq:reparam-Rp} gives
\begin{equation}
    \inf_{\mathbfcal{V}\in \Real^{r_1\times \ldots \times r_n \times \dimInd}} \Bigl\{ \sum_{\sumIndB_1, \ldots, \sumIndB_n, \sumIndB=1}^{\dimInd_1, \ldots, \dimInd_n, \dimInd} \beta(\kappa_{\sumIndB_1, \ldots, \sumIndB_n, \sumIndB})^2 \Bigl( \mathbfcal{V} \times_{\sumIndC =1}^n \mathbf{U}^\sumIndC \times_{n+1} \phi_{\mPoint} - \log_{\mPoint}\Tensor,  \Theta^{\sumIndB_1, \ldots, \sumIndB_n, \sumIndB}_{\mPoint} \Bigr)_{\mPoint}^2 \Bigr\}.
    \label{eq:inf-S-linearized-tangent-approximation-relaxed}
\end{equation}
A natural way of solving \cref{eq:inf-S-linearized-tangent-approximation-relaxed} is through first-order optimality conditions, i.e., a vanishing gradient. Let $f:\Real^{r_1\times \ldots\times r_n\times \dimInd} \to \Real$ be given by 
\begin{equation}
    f(\mathbfcal{V}):= \sum_{\sumIndB_1, \ldots, \sumIndB_n, \sumIndB=1}^{\dimInd_1, \ldots, \dimInd_n, \dimInd} \beta(\kappa_{\sumIndB_1, \ldots, \sumIndB_n, \sumIndB})^2 \Bigl( \mathbfcal{V} \times_{\sumIndC =1}^n \mathbf{U}^\sumIndC \times_{n+1} \phi_{\mPoint} - \log_{\mPoint}\Tensor,  \Theta^{\sumIndB_1, \ldots, \sumIndB_n, \sumIndB}_{\mPoint} \Bigr)_{\mPoint}^2
\end{equation}
The gradient of $f$ satisfies
\begin{equation}
    \nabla f (\mathbfcal{V})_{\sumIndA_1, \ldots, \sumIndA_n, \sumIndA} = 2  \sum_{\sumIndB_1, \ldots, \sumIndB_n, \sumIndB=1}^{\dimInd_1, \ldots, \dimInd_n, \dimInd} \beta(\kappa_{\sumIndB_1, \ldots, \sumIndB_n, \sumIndB})^2 \Bigl( \mathbfcal{V} \times_{\sumIndC =1}^n \mathbf{U}^\sumIndC \times_{n+1} \phi_{\mPoint} - \log_{\mPoint}\Tensor,  \Theta^{\sumIndB_1, \ldots, \sumIndB_n, \sumIndB}_{\mPoint} \Bigr)_{\mPoint} \Bigl( \mathbfcal{E}^{\sumIndA_1, \ldots, \sumIndA_n, \sumIndA} \times_{\sumIndC =1}^n \mathbf{U}^\sumIndC \times_{n+1} \phi_{\mPoint},  \Theta^{\sumIndB_1, \ldots, \sumIndB_n, \sumIndB}_{\mPoint})_{\mPoint}
\end{equation}
where
\begin{equation}
        (\mathbfcal{E}^{\sumIndA_1, \ldots, \sumIndA_n, \sumIndA} )_{\sumIndA_1', \ldots, \sumIndA_n', \sumIndA'} = \left\{\begin{matrix}
        1 & \text{if } \sumIndA_1' = \sumIndA_1, \ldots, \sumIndA_n' = \sumIndA_n, \sumIndA' = \sumIndA,  \\
        0 &  \text{otherwise}.\\
        \end{matrix}\right. 
        \label{eq:Eiiii}
\end{equation}
We want to solve $\nabla f (\mathbfcal{V}) = 0$. For that we will rewrite some terms. First, note that we may assume that 
\begin{equation}
        (\Theta^{\sumIndB_1, \ldots, \sumIndB_n, \sumIndB}_{\mPoint})_{\sumIndB_1', \ldots, \sumIndB_n'} = \left\{\begin{matrix}
        (\theta^{\sumIndB_1, \ldots, \sumIndB_n}_{\mPoint})_{\sumIndB} & \text{if } \sumIndB_1' = \sumIndB_1, \ldots, \sumIndB_n' = \sumIndB_n,  \\
        0 &  \text{otherwise},\\
        \end{matrix}\right. 
        \label{eq:assumption-digonalizing-basis}
\end{equation}
where
\begin{equation}
    \curvature((\theta^{\sumIndB_1, \ldots, \sumIndB_n}_{\mPoint})_{\sumIndB}, \log_\mPoint \Tensor_{\sumIndB_1, \ldots, \sumIndB_n})\log_\mPoint \Tensor_{\sumIndB_1, \ldots, \sumIndB_n} = \kappa_{\sumIndB_1, \ldots, \sumIndB_n, \sumIndB} (\theta^{\sumIndB_1, \ldots, \sumIndB_n}_{\mPoint})_{\sumIndB} \in \tangent_\mPoint\manifold.
    \label{eq:def-theta_j}
\end{equation}
Next, we define the real-valued tensor $\mathbfcal{B} \in \Real^{\dimInd_1 \times \ldots\times \dimInd_n \times \dimInd \times r_1 \times \ldots \times r_n \times \dimInd}$ given by
\begin{equation}
    \mathbfcal{B}_{\sumIndB_1, \ldots, \sumIndB_n, \sumIndB, \sumIndA_1, \ldots, \sumIndA_n, \sumIndA} := (\mathbf{U}^1)_{\sumIndB_{1}, \sumIndA_{1}} \cdot \ldots \cdot (\mathbf{U}^n)_{\sumIndB_{n}, \sumIndA_{n}} \bigl((\phi_{\mPoint})_{\sumIndA}, (\theta_{\mPoint}^{\sumIndB_1, \ldots, \sumIndB_n})_{\sumIndB} \bigr)_{\mPoint}, 
    \label{eq:B-tensor}
\end{equation}
and define the real-valued tensor $\mathbfcal{A}\in \Real^{r_1 \times \ldots \times r_n \times \dimInd \times r_1 \times \ldots \times r_n \times \dimInd}$ given by
\begin{equation}
    \mathbfcal{A}_{\sumIndA_1', \ldots, \sumIndA_n', \sumIndA', \sumIndA_1, \ldots, \sumIndA_n, \sumIndA} := \sum_{\sumIndB_1, \ldots, \sumIndB_n, \sumIndB=1}^{\dimInd_1, \ldots,\dimInd_n, \dimInd}  \beta(\kappa_{\sumIndB_1, \ldots, \sumIndB_n, \sumIndB})^2 \mathbfcal{B}_{\sumIndB_1, \ldots, \sumIndB_n, \sumIndB, \sumIndA_1', \ldots, \sumIndA_n', \sumIndA'} \mathbfcal{B}_{\sumIndB_1, \ldots, \sumIndB_n, \sumIndB, \sumIndA_1, \ldots, \sumIndA_n, \sumIndA}.
    \label{eq:A-tensor}
\end{equation}
Then, it is easy to see that $\nabla f (\mathbfcal{V}) = 0$ is equivalent to solving the linear system
\begin{equation}
    \sum_{\sumIndA_1, \ldots, \sumIndA_n, \sumIndA = 1}^{r_1, \ldots, r_n, \dimInd} \mathbfcal{A}_{\sumIndA_1', \ldots, \sumIndA_n', \sumIndA', \sumIndA_1, \ldots, \sumIndA_n, \sumIndA} \mathbfcal{V}_{\sumIndA_1, \ldots, \sumIndA_n, \sumIndA} = \sum_{\sumIndB_1, \ldots, \sumIndB_n, \sumIndB=1}^{\dimInd_1, \ldots,\dimInd_n, \dimInd} \beta(\kappa_{\sumIndB_1, \ldots, \sumIndB_n, \sumIndB})^2 \mathbfcal{B}_{\sumIndB_1, \ldots, \sumIndB_n, \sumIndB, \sumIndA_1', \ldots, \sumIndA_n', \sumIndA'} \Bigl( \log_\mPoint \Tensor_{\sumIndB_1,\ldots,\sumIndB_n},  (\theta_{\mPoint}^{\sumIndB_1, \ldots, \sumIndB_n})_{\sumIndB} \Bigr)_{\mPoint}
    \label{eq:curvature-correction-system}
\end{equation}
Finally, the above system can be solved uniquely under familiar conditions. The full process gives the curvature corrected truncated HOSVD, which is is summarized in \cref{alg:curvature-corrected-HOSVD}.
\begin{proposition}
    The linear operator $\mathbfcal{A}:\Real^{r_1 \times \ldots \times r_n \times \dimInd} \to \Real^{r_1 \times \ldots \times r_n \times \dimInd}$ given by
    \begin{equation}
        (\mathbfcal{A}(\mathbfcal{V}))_{\sumIndA_1', \ldots, \sumIndA_n', \sumIndA'}:= \sum_{\sumIndA_1, \ldots, \sumIndA_n, \sumIndA = 1}^{r_1, \ldots, r_n, \dimInd} \mathbfcal{A}_{\sumIndA_1', \ldots, \sumIndA_n', \sumIndA', \sumIndA_1, \ldots, \sumIndA_n, \sumIndA} \mathbfcal{V}_{\sumIndA_1, \ldots, \sumIndA_n, \sumIndA}, \quad \text{where $\mathbfcal{A}_{\sumIndA_1', \ldots, \sumIndA_n', \sumIndA', \sumIndA_1, \ldots, \sumIndA_n, \sumIndA}$ from \cref{eq:A-tensor}} 
    \end{equation}
    is symmetric and positive definite if $\kappa_{\sumIndB_1, \ldots, \sumIndB_n, \sumIndB}< \pi^2$ for all $\sumIndB_1 = 1, \ldots, \dimInd_1$ through $\sumIndB_n = 1, \dots, \dimInd_n$ and $\sumIndB= 1, \ldots, \dimInd$. In particular, under these assumptions the linear system \cref{eq:curvature-correction-system} has a unique solution.
\end{proposition}

\begin{proof}
    Symmetry of $\mathbfcal{A}$ follows directly from \cref{eq:A-tensor}. For positive definiteness we will first show that if $\kappa_{\sumIndB_1, \ldots, \sumIndB_n, \sumIndB} \leq \pi^2$
    \begin{equation}
        \mathbfcal{A} \succcurlyeq \beta(\kappa_{\max})^2 \operatorname{Id}, \quad \text{where $\operatorname{Id}: \Real^{r_1 \times \ldots \times r_n \times \dimInd} \to \Real^{r_1 \times \ldots \times r_n \times \dimInd}$ is the identity operator.}
    \end{equation}
    Then, positive definiteness follows from the assumption $\kappa_{\sumIndB_1, \ldots, \sumIndB_n, \sumIndB}< \pi^2$ by definition of $\beta$ in \cref{eq:beta}.

    So it remains to show that if $\kappa_{\sumIndB_1, \ldots, \sumIndB_n, \sumIndB} \leq \pi^2$
    \begin{equation}
        \bigl(\mathbfcal{V}, (\mathbfcal{A} - \beta(\kappa_{\max})^2 \operatorname{Id})(\mathbfcal{V})\bigr) \geq 0, \quad \text{for any $\mathbfcal{V} \in \Real^{r_1 \times \ldots \times r_n \times \dimInd}$.}
        \label{eq:claim-prop-inv-A}
    \end{equation}

     First, we define the linear operator $\mathbfcal{B}: \Real^{r_1 \times \ldots \times r_n \times \dimInd} \to \Real^{\dimInd_1 \times \ldots \times \dimInd_n \times \dimInd}$
    \begin{equation}
        (\mathbfcal{B}(\mathbfcal{V}))_{\sumIndB_1, \ldots, \sumIndB_n, \sumIndB}:= \sum_{\sumIndA_1, \ldots, \sumIndA_n, \sumIndA = 1}^{r_1, \ldots, r_n, \dimInd} \mathbfcal{B}_{\sumIndB_1, \ldots, \sumIndB_n, \sumIndB, \sumIndA_1, \ldots, \sumIndA_n, \sumIndA} \mathbfcal{V}_{\sumIndA_1, \ldots, \sumIndA_n, \sumIndA}, \quad \text{where $\mathbfcal{B}_{\sumIndB_1, \ldots, \sumIndB_n, \sumIndB, \sumIndA_1, \ldots, \sumIndA_n, \sumIndA}$ from \cref{eq:B-tensor}},
    \end{equation}
    Next, we note that
    \begin{multline}
        \sum_{\sumIndB_1, \ldots, \sumIndB_n, \sumIndB=1}^{\dimInd_1, \ldots,\dimInd_n, \dimInd}  \beta(\kappa_{\max})^2 \mathbfcal{B}_{\sumIndB_1, \ldots, \sumIndB_n, \sumIndB, \sumIndA_1', \ldots, \sumIndA_n', \sumIndA'} \mathbfcal{B}_{\sumIndB_1, \ldots, \sumIndB_n, \sumIndB, \sumIndA_1, \ldots, \sumIndA_n, \sumIndA} \\
        = \beta(\kappa_{\max})^2 \sum_{\sumIndB_1, \ldots, \sumIndB_n=1}^{\dimInd_1, \ldots,\dimInd_n} \sum_{\sumIndB=1}^\dimInd (\mathbf{U}^1)_{\sumIndB_{1}, \sumIndA_{1}'} \cdot \ldots \cdot (\mathbf{U}^n)_{\sumIndB_{n}, \sumIndA_{n}'} \bigl((\phi_{\mPoint})_{\sumIndA'}, (\theta_{\mPoint}^{\sumIndB_1, \ldots, \sumIndB_n})_{\sumIndB} \bigr)_{\mPoint} (\mathbf{U}^1)_{\sumIndB_{1}, \sumIndA_{1}} \cdot \ldots \cdot (\mathbf{U}^n)_{\sumIndB_{n}, \sumIndA_{n}} \bigl((\phi_{\mPoint})_{\sumIndA}, (\theta_{\mPoint}^{\sumIndB_1, \ldots, \sumIndB_n})_{\sumIndB} \bigr)_{\mPoint}
        \\
        = \beta(\kappa_{\max})^2 \sum_{\sumIndB_1, \ldots, \sumIndB_n=1}^{\dimInd_1, \ldots,\dimInd_n} (\mathbf{U}^1)_{\sumIndB_{1}, \sumIndA_{1}'} \cdot \ldots \cdot (\mathbf{U}^n)_{\sumIndB_{n}, \sumIndA_{n}'} (\mathbf{U}^1)_{\sumIndB_{1}, \sumIndA_{1}} \cdot \ldots \cdot (\mathbf{U}^n)_{\sumIndB_{n}, \sumIndA_{n}} \sum_{\sumIndB=1}^\dimInd \bigl((\phi_{\mPoint})_{\sumIndA'}, (\theta_{\mPoint}^{\sumIndB_1, \ldots, \sumIndB_n})_{\sumIndB} \bigr)_{\mPoint}  \bigl((\phi_{\mPoint})_{\sumIndA}, (\theta_{\mPoint}^{\sumIndB_1, \ldots, \sumIndB_n})_{\sumIndB} \bigr)_{\mPoint}\\
        = \beta(\kappa_{\max})^2 \sum_{\sumIndB_1, \ldots, \sumIndB_n=1}^{\dimInd_1, \ldots,\dimInd_n} (\mathbf{U}^1)_{\sumIndB_{1}, \sumIndA_{1}'} \cdot \ldots \cdot (\mathbf{U}^n)_{\sumIndB_{n}, \sumIndA_{n}'} (\mathbf{U}^1)_{\sumIndB_{1}, \sumIndA_{1}} \cdot \ldots \cdot (\mathbf{U}^n)_{\sumIndB_{n}, \sumIndA_{n}} \bigl((\phi_{\mPoint})_{\sumIndA'}, (\phi_{\mPoint})_{\sumIndA} \bigr)_{\mPoint}  \\
        = \beta(\kappa_{\max})^2 \delta_{\sumIndA',\sumIndA} \sum_{\sumIndB_1, \ldots, \sumIndB_n=1}^{\dimInd_1, \ldots,\dimInd_n} (\mathbf{U}^1)_{\sumIndB_{1}, \sumIndA_{1}'} \cdot \ldots \cdot (\mathbf{U}^n)_{\sumIndB_{n}, \sumIndA_{n}'} (\mathbf{U}^1)_{\sumIndB_{1}, \sumIndA_{1}} \cdot \ldots \cdot (\mathbf{U}^n)_{\sumIndB_{n}, \sumIndA_{n}} \\
        =  \beta(\kappa_{\max})^2 \delta_{\sumIndA',\sumIndA} \sum_{\sumIndB_1=1}^{\dimInd_1} (\mathbf{U}^1)_{\sumIndB_{1}, \sumIndA_{1}'} \cdot (\mathbf{U}^1)_{\sumIndB_{1}, \sumIndA_{1}} \cdot \ldots \cdot \sum_{\sumIndB_n=1}^{\dimInd_n}(\mathbf{U}^n)_{\sumIndB_{n}, \sumIndA_{n}'} \cdot (\mathbf{U}^n)_{\sumIndB_{n}, \sumIndA_{n}}\\
        = \beta(\kappa_{\max})^2 \delta_{\sumIndA_1',\sumIndA_1} \cdot \ldots \cdot \delta_{\sumIndA_n',\sumIndA_n}  \cdot \delta_{\sumIndA',\sumIndA},
        \label{eq:prop-inv-B-reduction}
    \end{multline}
    which allows us to evaluate
    \begin{multline}
        \|\beta(\kappa_{\max}) \mathbfcal{B}(\mathbfcal{V})\|_{F}^2 \\
        = \sum_{\sumIndB_1, \ldots, \sumIndB_n, \sumIndB=1}^{\dimInd_1, \ldots,\dimInd_n, \dimInd} \sum_{\sumIndA_1, \ldots, \sumIndA_n, \sumIndA = 1}^{r_1, \ldots, r_n, \dimInd} \sum_{\sumIndA'_1, \ldots, \sumIndA'_n, \sumIndA' = 1}^{r_1, \ldots, r_n, \dimInd}  \beta(\kappa_{\max})^2 \mathbfcal{B}_{\sumIndB_1, \ldots, \sumIndB_n, \sumIndB, \sumIndA_1', \ldots, \sumIndA_n', \sumIndA'} \mathbfcal{B}_{\sumIndB_1, \ldots, \sumIndB_n, \sumIndB, \sumIndA_1, \ldots, \sumIndA_n, \sumIndA} \mathbfcal{V}_{\sumIndA_1, \ldots, \sumIndA_n, \sumIndA} \mathbfcal{V}_{\sumIndA'_1, \ldots, \sumIndA'_n, \sumIndA'}\\
        \overset{\text{\cref{eq:prop-inv-B-reduction}}}{=} \sum_{\sumIndA_1, \ldots, \sumIndA_n, \sumIndA = 1}^{r_1, \ldots, r_n, \dimInd} \sum_{\sumIndA'_1, \ldots, \sumIndA'_n, \sumIndA' = 1}^{r_1, \ldots, r_n, \dimInd} \beta(\kappa_{\max})^2 \delta_{\sumIndA_1',\sumIndA_1} \cdot \ldots \cdot \delta_{\sumIndA_n',\sumIndA_n}  \cdot \delta_{\sumIndA',\sumIndA}\mathbfcal{V}_{\sumIndA_1, \ldots, \sumIndA_n, \sumIndA} \mathbfcal{V}_{\sumIndA'_1, \ldots, \sumIndA'_n, \sumIndA'} \\
        = \sum_{\sumIndA'_1, \ldots, \sumIndA'_n, \sumIndA' = 1}^{r_1, \ldots, r_n, \dimInd}\beta(\kappa_{\max})^2 \mathbfcal{V}_{\sumIndA_1', \ldots, \sumIndA_n', \sumIndA'}^2 = \|\beta(\kappa_{\max})\mathbfcal{V}\|_{F}^2.
        \label{eq:prop-invertible-A-Id-rewrite}
    \end{multline}
    Additionally, we can evaluate
    \begin{equation}
        \|\beta(\kappa) \odot \mathbfcal{B}(\mathbfcal{V})\|_{F}^2 :=  \sum_{\sumIndB_1, \ldots, \sumIndB_n, \sumIndB=1}^{\dimInd_1, \ldots,\dimInd_n, \dimInd} \beta(\kappa_{\sumIndB_1, \ldots, \sumIndB_n, \sumIndB})^2 (\mathbfcal{B}(\mathbfcal{V}))_{\sumIndB_1, \ldots, \sumIndB_n, \sumIndB}^2  = (\mathbfcal{V}, \mathbfcal{A}(\mathbfcal{V})),
        \label{eq:prop-invertible-A-rewrite}
    \end{equation}
    where $\odot$ denotes the Hadamard product or entry-wise multiplication, and where $\beta(\kappa) \in \Real^{\dimInd_1 \times \ldots \times \dimInd_n \times \dimInd}$ is given by
    \begin{equation}
        \beta(\kappa)_{\sumIndB_1, \ldots, \sumIndB_n, \sumIndB} := \beta(\kappa_{\sumIndB_1, \ldots, \sumIndB_n, \sumIndB}).
    \end{equation}
    
    Finally, under the assumption $\kappa_{\sumIndB_1, \ldots, \sumIndB_n, \sumIndB} \leq \pi^2$ we have that $\beta(\kappa_{\sumIndB_1, \ldots, \sumIndB_n, \sumIndB}) \geq \beta(\kappa_{\max})$, because the function $\beta$ from \cref{eq:beta} is decreasing on $(-\infty, \pi^2]$. 
    
    Bringing everything together gives
    \begin{multline}
        (\mathbfcal{V}, (\mathbfcal{A} - \beta(\kappa_{\max})^2 \operatorname{Id})(\mathbfcal{V})) = \bigl(\mathbfcal{V}, \mathbfcal{A}(\mathbfcal{V})\bigr) - \|\beta(\kappa_{\max})\mathbfcal{V}\|_{F}^2 \overset{\text{\cref{eq:prop-invertible-A-Id-rewrite,eq:prop-invertible-A-rewrite}}}{=} \|\beta(\kappa) \odot \mathbfcal{B}(\mathbfcal{V})\|_{F}^2 - \|\beta(\kappa_{\max}) \mathbfcal{B}(\mathbfcal{V})\|_{F}^2\\
        = \sum_{\sumIndB_1, \ldots, \sumIndB_n, \sumIndB=1}^{\dimInd_1, \ldots,\dimInd_n, \dimInd} (\beta(\kappa_{\sumIndB_1, \ldots, \sumIndB_n, \sumIndB})^2 -\beta(\kappa_{\max})^2 ) (\mathbfcal{B}(\mathbfcal{V}))_{\sumIndB_1, \ldots, \sumIndB_n, \sumIndB}^2
        = \|\sqrt{\beta(\kappa)^2 -\beta(\kappa_{\max})^2 }\odot \mathbfcal{B}(\mathbfcal{V})\|_{F}^2 \geq 0,
    \end{multline}
    where $\sqrt{\beta(\kappa)^2 -\beta(\kappa_{\max})^2 } \in \Real^{\dimInd_1 \times \ldots\times \dimInd_n \times \dimInd}$ given by
    \begin{equation}
        (\sqrt{\beta(\kappa)^2 -\beta(\kappa_{\max})^2 })_{\sumIndB_1, \ldots, \sumIndB_n, \sumIndB} := \sqrt{\beta(\kappa_{\sumIndB_1, \ldots, \sumIndB_n, \sumIndB})^2 -\beta(\kappa_{\max})^2 },
    \end{equation}
    which yields the claim \cref{eq:claim-prop-inv-A}. The uniqueness of the solution of the linear system follows from the fact that symmetry with positive definiteness implies invertibility.
\end{proof}

\begin{algorithm}[h!]
\caption{Curvature corrected truncated HOSVD (CC-tHOSVD)}
\label{alg:curvature-corrected-HOSVD}
\begin{algorithmic}
\STATE{\textit{Initialisation}: $\Tensor\in \manifold^{\dimInd_1\times  \ldots\times \dimInd_n}$, $\mPoint\in\manifold$, $\mathbf{r}\in \Natural^n$ }
\STATE Compute the tangent space HOSVD $\log_{\mPoint} \Tensor = \tilde{\mathbfcal{R}}_{\mPoint} \times_{\sumIndC=1}^n \tilde{\mathbf{U}}^\sumIndC$ using  \cref{alg:tangent-space-HOSVD}
\FOR{$\sumIndC = 1, \ldots, n$}
\STATE Construct the truncated $\mathbf{U}^\sumIndC \in \Real^{\dimInd_\sumIndC \times r_\sumIndC}$ from the first $r_\sumIndC$ columns of $\tilde{\mathbf{U}}^\sumIndC \in \Real^{\dimInd_\sumIndC \times R_\sumIndC}$
\ENDFOR
\STATE Construct the eigenpairs $(\kappa_{\sumIndB_1, \ldots, \sumIndB_n, \sumIndB},(\theta_{\mPoint}^{\sumIndB_1, \ldots, \sumIndB_n})_{\sumIndB})$ satisfying \cref{eq:def-theta_j}
\STATE Compute $\mathbfcal{V}$ through solving the linear system \cref{eq:curvature-correction-system} for any orthonormal basis  $\{\phi_{\mPoint}^{\sumIndA}\}_{\sumIndA=1}^\dimInd\subset \tangent_{\mPoint}\manifold$
\STATE $\mathbfcal{R}_{\mPoint} := \mathbfcal{V} \times_{n+1} \phi_{\mPoint} $
\RETURN $\mathbfcal{R}_{\mPoint} \in \tangent_{\mPoint} \manifold^{r_1 \times \ldots \times r_n}$ and $\{\mathbf{U}^\sumIndC\}_{\sumIndC=1}^n$ with $\mathbf{U}^\sumIndC \in \Real^{\dimInd_\sumIndC \times r_\sumIndC}$
\end{algorithmic}
\end{algorithm}

\begin{remark}
\label{rem:using-lower-bound}
    By \cref{cor:approximation-error-lower-bound} we know that the minimizing energy of the full problem \cref{eq:inf-Sr-linearized-tangent-approximation} is bounded from below. Additionally, \cref{eq:inf-Sr-linearized-tangent-approximation} is bounded from above by the minimizing energy of the restricted problem \cref{eq:inf-Sr-linearized-tangent-approximation-only-R}. In other words, a small discrepancy between this upper and lower bound can in practice be used as numerical validation of the assumption on the unimportance of also optimizing over $\mathbf{U}^\sumIndC$.
\end{remark}

\begin{remark}
    Although not being the focus of this work, we would like to point out that a similar scheme can be constructed for low CP-rank approximation. In particular, after having computed and truncated all $\mathbf{a}^{\sumIndC, \sumIndA}$, only the tangent vectors $\eta_{\mPoint}^\sumIndA$ would be corrected for curvature.
\end{remark}


\subsection{Metric corrected low-multilinear rank approximation.}
\label{sec:metric-corrected-hosvd}

As mentioned in the introduction we will conclude this section with constructing the metric corrected truncated HOSVD (MC-tHOSVD). This scheme is similar in spirit to the CC-tHOSVD, but approximates the minimizer of the global approximation error instead of the curvature corrected error. That is, instead of solving
\begin{equation}
    \inf_{\mTVector_{\mPoint}\in \mathcal{S}^{\mathbf{r}}_{\mPoint}(\manifold^{\dimInd_1\times \ldots \times \dimInd_n})} \Bigl\{ \distance_{\manifold^{\dimInd_1 \times \cdots \times \dimInd_{n}}}(\Tensor, \exp_{\mPoint}(\mTVector_{\mPoint}))^2 \Bigr\},
    \label{eq:global-approximation-error-minimization-low-rank}
\end{equation}
we will solve the restricted problem
\begin{equation}
    \inf_{\mathbfcal{V}\in \Real^{r_1\times \ldots \times r_n \times \dimInd}} \Bigl\{ \distance_{\manifold^{\dimInd_1 \times \cdots \times \dimInd_{n}}}(\Tensor, \exp_{\mPoint}(\mathbfcal{V} \times_{\sumIndC =1}^n \mathbf{U}^\sumIndC \times_{n+1} \phi_{\mPoint}))^2 \Bigr\},
    \label{eq:global-approximation-error-minimization-relaxed}
\end{equation}
where $\mathbf{U}^\sumIndC \in \Real^{\dimInd_\sumIndC \times r_\sumIndC}$ are the truncated columns of the HOSVD solution and $\{\phi_{\mPoint}^{\sumIndA}\}_{\sumIndA=1}^\dimInd\subset \tangent_{\mPoint}\manifold$ is any orthonormal basis.

Unlike the restricted problem \cref{eq:inf-S-linearized-tangent-approximation-relaxed}, the problem \cref{eq:global-approximation-error-minimization-relaxed} above does not have a closed-form solution. Instead, we resort to solving it with gradient descent, which can be computed reasonably efficiently under the symmetric Riemannian manifold assumption. That is, let $g:\Real^{r_1\times \ldots\times r_n\times \dimInd} \to \Real$ be given by 
\begin{equation}
    g(\mathbfcal{V}):= \distance_{\manifold^{\dimInd_1 \times \cdots \times \dimInd_{n}}}(\Tensor, \exp_{\mPoint}(\mathbfcal{V} \times_{\sumIndC =1}^n \mathbf{U}^\sumIndC \times_{n+1} \phi_{\mPoint}))^2.
\end{equation}
The gradient of $g$ can be evaluated on symmetric Riemannian manifolds, i.e.,
\begin{multline}
    \nabla g (\mathbfcal{V})_{\sumIndA_1, \ldots, \sumIndA_n, \sumIndA} \overset{\text{chain rule}}{=} -2\Bigl(\log_{\TensorB}\Tensor, D_{\log_\mPoint \TensorB} \exp_{\mPoint} [\mathbfcal{E}^{\sumIndA_1, \ldots, \sumIndA_n, \sumIndA} \times_{\sumIndC =1}^n \mathbf{U}^\sumIndC \times_{n+1} \phi_{\mPoint}]\Bigr)_{\TensorB}\\
    \overset{\text{\cite[Lemma~1]{bergmann2019recent}}}{=} - 2  \sum_{\sumIndB_1, \ldots, \sumIndB_n, \sumIndB=1}^{\dimInd_1, \ldots, \dimInd_n, \dimInd} \beta(\lambda_{\sumIndB_1, \ldots, \sumIndB_n, \sumIndB}) \Bigl(\log_{\TensorB}\Tensor, \Psi^{\sumIndB_1, \ldots, \sumIndB_n, \sumIndB}_{\TensorB}\Bigr)_{\TensorB}\Bigl( \mathbfcal{E}^{\sumIndA_1, \ldots, \sumIndA_n, \sumIndA} \times_{\sumIndC =1}^n \mathbf{U}^\sumIndC \times_{n+1} \phi_{\mPoint},  \Psi^{\sumIndB_1, \ldots, \sumIndB_n, \sumIndB}_{\mPoint}\Bigr)_{\mPoint},
    \label{eq:g-grad-symmetric}
\end{multline}
where $\TensorB:= \exp_\mPoint(\mathbfcal{V} \times_{\sumIndC=1}^n \mathbf{U}^\sumIndC \times_{n+1}\phi_\mPoint)$, $\mathbfcal{E}^{\sumIndA_1, \ldots, \sumIndA_n, \sumIndA}$ as in \cref{eq:Eiiii}, $\Psi^{\sumIndB_1, \ldots, \sumIndB_n, \sumIndB}_{\mPoint}$ such that 
\begin{equation}
    \curvature (\Psi^{\sumIndB_1, \ldots, \sumIndB_n, \sumIndB}_{\mPoint}, \log_{\mPoint} \TensorB)\log_{\mPoint} \TensorB = \lambda_{\sumIndB_1, \ldots, \sumIndB_n, \sumIndB} \Psi^{\sumIndB_1, \ldots, \sumIndB_n, \sumIndB}_{\mPoint},
\end{equation}
and $\Psi^{\sumIndB_1, \ldots, \sumIndB_n, \sumIndB}_{\TensorB} := \mathcal{P}_{\TensorB\leftarrow\mPoint}\Psi^{\sumIndB_1, \ldots, \sumIndB_n, \sumIndB}_\mPoint$.



The full scheme is summarized in \cref{alg:geometry-corrected-HOSVD}. 

\begin{algorithm}[h!]
\caption{Metric corrected truncated HOSVD (MC-tHOSVD)}
\label{alg:geometry-corrected-HOSVD}
\begin{algorithmic}
\STATE{\textit{Initialisation}: $\mathbf{r}\in \Natural^n$ $\Tensor\in \manifold^{\dimInd_1 \times \ldots\times \dimInd_n}$, $\mPoint\in\manifold$, $\tau>0$ $\iterInd:=0$}
\STATE Compute the tangent space HOSVD $\log_{\mPoint} \Tensor = \tilde{\mathbfcal{R}}_{\mPoint} \times_{\sumIndC=1}^n \tilde{\mathbf{U}}^\sumIndC$ using  \cref{alg:tangent-space-HOSVD}
\FOR{$\sumIndC = 1, \ldots, n$}
\STATE Construct the truncated $\mathbf{U}^\sumIndC \in \Real^{\dimInd_\sumIndC \times r_\sumIndC}$ from the first $r_\sumIndC$ columns of $\tilde{\mathbf{U}}^\sumIndC \in \Real^{\dimInd_\sumIndC \times R_\sumIndC}$
\ENDFOR
\STATE Initialize $\mathbfcal{V}^0$ as $\mathbfcal{V}^0_{\sumIndA_1, \ldots, \sumIndA_n, \sumIndA} := ((\tilde{\mathbfcal{R}}_{\mPoint})_{\sumIndA_1, \ldots, \sumIndA_n}, (\phi_{\mPoint})_{\sumIndA})_{\mPoint}$ for any orthonormal basis  $\{\phi_{\mPoint}^{\sumIndA}\}_{\sumIndA=1}^\dimInd\subset \tangent_{\mPoint}\manifold$
\WHILE{not converged}
\STATE Compute $\nabla g (\mathbfcal{V}^{\iterInd})$ using \cref{eq:g-grad-symmetric} 
\STATE $\mathbfcal{V}^{\iterInd+1} := \mathbfcal{V}^{\iterInd} - \tau \nabla g (\mathbfcal{V}^{\iterInd})$
\STATE $\iterInd := \iterInd +1$
\ENDWHILE
\STATE $\mathbfcal{R}_{\mPoint} := \mathbfcal{V}^\iterInd \times_{n+1} \phi_{\mPoint} $
\RETURN $\mathbfcal{R}_{\mPoint} \in \tangent_{\mPoint} \manifold^{r_1 \times \ldots \times r_n}$ and $\{\mathbf{U}^\sumIndC\}_{\sumIndC=1}^n$ with $\mathbf{U}^\sumIndC \in \Real^{\dimInd_\sumIndC \times r_\sumIndC}$
\end{algorithmic}
\end{algorithm}

\begin{remark}
    We note that \cref{alg:geometry-corrected-HOSVD} does not necessarily converge because of the non-convexity due to the non-linearity of the function $g$. In practice, a small step-size will be necessary to promote numerical stability, but will deteriorate the run time. In other words, this method runs into the same issue as methods such as \cite{chakraborty2016efficient,lazar2017scale,said2007exact,sommer2014optimization}. That is, MC-tHOSVD can do well in terms of accounting for global geometry and even general applicability, but is very likely to be slow -- even with the closed-form gradient that was not used in aforementioned work.
\end{remark}


\section{Numerics.} 
\label{sec:manifold-valued-tensors-numerics}
In this section the proposed methods from \cref{sec:cc-lra} are tested through several numerical experiments with both 1D synthetic and 2D real data living in symmetric Riemannian manifolds. Our overall goal is to test the extent to which the CC-tHOSVD (\cref{alg:curvature-corrected-HOSVD}) is generally applicable, computationally feasible and global-geometry aware in realistic settings compared to baseline-like approaches such as tHOSVD and MC-tHOSVD (\cref{alg:tangent-space-HOSVD,alg:geometry-corrected-HOSVD}). In particular, we expect that CC-tHOSVD will be a compromise between tHOSVD and MC-tHOSVD in the sense that it will be similar -- and only slightly slower -- to tHOSVD in terms of run time and similar -- but slightly less optimal -- to MC-tHOSVD in terms of global-geometry awareness. All of the experiments are implemented in Julia 1.7.2 and run on a 2 GHz Quad-Core Intel Core i5 with 16GB RAM. We use \texttt{Manifolds.jl} \cite{axen2021manifolds} for all manifold mappings and use the export functionalities in \texttt{Manopt.jl} \cite{bergmann2022manopt} for generating the visualizations of the data. All reported time measurements have been made using the \texttt{@benchmark} macro from \texttt{BenchmarkTools.jl}. 

\paragraph{Overview of the experiments.}
To get some basic intuition for the algorithms' performance, we consider the effect of curvature on low rank approximation of synthetic 1D manifold-valued signals in \cref{sec:numerics-low-rank-1D}. In particular, we focus on the effects of positive curvature through signals on the 6-sphere ($\Sphere^6$) and on the effects of negative curvature through signals on the space of $3 \times 3$ symmetric positive definite matrices ($\mathcal{P}(3)$). Besides, getting an idea how the proposed methods perform on these spaces, our main goal in this experiment is to do a comparison in the performance between \cref{alg:tangent-space-HOSVD,alg:curvature-corrected-HOSVD,alg:geometry-corrected-HOSVD} in terms of numerical feasibility and global-geometry awareness -- as all of these methods are generally applicable to general symmetric manifolds.

The experiments in \cref{sec:numerics-low-rank-1D} then suggests that negative curvature is more problematic for low rank approximation than positive curvature \textit{if uncorrected for}, and that the proposed CC-tHOSVD outperforms the two other methods in terms of having both global-geometry awareness and speed. That is, this is contrary to tHOSVD, which only does well in terms of speed, and MC-tHOSVD, which only does well in terms of global-geometry awareness. So naturally, in \cref{sec:numerics-low-rank-2D} we will focus on \cref{alg:curvature-corrected-HOSVD} and move to real-world diffusion tensor MRI data, which is also $\mathcal{P}(3)$-valued. This time we focus on the limitations of CC-tHOSVD as this data set does not necessarily have low rank.



\paragraph{General error metrics.}
Throughout all of the experiments, given a (multi-linear) rank-$r$ approximation $\mTVector^r_{\mPoint}\in \tangent_{\mPoint}\manifold^{\dimInd_1 \times \ldots\times \dimInd_n}$ with $r \in \Natural^n$ and $n \in \{1,2\}$ obtained from one of the proposed methods, we consider the \emph{relative error}
\begin{equation}
    \varepsilon_{rel} (r) := \frac{d_{\manifold^{d_1 \times \cdots \times d_{n}}}(\Tensor, \exp_{\mPoint}(\mTVector^r_{\mPoint}))^2}{d_{\manifold^{d_1 \times \cdots \times d_{n}}}(\Tensor, \mPoint)^2} =  \frac{\sum_{i_1,\ldots, i_{n}} d_{\manifold}(\Tensor_{i_1,\ldots, i_{n}}, \exp_{\mPoint}((\mTVector^r_{\mPoint})_{i_1,\ldots, i_{n}}))^2}{\sum_{i_1,\ldots, i_{n}} d_{\manifold}(\Tensor_{i_1,\ldots, i_{n}}, \mPoint)^2},
    \label{eq:numerics-relative-error}
\end{equation}
as our main metric for measuring whether the tangent space-based schemes actually give a low error on the manifold. A useful tool to complement interpreting the quality of approximation through relative error is the lower bound \cref{eq:cor-lower-bound} from \cref{cor:approximation-error-lower-bound} as already noted in \cref{rem:using-lower-bound}. In 1D the truncated tangent space (HO)SVD is optimal by Eckart-Young. In other words, we have access to the lower bound in the corollary. However, in 2D the truncated tangent space HOSVD will in general not be the minimizer, 
as there is no Eckart-Young result for the HOSVD. Nevertheless, tHOSVD can still give insight as a pseudo-lower bound in the latter case. It should also be noted that \cref{cor:approximation-error-lower-bound} lower bounds the curvature corrected approximation error \cref{eq:inf-Sr-linearized-tangent-approximation} and not necessarily the global approximation error \cref{eq:global-approximation-error-minimization-low-rank}, on which the relative error is based. The lower bound will hold for the relative error as the tangent space error vanishes, since that results in vanishing discrepancy. So to avoid confusion the lower bound from \cref{cor:approximation-error-lower-bound} is referred to as the zero-discrepancy (pseudo) lower bound or \emph{zero-$\delta$ (pseudo) lower bound} throughout the rest of this section.

Additional to the relative error and only for CC-tHOSVD, we consider the \emph{relative discrepancy}
\begin{equation}
    \delta_{rel} (r) := \frac{\sum_{\sumIndB_1, \ldots, \sumIndB_n, \sumIndB=1}^{\dimInd_1, \ldots, \dimInd_n, \dimInd} \beta(\kappa_{\sumIndB_1, \ldots, \sumIndB_n, \sumIndB})^2 \Bigl( \mTVector^r_{\mPoint} - \log_{\mPoint}\Tensor,  \Theta^{\sumIndB_1, \ldots, \sumIndB_n, \sumIndB}_{\mPoint} \Bigr)_{\mPoint}^2- d_{\manifold^{d_1 \times \cdots \times d_{n}}}(\Tensor, \exp_{\mPoint}(\mTVector^r_{\mPoint}))^2 }{\|\mTVector^r_{\mPoint} - \log_{\mPoint}\Tensor\|^3_\mPoint },
\end{equation}
as a basic check, as we expect $\delta_{rel}(r) \in \mathcal{O}(1)$ by our main result \cref{eq:general-error-symmetric spaces}.

\subsection{Synthetic low-rank 1D signals.} 
\label{sec:numerics-low-rank-1D}

We consider low rank 1D signals on symmetric Riemannian manifolds with both positive and negative curvature. As these data are tensors of order 1, we will more naturally refer to \cref{alg:tangent-space-HOSVD,alg:curvature-corrected-HOSVD,alg:geometry-corrected-HOSVD} as tSVD, CC-tSVD and MC-tSVD respectively, i.e., we drop the \enquote{higher-order} (HO) from their names.



\subsubsection{Spherical data.}
Starting with spherical signals -- which is a symmetric manifold under the standard Euclidean pull-back metric --, consider the order-1 manifold-valued tensor $\Tensor'\in (\Sphere^6)^{100}$ defined entry-wise as
\begin{equation}
    \Tensor'_\sumIndA:= (\cos(\varphi_\sumIndA), \sin(\varphi_\sumIndA), 0, 0, 0, 0, 0)^\top, \quad \varphi_\sumIndA := 2 \pi \frac{\sumIndA-1}{100}, \; \sumIndA=1, \ldots, 100,
\end{equation}
and $\Tensor\in (\Sphere^6)^{100}$ defined entry-wise as
\begin{equation}
    \Tensor_\sumIndA:= \exp_{ \Tensor'_\sumIndA} (\eta_{ \Tensor'_\sumIndA}^{\sumIndA}), \quad \text{where } \eta_{ \Tensor'_\sumIndA}^{\sumIndA}\in \tangent_{\Tensor'_\sumIndA} \Sphere^6 \text{ is a centered Gaussian random tangent vector with variance $0.05$.}
    \label{eq:1D-S6-data-set}
\end{equation}
The data set can be visualized through projecting its first three coordinates onto $\Sphere^2$ as done in \cref{fig:1d-S6-data}.

\begin{figure}[h!]
    \centering
    \includegraphics[width=0.4\linewidth]{"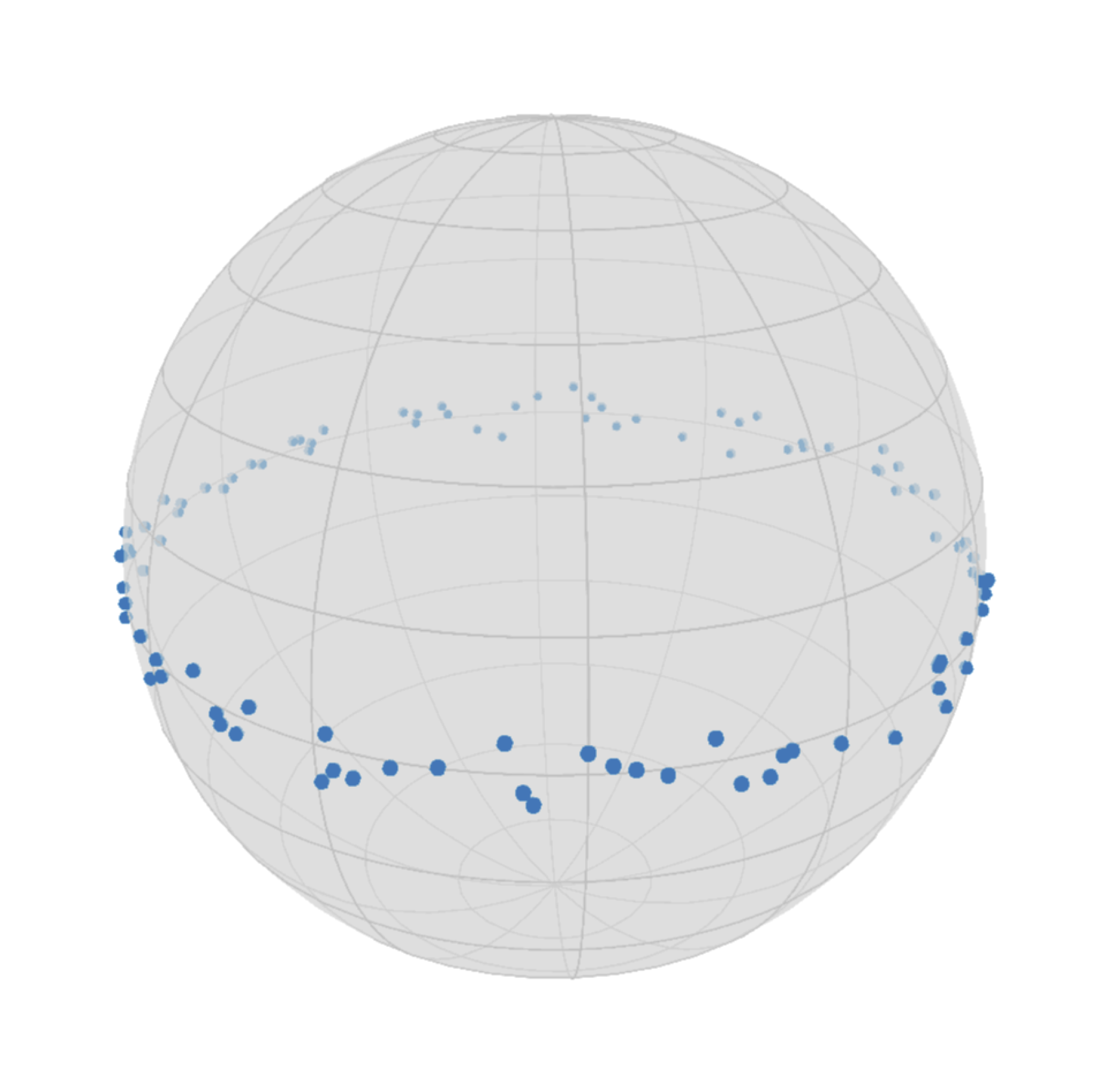"}
    \caption{The first three coordinates of the $\Sphere^6$-valued 1D synthetic data set defined in \cref{eq:1D-S6-data-set} projected on $\Sphere^2$.}
    \label{fig:1d-S6-data}
\end{figure}

Even though all mappings needed for executing \cref{alg:tangent-space-HOSVD,alg:curvature-corrected-HOSVD,alg:geometry-corrected-HOSVD} are provided in \texttt{Manifolds.jl} \cite{axen2021manifolds}, it is worth highlighting that we can construct the diagonalizing orthonormal basis at some point $\mPoint\in \Sphere^6$ \cite{bacak2016second} of the form \cref{eq:assumption-digonalizing-basis} by choosing $(\theta^{\sumIndB_1}_{\mPoint})_{1} := \frac{\log_{\mPoint}\Tensor_{\sumIndB_1}}{\|\log_{\mPoint}\Tensor_{\sumIndB_1}\|_{\mPoint}}\in \tangent_\mPoint  \Sphere^6$ and complete this to an orthonormal basis $\{(\theta^{\sumIndB_1}_{\mPoint})_{1}, \ldots, (\theta^{\sumIndB_1}_{\mPoint})_{6}\}\subset \tangent_\mPoint \mathbb{S}^6$.
The corresponding eigenvalues are then given by
\begin{equation}
    \kappa_{\sumIndB_1, \sumIndB} = \left \{\begin{matrix}
        0, & \text{if $\sumIndB=1$}, \\
        d_{\mathbb{S}^6}(\mPoint, \Tensor_{\sumIndB_1})^2, & \text{if $\sumIndB\neq 1$}.
        \end{matrix}\right.
        \label{eq:Sd-kappas}
\end{equation}



Next, in the Euclidean case a natural choice for the base point $\mPoint$ would be the mean of the data set. This can be generalized in the manifold setting by passing to the Riemannian barycentre, but as noted in work such as \cite{jung2012analysis} this can be problematic in the case of spherical data. Indeed, for the data set defined in \cref{eq:1D-S6-data-set} the barycentre will be close to a pole of $\Sphere^6$ and the data will look 2-dimensional rather than 1-dimensional as \cref{fig:1d-S6-data} suggests. For getting intuition for how the proposed methods behave in both problematic and non-problematic choices of base point we will take a closer look at this behavior
by considering two base points. That is, we consider
\begin{equation}
    \mPoint^1 := \hat{\mPoint}^1 \in \argmin_{\mPoint\in \Sphere^6} \sum_{\sumIndA=1}^{100} d_{\mathbb{S}^6}(\mPoint, \Tensor_{\sumIndA})^2, \quad \text{and} \quad \mPoint^2 := \hat{\mPoint}^2 \in \argmin_{\mPoint\in \Tensor} \sum_{\sumIndA=1}^{100} d_{\mathbb{S}^6}(\mPoint, \Tensor_{\sumIndA})^2.
    \label{eq:p1-p2}
\end{equation}

So now given the manifold-valued tensor $\Tensor\in (\Sphere^6)^{100}$ and base points $\mPoint^1,\mPoint^2\in \Sphere^6$, we compare \cref{alg:tangent-space-HOSVD,alg:curvature-corrected-HOSVD,alg:geometry-corrected-HOSVD} in terms of global-geometry awareness and computational feasibility. For the former we will consider the relative error \cref{eq:numerics-relative-error} for all ranks except full rank, i.e., $r\in \{1, \ldots, 5\}$, and compare these to the zero-$\delta$ lower bound. For the latter, we consider the averaged run times for all ranks using the \texttt{@benchmark} macro as mentioned above. \Cref{alg:tangent-space-HOSVD,alg:curvature-corrected-HOSVD} do not need extra parameters, but \cref{alg:geometry-corrected-HOSVD} needs a stopping criterion and a step size. For a fair comparison, we pick the largest step size  -- up to a factor 2 -- under which the iterates converge. This gives $\tau_1 = 2^{-6} = 1/64$ and $\tau_2 = 2^{-4} = 1/16$. We let the algorithm run until the mild relative gradient condition $\|\nabla g (\mathbfcal{V}^{\iterInd})\|_2/\|\nabla g (\mathbfcal{V}^{0})\|_2< 10^{-2}$ is satisfied. Finally and additionally, to check whether the curvature corrected approximation error really is a good proxy for the global approximation error on $\Sphere^6$ we will also consider the relative discrepancy.

The relative error and the relative discrepancy results for low rank approximation at $\mPoint_1$ and $\mPoint_2$ are shown in \cref{fig:p1-error-and-discr,fig:p2-error-and-discr}. The averaged run times are provided in \cref{tab:p1-run-times,tab:p2-run-times}
We note that MC-tSVD converges and we refer the reader for more details to the convergence plots in \cref{fig:p1-convergence,fig:p2-convergence} in \cref{app:extra-results}.

\paragraph{Low rank approximation at $\mPoint_1\in \Sphere^6$.}
The relative error results in \cref{fig:p1-errors} are in line with our expectations, i.e., that both CC-tSVD and MC-tSVD pick up on the global geometry better than tSVD, although it has to be emphasized that the difference is minor and is only really visible for rank 1 approximation\footnote{Indeed, for the other ranks the error vanishes and the logarithmic-scaled plot in \cref{fig:p1-log-errors} in \cref{app:extra-results} suggests that all methods perform about equally well.}. This was indeed to be expected, as the effect of positive curvature on the global approximation error is less detrimental than the effect of negative curvature according to \cref{eq:general-error-symmetric spaces}. The figure also confirms that even though the data has rank 1 for a point on the equator, it looks like a rank-2 data set from $\mPoint_1$. More interestingly, the results do not come close to the zero-$\delta$ lower bound, which very likely because of the lower bound being overly optimistic. Indeed, having already one point far away from the base point will give a large $\kappa$ -- as can be seen in \cref{eq:Sd-kappas} --, and subsequently will give an unrealistically low $\beta(\kappa_{\max})$-value in the lower bound \cref{eq:cor-lower-bound}. Moving on, the relative discrepancy in \cref{fig:p1-discr} makes an encouraging case for the accuracy of \cref{eq:general-error-symmetric spaces}. Indeed, we know that $\delta_{rel}(r) \in \mathcal{O}(1)$. Now, the figure suggests that the constant is very small. In other words, the curvature corrected approximation error is global-geometry aware for positive curvature symmetric Riemannian manifolds.

Finally, \cref{tab:p1-run-times} tells us that the proposed CC-tSVD algorithm is not only global-geometry aware, but also computationally efficient. Indeed, the CC-tSVD run time is similar to that of tSVD, whereas only one iteration of MC-tSVD already takes up more time. It should be noted here that all three algorithms are implemented in a similar way, e.g., using the same type of vectorization and memory allocation. So these discrepancies are very likely caused by the gradients in \cref{alg:geometry-corrected-HOSVD} being more expensive to compute, possibly due to the extra 100 logarithmic and exponential map evaluations and the additional 600 parallel transport evaluations per gradient call.

\begin{figure}[h!]
    \centering
    \begin{subfigure}{0.48\linewidth}
    \includegraphics[width=\linewidth]{"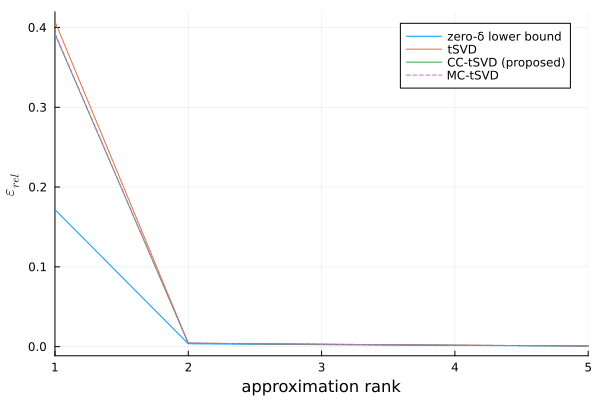"}
    \caption{}
    \label{fig:p1-errors}
    \end{subfigure}
    \hfill
    \begin{subfigure}{0.48\linewidth}
    \includegraphics[width=\linewidth]{"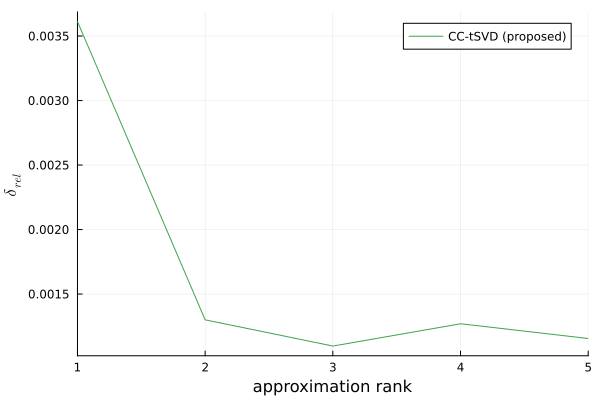"}
    \caption{}
    \label{fig:p1-discr}
    \end{subfigure}
    \caption{The progressions of the relative error (a) and the relative discrepancy (b) of \cref{alg:tangent-space-HOSVD,alg:curvature-corrected-HOSVD,alg:geometry-corrected-HOSVD} for an increasing approximation rank of 1D data on $\Sphere^6$, linearized at $\mPoint_1$ from \cref{eq:p1-p2}. (a) suggests that positive curvature does not deteriorate the tSVD approximation, but shows that CC-tSVD and MC-tSVD attain lower errors through picking up on global geometry. (b) shows that discrepancy between the curvature corrected approximation error and the exact global approximation error is very tight, supporting the accuracy of \cref{eq:general-error-symmetric spaces}.}
    \label{fig:p1-error-and-discr}
\end{figure}

\begin{table}[h!]
    \centering
    \caption{The average run times in seconds of \cref{alg:tangent-space-HOSVD,alg:curvature-corrected-HOSVD,alg:geometry-corrected-HOSVD} per approximation rank of 1D data on $\Sphere^6$ linearized at $\mPoint_1$ from \cref{eq:p1-p2}. tSVD and CC-tSVD have similar runtimes, but MC-tSVD is significantly slower.}
    \label{tab:p1-run-times}
    \begin{tabular}{lccccc}
            \toprule
            Method & $\operatorname{rank} = 1$  & $\operatorname{rank} = 2$  & $\operatorname{rank} = 3$  & $\operatorname{rank} = 4$  & $\operatorname{rank} = 5$   \\ 
            \midrule
            tSVD & $0.0008435$ & $0.0007377$ & $0.0008345$ & $0.001223$ & $0.00087$ \\ 
            CC-tSVD (proposed) & $0.00193$ & $0.002363$ & $0.002367$ & $0.002722$ & $0.002932$ \\ 
            MC-tSVD (1 iteration) & $0.00629$ & $0.006916$ & $0.00837$ & $0.00856$ & $0.008385$\\ 
            MC-tSVD & $0.3813$ & $0.3687$ & $0.3984$ & $0.4485$ & $0.508$\\ 
            \bottomrule
    \end{tabular}
\end{table}

\paragraph{Low rank approximation at $\mPoint_2\in \Sphere^6$.}
Most of the discussion on low rank approximation at $\mPoint_1$ in terms of global-geometry awareness applies to $\mPoint_2$. That is, even though the positive curvature is again not strongly deteriorating the naive tSVD approximation, both CC-tSVD and MC-tSVD do pick up on the global geometry slightly better than the naive method for the rank-1 approximation (\cref{fig:p2-errors}). For the remaining ranks CC-tHOSVD performs slightly worse, i.e., even compared to tSVD\footnote{The logarithmic-scaled plot in \cref{fig:p2-log-errors} in \cref{app:extra-results} makes this easier to see.}. The latter observation is most likely due to an increase in the -- still very small -- relative discrepancy (\cref{fig:p2-discr}) compared to the $\mPoint_1$ case (\cref{fig:p1-discr}). This increase makes it more challenging to get a very accurate correction for already almost-zero loss. However, in practice these rank-2 through rank-5 cases are less important than the rank-1 case. So overall, the curvature corrected approximation error does well once again in terms of global-geometry awareness for positive curvature symmetric Riemannian manifolds.

Regarding computational feasibility, we also observe the same as before. That is, \cref{tab:p2-run-times} shows that tSVD and CC-tSVD have similar run times and MC-tSVD is significantly slower.


\begin{figure}[h!]
    \centering
    \begin{subfigure}{0.48\linewidth}
    \includegraphics[width=\linewidth]{"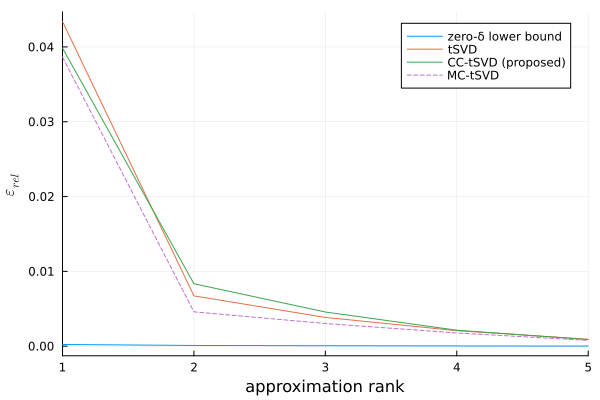"}
    \caption{}
    \label{fig:p2-errors}
    \end{subfigure}
    \hfill
    \begin{subfigure}{0.48\linewidth}
    \includegraphics[width=\linewidth]{"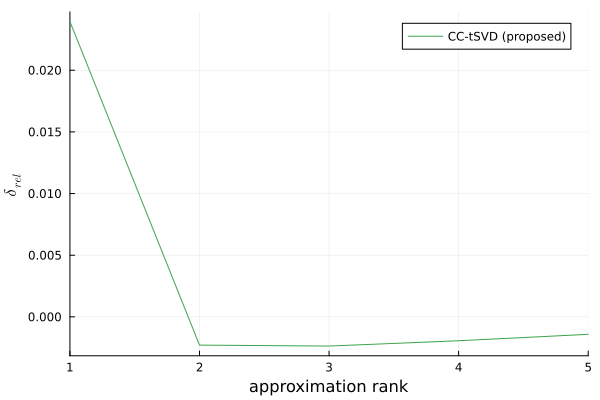"}
    \caption{}
    \label{fig:p2-discr}
    \end{subfigure}
    \caption{The progressions of the relative error (a) and the relative discrepancy (b) of \cref{alg:tangent-space-HOSVD,alg:curvature-corrected-HOSVD,alg:geometry-corrected-HOSVD} for an increasing approximation rank of 1D data on $\Sphere^6$, linearized at $\mPoint_2$ from \cref{eq:p1-p2}. (a) shows that CC-tSVD and MC-tSVD attain slightly lower -- but similar -- errors through picking up on global geometry for rank 1. For higher ranks both tSVD and MC-tSVD perform slightly better. (b) shows that discrepancy between the curvature corrected approximation error and the exact global approximation error is less tight than in \cref{fig:p1-error-and-discr}, which is the likely cause of the difference in performance for the higher rank cases as a more accurate approximation is likely needed.}
    \label{fig:p2-error-and-discr}
\end{figure}

\begin{table}[h!]
    \centering
    \caption{The average run times in seconds of \cref{alg:tangent-space-HOSVD,alg:curvature-corrected-HOSVD,alg:geometry-corrected-HOSVD} per approximation rank of 1D data on $\Sphere^6$ linearized at $\mPoint_2$ from \cref{eq:p1-p2}. Similarly as in \cref{tab:p2-run-times} tSVD and CC-tSVD are much faster than MC-tSVD.}
    \label{tab:p2-run-times}
    \begin{tabular}{lccccc}
            \toprule
            Method & $\operatorname{rank} = 1$  & $\operatorname{rank} = 2$  & $\operatorname{rank} = 3$  & $\operatorname{rank} = 4$  & $\operatorname{rank} = 5$   \\ 
            \midrule
            tSVD & $0.0008183$ & $0.0008283$ & $0.0009403$ & $0.0008883$ & $0.000927$\\ 
            CC-tSVD (proposed) & $0.001951$ & $0.002422$ & $0.002525$ & $0.003294$ & $0.003426$\\ 
            MC-tSVD (1 iteration) & $0.00653$ & $0.00752$ & $0.00781$ & $0.008286$ & $0.00883$\\ 
            MC-tSVD & $0.2112$ & $1.524$ & $1.911$ & $2.178$ & $1.934$\\ 
            \bottomrule
    \end{tabular}
\end{table}

\paragraph{Low rank approximation on $(\Sphere^6)^{100}$.}

Overall, the above results strongly support our claim that CC-tSVD is a good compromise between global-geometry awareness and computational feasibility in the case of positive curvature. That is, compared to tSVD, which is fast but not global-geometry aware, and compared to MC-tSVD, which is global-geometry aware but slow. It should be emphasized once more that the relative errors are very comparable, which is to be expected as curvature effects due to positive curvature are less of a problem than those due to negative curvature according to \cref{eq:general-error-symmetric spaces}. Next, we will get a better idea to what extent the proposed method can account for the more problematic case of negative curvature.


\subsubsection{Symmetric positive definite matrix data.}
\label{sec:P3-1D-tensor}
Next, we move on to $\mathcal{P}(3)$ -- which is a symmetric manifold under the bi-invariant metric -- and consider the order-1 manifold-valued tensor $\Tensor'\in (\mathcal{P}(3))^{100}$ defined entry-wise as
\begin{equation}
    \Tensor'_\sumIndA:= \exp_{I_3} (\tau_\sumIndA \begin{pmatrix}
0 & 0 & 0 \\
0 & 1 & 0 \\
0 & 0 & 0 \\
\end{pmatrix}), \quad \tau_\sumIndA\sim \mathcal{N}(0,2),\; \sumIndA=1, \ldots, 100,
\end{equation}
and $\Tensor\in \mathcal{P}(3)^{100}$ defined entry-wise as
\begin{equation}
    \Tensor_\sumIndA:= \exp_{\Tensor'_\sumIndA} (\eta_{ \Tensor'_\sumIndA}^{\sumIndA}), \quad \text{where } \eta_{ \Tensor'_\sumIndA}^{\sumIndA}\in \tangent_{\Tensor'_\sumIndA} \mathcal{P}(3) \text{ is a centered Gaussian random tangent vector with variance $0.05$.}
    \label{eq:1D-P3-data-set}
\end{equation}
A subset of this data set is visualized in \cref{fig:1d-P3-data}, from which it is clear that the data set has rank 1.

\begin{figure}[h!]
    \centering
    \includegraphics[width=0.4\linewidth]{"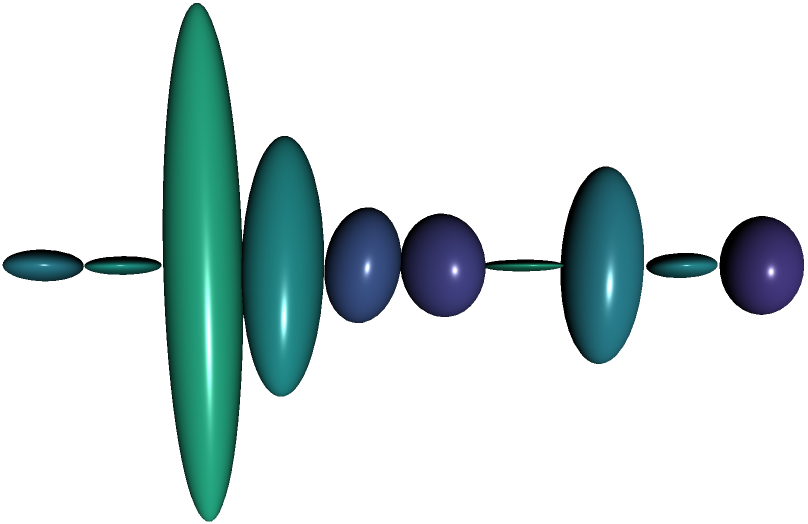"}
    \caption{A visualization of the 1D $\mathcal{P}(3)$-valued data set defined in \cref{eq:1D-P3-data-set}. }
    \label{fig:1d-P3-data}
\end{figure}

Again, for the sake of completeness, it is worth highlighting that we can construct the diagonalizing orthonormal basis at some point $\mPoint\in \mathcal{P}(3)$ \cite{bacak2016second} of the form \cref{eq:assumption-digonalizing-basis} by first considering the eigendecomposition of $\frac{\log_{\mPoint}\Tensor_{\sumIndB_1}}{\|\log_{\mPoint}\Tensor_{\sumIndB_1}\|_{\mPoint}}\in \tangent_\mPoint  \mathcal{P}(3)$, which is a symmetric $3 \times 3$ matrix, i.e.,
\begin{equation}
    \frac{\log_{\mPoint}\Tensor_{\sumIndB_1}}{\|\log_{\mPoint}\Tensor_{\sumIndB_1}\|_{\mPoint}} = \sum_\sumIndC^3  \lambda_\sumIndC^{\sumIndB_1} \mathbf{v}_\sumIndC^{\sumIndB_1} (\mathbf{v}_\sumIndC^{\sumIndB_1})^\top.
\end{equation}
Next, defining the index set
\begin{equation}
    \mathcal{J} := \{ \sumIndB = (\sumIndC, \sumIndD)\: \mid \: \sumIndC = 1,\dots , 3, \; \sumIndD = \sumIndC, \ldots, 3\},
\end{equation}
we obtain an orthonormal basis $\{(\theta^{\sumIndB_1}_{\mPoint})_{(1,1)}, (\theta^{\sumIndB_1}_{\mPoint})_{(1,2)}, (\theta^{\sumIndB_1}_{\mPoint})_{(1,3)}, (\theta^{\sumIndB_1}_{\mPoint})_{(2,2)}, (\theta^{\sumIndB_1}_{\mPoint})_{(2,3)}, (\theta^{\sumIndB_1}_{\mPoint})_{(3,3)}\}\subset \tangent_\mPoint \mathcal{P}(3)$ given by
\begin{equation}
    (\theta^{\sumIndB_1}_{\mPoint})_{\sumIndB} := \left\{\begin{matrix}
\frac{1}{2}(\mathbf{v}^{\sumIndB_1}_\sumIndC (\mathbf{v}^{\sumIndB_1}_\sumIndD)^\top +\mathbf{v}_\sumIndD^{\sumIndB_1} (\mathbf{v}_\sumIndC^{\sumIndB_1})^\top), & \sumIndB =(\sumIndC, \sumIndD)\in \mathcal{J} & \text{if $\sumIndC=\sumIndD$}, \\
\frac{1}{\sqrt{2}}(\mathbf{v}_\sumIndC^{\sumIndB_1} (\mathbf{v}_\sumIndD^{\sumIndB_1})^\top + \mathbf{v}_\sumIndD^{\sumIndB_1} (\mathbf{v}_\sumIndC^{\sumIndB_1})^\top), & \sumIndB =(\sumIndC, \sumIndD)\in \mathcal{J} & \text{if $\sumIndC\neq \sumIndD$},
\label{eq:Pr-basis-diagonalizing}
\end{matrix}\right. .
\end{equation}
The corresponding eigenvalues are then given by
\begin{equation}
    \kappa_{\sumIndB_1, \sumIndB} = -\frac{1}{4} (\lambda_\sumIndC^{\sumIndB_1} - \lambda_\sumIndD^{\sumIndB_1})^2 d_{\mathcal{P}(3)}(\mPoint, \Tensor_{\sumIndB_1})^2, \quad \sumIndB =(\sumIndC, \sumIndD)\in \mathcal{J}.
\end{equation}
Unlike the spherical case, we do not expect problematic behavior from picking the base point as the barycentre of the data. That is, we consider
\begin{equation}
    \mPoint := \argmin_{\mPoint\in\mathcal{P}(3)} \sum_{\sumIndA=1}^{100} d_{\mathcal{P}(3)}(\mPoint, \Tensor_{\sumIndA})^2.
    \label{eq:p-P3}
\end{equation}
Note that we do indeed get a unique minimizer due to the fact that the barycentre problem is strongly geodesically convex on the Hadamard manifold $\mathcal{P}(3)$.

We will run the same experimental setup as for spherical data. That is, we use the same stopping criterion for  \cref{alg:geometry-corrected-HOSVD} as before and again  use the largest step size -- up to a factor -- that does not result in a diverging solver, which is $\tau = 2^{-16}$.

The relative error and the relative discrepancy results for low rank approximation at $\mPoint$ are shown in \cref{fig:p-error-and-discrl-P3}. The averaged run times are provided in \cref{tab:p-run-times-P3}. We note that MC-tSVD converges and we refer the reader for more details to the convergence plots in \cref{fig:p-convergence-P3} in \cref{app:extra-results}.

\paragraph{Low rank approximation at $\mPoint\in \mathcal{P}(3)$.}
For the data set in \cref{eq:1D-P3-data-set} the effect of curvature is more clearly pronounced than for the positive curvature case. Indeed, there is a large difference between the relative errors of tSVD and those of CC-tSVD and MC-tSVD (\cref{fig:p-errors-P3}). Moreover, both CC-tSVD and MC-tSVD almost reach the zero-$\delta$ lower bound\footnote{This is visible more clearly in the logarithmic-scaled plot in \cref{fig:p-log-errors-P3}.}, which makes a strong case for the global-geometry awareness of CC-tSVD for symmetric Riemannian manifolds with negative curvature. This also suggests that the lower bound is actually reasonably tight in the negative curvature case, contrary to the positive curvature case. Additional evidence of the global-geometry awareness of CC-tSVD follows from \cref{fig:p-discr-P3}, which also suggests that the constant we expect from $\delta_{rel}(r) \in \mathcal{O}(1)$ is once again very small. Summarized, we conclude that the curvature corrected approximation error and CC-tSVD are global-geometry aware for negative curvature as well. 

Regarding computational feasibility, we make similar observations as for the positive curvature case. That is, \cref{tab:p-run-times-P3} shows that tSVD and CC-tSVD have similar run times and MC-tSVD is on average more than an order of magnitude slower.

\begin{figure}[h!]
    \centering
    \begin{subfigure}{0.48\linewidth}
    \includegraphics[width=\linewidth]{"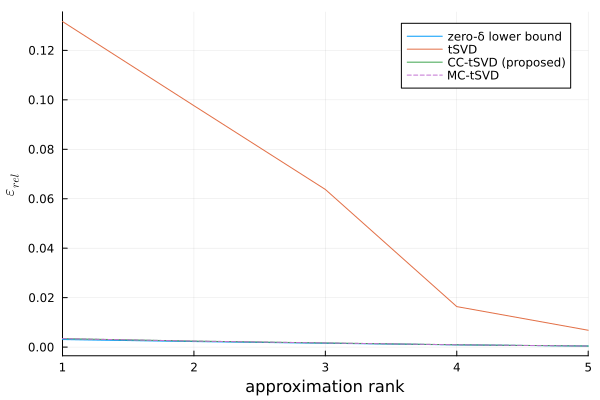"}
    \caption{}
    \label{fig:p-errors-P3}
    \end{subfigure}
    \hfill
    \begin{subfigure}{0.48\linewidth}
    \includegraphics[width=\linewidth]{"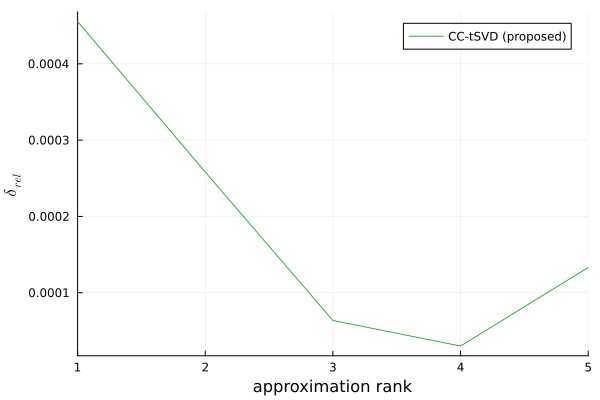"}
    \caption{}
    \label{fig:p-discr-P3}
    \end{subfigure}
    \caption{The progressions of the relative error (a) and the relative discrepancy (b) of \cref{alg:tangent-space-HOSVD,alg:curvature-corrected-HOSVD,alg:geometry-corrected-HOSVD} for an increasing approximation rank of 1D data on $\mathcal{P}(3)$, linearized at $\mPoint$ from \cref{eq:p-P3}. (a) shows that CC-tSVD and MC-tSVD attain significantly lower errors through picking up on global geometry for all ranks. (b)  shows that discrepancy between the curvature corrected approximation error and the exact global approximation error is very tight, supporting again the accuracy of \cref{eq:general-error-symmetric spaces}.}
    \label{fig:p-error-and-discrl-P3}
\end{figure}

\begin{table}[h!]
    \centering
    \caption{The average run times in seconds of \cref{alg:tangent-space-HOSVD,alg:curvature-corrected-HOSVD,alg:geometry-corrected-HOSVD} per approximation rank of 1D data on $\mathcal{P}(3)$ linearized at $\mPoint$ from \cref{eq:p-P3}. Similarly to \cref{tab:p1-run-times,tab:p2-run-times} tSVD and CC-tSVD are much faster than MC-tSVD.}
    \label{tab:p-run-times-P3}
    \begin{tabular}{lccccc}
            \toprule
            Method & $\operatorname{rank} = 1$  & $\operatorname{rank} = 2$  & $\operatorname{rank} = 3$  & $\operatorname{rank} = 4$  & $\operatorname{rank} = 5$   \\ 
            \midrule
            tSVD & $0.01811$ & $0.01816$ & $0.0177$ & $0.02174$ & $0.01868$\\ 
            CC-tSVD (proposed) & $0.03735$ & $0.0394$ & $0.04367$ & $0.05252$ & $0.05743$\\ 
            MC-tSVD (1 iteration) & $0.1534$ & $0.1726$ & $0.1981$ & $0.2029$ & $0.2186$\\ 
            MC-tSVD & $1.091$ & $0.933$ & $0.5923$ & $0.3997$ & $0.424$\\ 
            \bottomrule
    \end{tabular}
\end{table}

\paragraph{Low rank approximation on $\mathcal{P}(3)^{100}$.}

Overall, the above results also strongly support our claim that CC-tSVD is a good compromise between global-geometry awareness and computational feasibility in the case of positive curvature. It should be emphasized this time that the relative errors between the naive tSVD and CC-tSVD is very significant due to the negative curvature, but that the difference between MC-tSVD and CC-tSVD is very slight.

\subsection{Real 2D DT-MRI data.} 
\label{sec:numerics-low-rank-2D}
An important application of $\mathcal{P}(3)$-valued data processing is in \emph{diffusion tensor magnetic resonance imaging} (DT-MRI). In the following, we apply the proposed CC-tHOSVD on a DT-MRI data set\footnote{Follow the tutorial at \href{http://camino.cs.ucl.ac.uk/index.php?n=Tutorials.DTI}{http://camino.cs.ucl.ac.uk/index.php?n=Tutorials.DTI}.} of a human head from the Camino project \cite{cook2005camino} and compare performance with the naive tHOSVD.

The Camino data set gives us a manifold-valued tensor in $\Tensor' \in \mathcal{P}(3)^{112\times 112\times 50}$ of order 3. From this manifold-valued tensor, we construct two manifold-valued tensors of order 2. That is, we construct the manifold-valued tensors $\Tensor \in \mathcal{P}(3)^{40\times 40}$ and $\Tensor^{patch} \in \mathcal{P}(3)^{10\times 10}$ defined as
\begin{equation}
    \Tensor := \Tensor'_{38:77,36:75,15}, \quad \Tensor^{patch} := \Tensor'_{38:47,36:45,15} = \Tensor_{1:10,1:10}.
    \label{eq:2D-data-DTMRI}
\end{equation}
Both data sets are visualized in \cref{fig:camino-data}. 

\begin{figure}[h!]
    \centering
    \begin{subfigure}{0.48\linewidth}
    \centering
    \includegraphics[width=0.8\linewidth]{"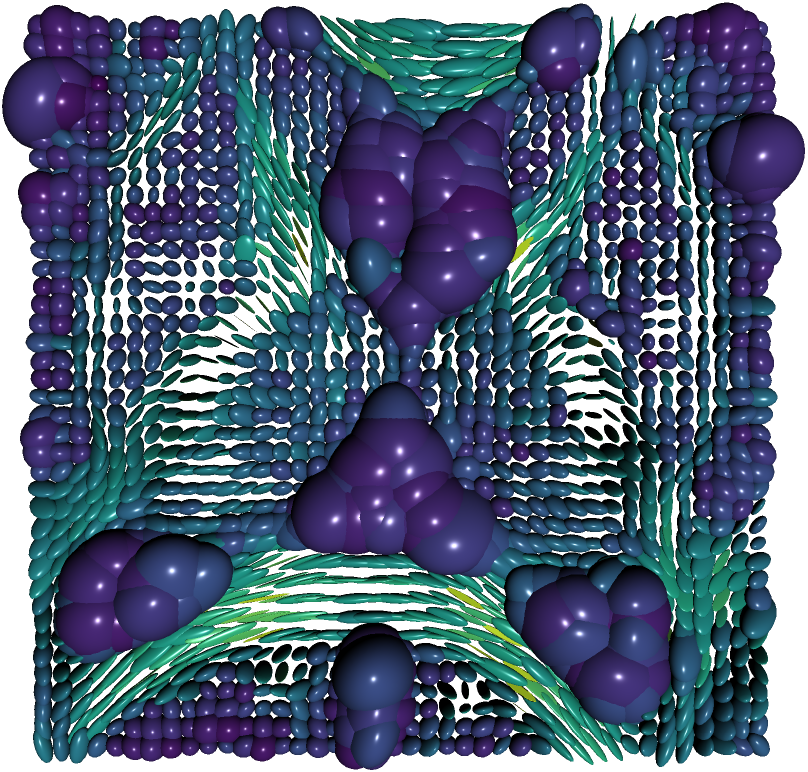"}
    \caption{}
    \label{fig:camino-data-full}
    \end{subfigure}
    \hfill
    \begin{subfigure}{0.48\linewidth}
    \centering
    \includegraphics[width=0.8\linewidth]{"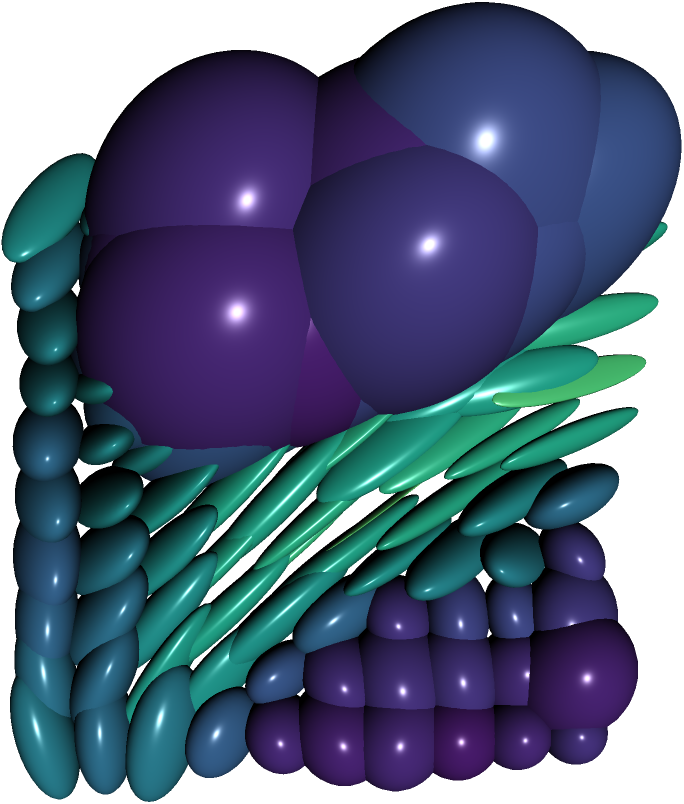"}
    \caption{}
    \label{fig:camino-data-patch}
    \end{subfigure}
    \caption{A visualization of the 2D $\mathcal{P}(3)$-valued slices $\Tensor \in \mathcal{P}(3)^{40\times 40}$ and $\Tensor^{patch} \in \mathcal{P}(3)^{10\times 10}$ of DT-MRI data from the Camino project \cite{cook2005camino}. Note that (b) is the lower left $10 \times 10$ patch of (a).}
    \label{fig:camino-data}
\end{figure}

The curvature tensor diagonalizing orthonormal basis at some point $\mPoint\in \mathcal{P}(3)$ of the form \cref{eq:assumption-digonalizing-basis} can be constructed analogously to the construction in \cref{sec:P3-1D-tensor}. Also analogously to above, we pick each data's respective barycentre as the base point for the low multi-linear rank approximation.


For exploring the behavior of the relative error and the relative discrepancy for approximating $\Tensor$ we compute a multi-linear rank-$(r,r)$ approximation using \cref{alg:tangent-space-HOSVD,alg:curvature-corrected-HOSVD} for $r$ ranging from 1 to 40 with intervals of 3. For the run time comparison, we take a coarser partition and use intervals of 7. Subsequently, for $\Tensor^{patch}$ we compute a multi-linear rank-$(r,r)$ approximation for all $r$ from 1 to 9, i.e., except full rank. For the run time comparison, only the odd ranks are benchmarked.

The relative error and the relative discrepancy results for low multi-linear rank-$(r,r)$ approximation of $\Tensor$ and $\Tensor^{patch}$ are shown in \cref{fig:p-error-and-discr-P3-Camino,fig:p-error-and-discr-P3-Camino-patch} respectively. The averaged run times are provided in \cref{tab:p-run-times-P3-Camino,tab:p-run-times-P3-Camino-patch}.


\paragraph{Low multi-linear rank approximation on $\mathcal{P}(3)^{40\times 40}$.}

The data set in \cref{fig:camino-data-full} is clearly testing the limits of the CC-tHOSVD now that the data set does not have a low multi-linear rank, unlike the 1D examples from above. In particular, as a large relative error at multi-linear rank $(1,1)$ goes hand in hand with large tangent space error, the discrepancy between curvature approximation error and the global approximation error in \cref{eq:general-error-symmetric spaces} is large. This leaves CC-tHOSVD unable to take the global geometry into account and results in a worse relative error than that of tHOSVD (\Cref{fig:p-errors-P3-Camino}). However, for larger approximation rank, the error naturally decreases and CC-tHOSVD once again picks up on the global geometry and reaches the zero-$\delta$ lower bound\footnote{The logarithmic-scaled plot in \cref{fig:p-log-errors-P3-Camino} in \cref{app:extra-results} makes this easier to see.}. The latter claim is also supported by the relative discrepancy in \cref{fig:p-discr-P3-Camino} that decreases to a small value, indicating that the constant we expect from $\delta_{rel}(r) \in \mathcal{O}(1)$ is once again very small.

Regarding computational feasibility, we observe that the run time increases quadratically in the approximation rank (\cref{tab:p-run-times-P3-Camino}), or equivalently linearly in the core tensor size. Quadratic scaling indicates that the bottleneck is the construction of the real-valued tensor $\mathbfcal{B}$ from \cref{eq:B-tensor} and not computing $\mathbfcal{A}$ from \cref{eq:A-tensor} -- which would scale with a fourth power in the approximation rank or quadratically in the core tensor size -- or solving the linear system -- which would scale with a sixth power in rank and cubically in the core tensor size. In other words, this would not be resolved by passing to MC-tHOSVD, which is -- similarly to the 1D case -- even slower due to the additional logarithmic map, exponential map and parallel transport calls\footnote{
The sheer computational infeasibility was the main reason to exclude MC-tHOSVD from this experiment.}. 

So overall, CC-tHOSVD remains to be -- even in the current implementation -- a good compromise between global-geometry awareness and computational feasibility.





\begin{figure}[h!]
    \centering
    \begin{subfigure}{0.48\linewidth}
    \includegraphics[width=\linewidth]{"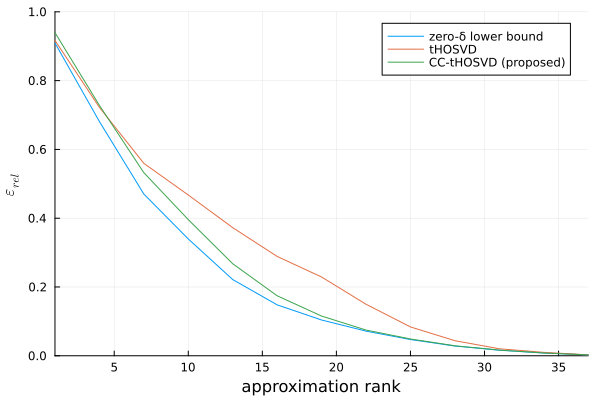"}
    \caption{}
    \label{fig:p-errors-P3-Camino}
    \end{subfigure}
    \hfill
    \begin{subfigure}{0.48\linewidth}
    \includegraphics[width=\linewidth]{"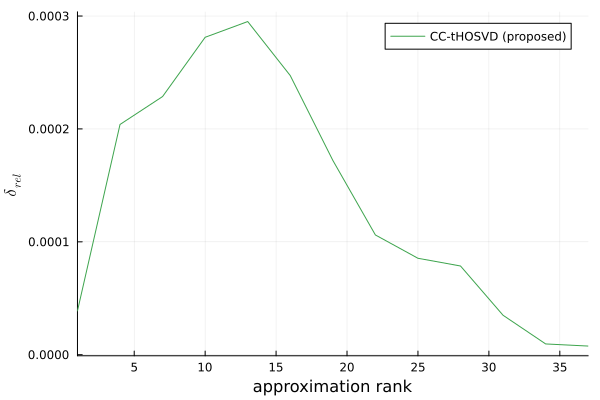"}
    \caption{}
    \label{fig:p-discr-P3-Camino}
    \end{subfigure}
    \caption{The progressions of the relative error (a) and the relative discrepancy (b) of \cref{alg:tangent-space-HOSVD,alg:curvature-corrected-HOSVD} for an increasing approximation rank of the $40\times 40$ $\mathcal{P}(3)$-valued tensor linearized at the data barycentre. (a) shows that CC-tHOSVD gruadually picks up on the global geometry as the relative error decreases, and in doing so attains better estimates for higher ranks. (b) shows that discrepancy between the curvature corrected approximation error and the exact global approximation error is very tight, supporting again the accuracy of \cref{eq:general-error-symmetric spaces}.}
    \label{fig:p-error-and-discr-P3-Camino}
\end{figure}

\begin{table}[h!]
    \centering
    \caption{The average run times in seconds of \cref{alg:tangent-space-HOSVD,alg:curvature-corrected-HOSVD} per multi-linear approximation rank of the $40\times 40$ Camino slice data on $\mathcal{P}(3)$ linearized at the data barycentre. The run time of CC-tHOSVD increases approximately quadratically with the rank.}
    \label{tab:p-run-times-P3-Camino}
    \begin{tabular}{lcccccc}
            \toprule
            Method & $\operatorname{rank} = 1$  & $\operatorname{rank} =8$  & $\operatorname{rank} = 15$  & $\operatorname{rank} = 22$  & $\operatorname{rank} = 29$ & $\operatorname{rank} = 36$   \\ 
            \midrule
            tHOSVD & $0.2888$ & $0.4062$ & $0.4438$ & $0.6294$ & $0.9385$ & $0.969$\\ 
            CC-tHOSVD (proposed) & $0.4683$ & $6.93$ & $24.92$ & $53.28$ & $99.06$ & $157.8$\\ 
            \bottomrule
    \end{tabular}
\end{table}

\paragraph{Low multi-linear rank approximation on $\mathcal{P}(3)^{10\times 10}$.}

The data set in \cref{fig:camino-data-patch} on the other hand behaves similarly to the 1D example in \cref{sec:P3-1D-tensor}. That is, for all approximation ranks the relative error reaches the zero-$\delta$ lower bound\footnote{The logarithmic-scaled plot in \cref{fig:p-log-errors-P3-Camino-patch} in \cref{app:extra-results} makes this easier to see.} and the relative discrepancy (\cref{fig:p-discr-P3-Camino-patch}) is very low at all times. Combined, these two observations indicate strong global-geometry awareness. Note that it is unsurprising that patch-wise compression of images works better than global compression as attempted above. That is, because the image's pixels tend to vary little on a local scale, but typically vary significantly when considered globally.

Regarding the run times in \cref{tab:p-run-times-P3-Camino-patch}, we observe quadratic scaling in rank once again, but note that this time the difference to tHOSVD is much more managable.

So, also for this data set we conclude that CC-tHOSVD is a good compromise between global-geometry awareness and computational feasibility. 

\begin{figure}[h!]
    \centering
    \begin{subfigure}{0.48\linewidth}
    \includegraphics[width=\linewidth]{"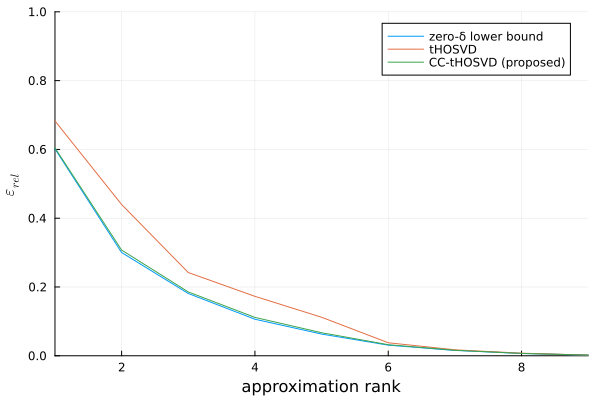"}
    \caption{}
    \label{fig:p-errors-P3-Camino-patch}
    \end{subfigure}
    \hfill
    \begin{subfigure}{0.48\linewidth}
    \includegraphics[width=\linewidth]{"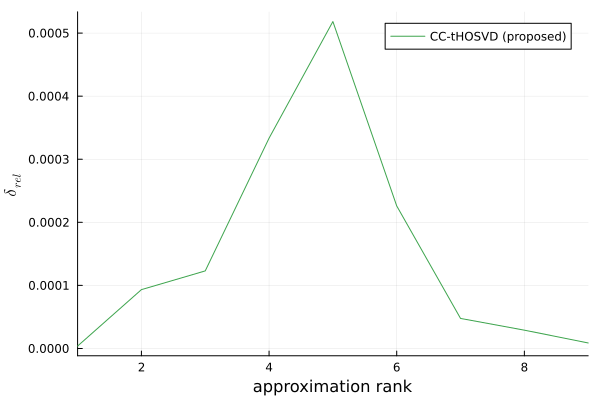"}
    \caption{}
    \label{fig:p-discr-P3-Camino-patch}
    \end{subfigure}
    \caption{The progressions of the relative error (a) and the relative discrepancy (b) \cref{alg:tangent-space-HOSVD,alg:curvature-corrected-HOSVD} for an increasing approximation rank of the $10\times 10$ $\mathcal{P}(3)$-valued tensor linearized at the data barycentre. (a) shows that CC-tHOSVD attains significantly lower errors through picking up on global geometry for all ranks. (b) shows that discrepancy between the curvature corrected approximation error and the exact global approximation error is very tight, supporting again the accuracy of \cref{eq:general-error-symmetric spaces}.}
    \label{fig:p-error-and-discr-P3-Camino-patch}
\end{figure}

\begin{table}[h!]
    \centering
    \caption{The average run times in seconds of \cref{alg:tangent-space-HOSVD,alg:curvature-corrected-HOSVD} per approximation rank of the $10\times 10$ Camino slice data on $\mathcal{P}(3)$ linearized at the data barycentre. The run time of CC-tHOSVD increases approximately quadratically with the rank.}
    \label{tab:p-run-times-P3-Camino-patch}
    \begin{tabular}{lccccc}
            \toprule
            Method & $\operatorname{rank} = 1$  & $\operatorname{rank} = 3$  & $\operatorname{rank} = 5$  & $\operatorname{rank} = 7$  & $\operatorname{rank} = 9$    \\ 
            \midrule
            tHOSVD & $0.00768$ & $0.009575$ & $0.00896$ & $0.00959$ & $0.010315$\\ 
            CC-tHOSVD (proposed) & $0.02036$ & $0.08746$ & $0.1853$ & $0.3186$ & $0.5483$\\ 
            \bottomrule
    \end{tabular}
\end{table}

\section{Conclusions.} 
\label{sec:manifold-valued-tensors-conclusions}

In this work we set out to construct a general framework for tangent space-based approximation of manifold-valued tensors that is applicable to general Riemannian manifolds, global-geometry aware and computationally feasible. 

\paragraph{General applicability.}
Through focusing on symmetric Riemannian manifolds -- covering a large class of manifolds that are being used in many practical settings -- we were able to show that the leading term of the global approximation error \cref{eq:manifold-approximation-error} -- the curvature corrected approximation error -- is fully determined by the tangent space approximation error and the curvature at the point of linearization. 

\paragraph{Global-geometry awareness.}
Subsequently, we showed that minimizing this curvature corrected approximation error -- already giving us general applicability -- instead of the global approximation error also gives us global-geometry awareness. In particular, we showed that under mild conditions minimizing 
the curvature corrected approximation error diminishes the discrepancy with the global approximation error. In other words, minimizers of the curvature corrected approximation error based methods will have a low global approximation error. Moreover and under mild conditions on the set of admissible approximators, we showed that as curvature effects vanish, minimizers of the curvature corrected approximation error converge to the minimizer of the zero-curvature case, i.e., the naive uncorrected tangent space-based approach. 

\paragraph{Computational feasibility.}
These latter observations were a strong motivation to consider numerical schemes that are initialized from naive tangent space-based methods and that are corrected for curvature from thereon. In particular, in a case study on low multi-linear rank approximation we used such a strategy and proposed the curvature corrected truncated higher-order singular value decomposition (CC-tHOSVD) scheme along with a generalization of the (truncated) HOSVD and metric-corrected tHOSVD (MC-tHOSVD) as a baseline approach that minimizes the global approximation error for benchmarking. The CC-tHOSVD -- already considerate of general applicability and global-geometry awareness by construction -- was constructed to account for computational feasibility through its main step being to solve a linear system, whose matrix representation is symmetric positive semi-definite and is positive definite under mild conditions.

Finally, in numerical experiments on synthetic and real tensor-valued data living in symmetric Riemannian manifolds that cover the cases of both positive and negative curvature, we tested the algorithm's performance and conclude that CC-tHOSVD is indeed generally applicable, global-geometry aware and computationally feasible.


\section*{Acknowledgments.}
Part of this research was performed while WD and DN were visiting the Institute for Pure and Applied Mathematics (IPAM), which is supported by the National Science Foundation (Grant No. DMS- 1925919). DN was partially supported by NSF-DMS 2011140.

\bibliographystyle{plain}
\bibliography{references} 

\clearpage
\appendix

\section{Supplementary results from \cref{sec:manifold-valued-tensors-numerics}.}
\label{app:extra-results}

\begin{figure}[h!]
    \centering
    \begin{subfigure}{0.48\linewidth}
    \includegraphics[width=\linewidth]{"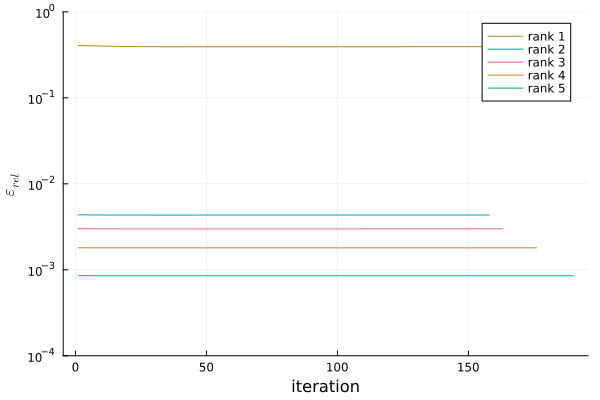"}
    \caption{ }
    \label{fig:p1-convergence}
    \end{subfigure}
    \hfill
    \begin{subfigure}{0.48\linewidth}
    \includegraphics[width=\linewidth]{"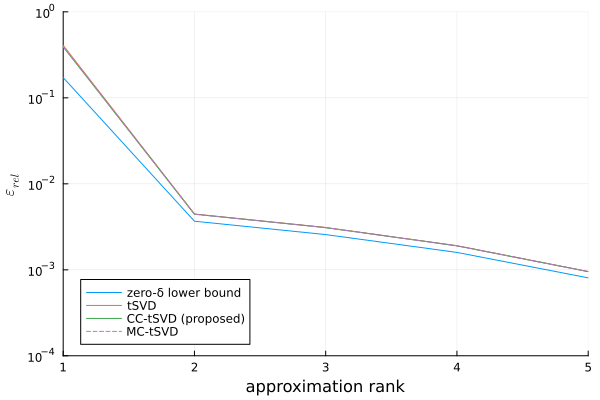"}
    \caption{ }
    \label{fig:p1-log-errors}
    \end{subfigure}
    \caption{The progressions of the relative error (a) of \cref{alg:geometry-corrected-HOSVD} with different approximation rank until convergence and (b) of \cref{alg:tangent-space-HOSVD,alg:curvature-corrected-HOSVD,alg:geometry-corrected-HOSVD} for an increasing approximation rank of the 1D data on $\Sphere^6$ from \cref{eq:1D-S6-data-set}, linearized at $\mPoint_1$ from \cref{eq:p1-p2}. Both y-axes are scaled logarithmically.}
    \label{fig:p1-extras}
\end{figure}

\begin{figure}[h!]
    \centering
    \begin{subfigure}{0.48\linewidth}
    \includegraphics[width=\linewidth]{"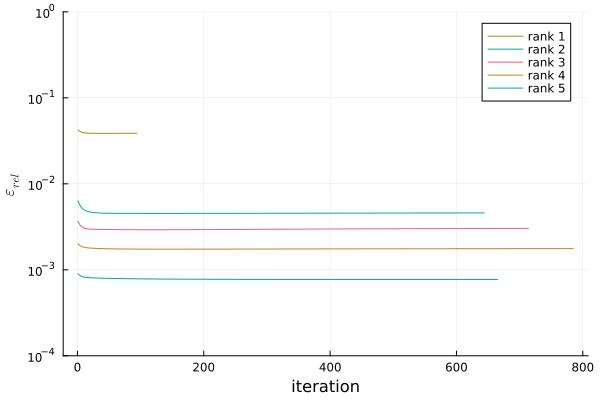"}
    \caption{ }
    \label{fig:p2-convergence}
    \end{subfigure}
    \hfill
    \begin{subfigure}{0.48\linewidth}
    \includegraphics[width=\linewidth]{"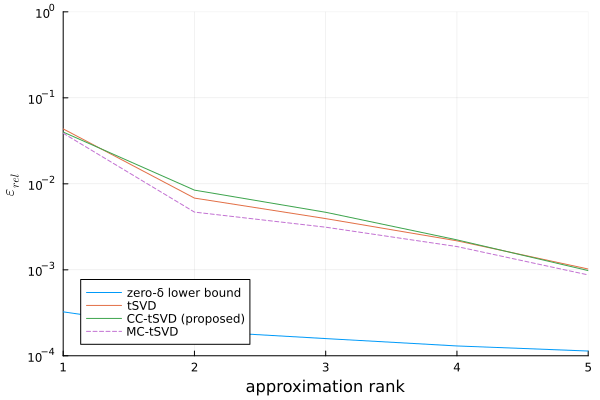"}
    \caption{ }
    \label{fig:p2-log-errors}
    \end{subfigure}
    \caption{The progressions of the relative error (a) of \cref{alg:geometry-corrected-HOSVD} with different approximation rank until convergence and (b) of \cref{alg:tangent-space-HOSVD,alg:curvature-corrected-HOSVD,alg:geometry-corrected-HOSVD} for an increasing approximation rank of the 1D data on $\Sphere^6$ from \cref{eq:1D-S6-data-set}, linearized at $\mPoint_2$ from \cref{eq:p1-p2}. Both y-axes are scaled logarithmically.}
    \label{fig:p2-extras}
\end{figure}

\begin{figure}[h!]
    \centering
    \begin{subfigure}{0.48\linewidth}
    \includegraphics[width=\linewidth]{"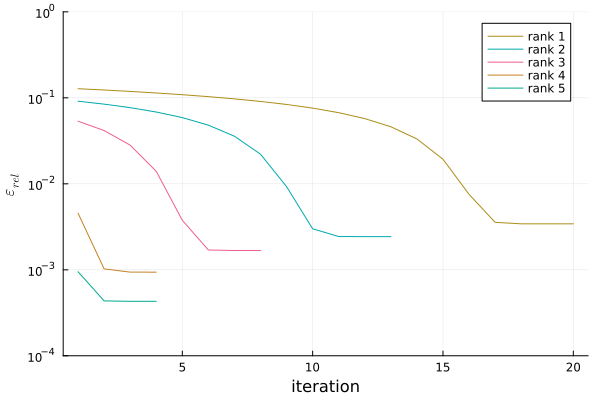"}
    \caption{ }
    \label{fig:p-convergence-P3}
    \end{subfigure}
    \hfill
    \begin{subfigure}{0.48\linewidth}
    \includegraphics[width=\linewidth]{"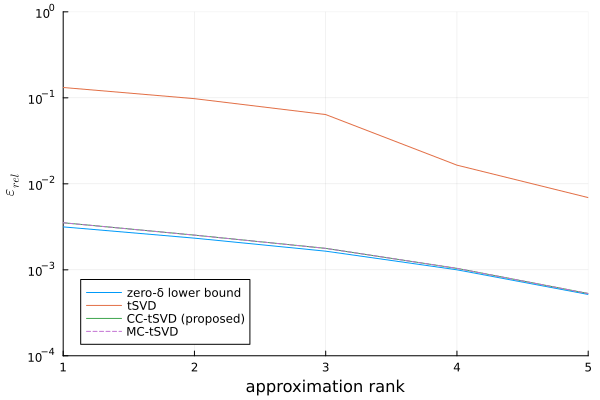"}
    \caption{ }
    \label{fig:p-log-errors-P3}
    \end{subfigure}
    \caption{The progressions of the relative error (a) of \cref{alg:geometry-corrected-HOSVD} with different approximation rank until convergence and (b) of \cref{alg:tangent-space-HOSVD,alg:curvature-corrected-HOSVD,alg:geometry-corrected-HOSVD} for an increasing approximation rank of the 1D data on $\mathcal{P}(3)$ from \cref{eq:1D-P3-data-set}, linearized at $\mPoint$ from \cref{eq:p-P3}. Both y-axes are scaled logarithmically.}
    \label{fig:p-extras-P3}
\end{figure}

\begin{figure}[h!]
    \centering
    \begin{subfigure}{0.48\linewidth}
    \includegraphics[width=\linewidth]{"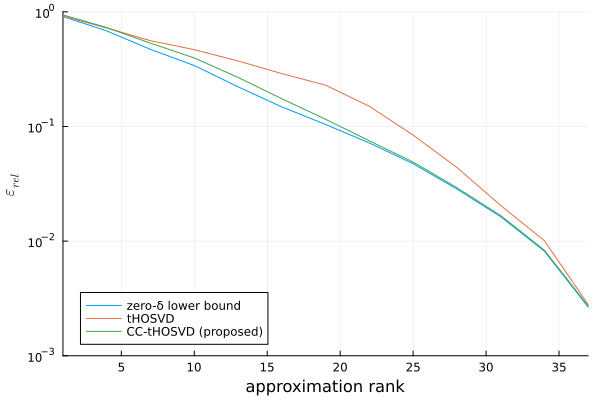"}
    \caption{ }
    \label{fig:p-log-errors-P3-Camino}
    \end{subfigure}
    \hfill
    \begin{subfigure}{0.48\linewidth}
    \includegraphics[width=\linewidth]{"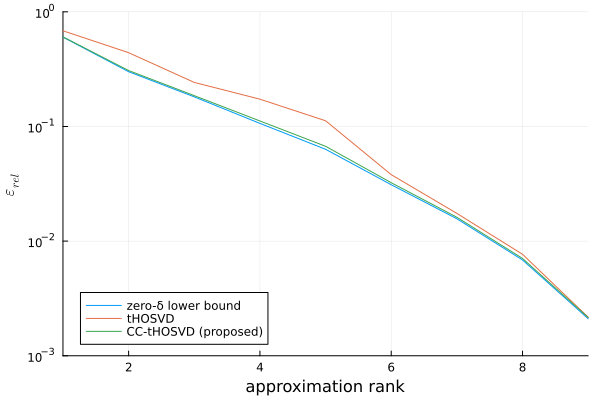"}
    \caption{ }
    \label{fig:p-log-errors-P3-Camino-patch}
    \end{subfigure}
    \caption{The progressions of the relative error of \cref{alg:tangent-space-HOSVD,alg:curvature-corrected-HOSVD} for an increasing approximation rank of (a) the $40\times 40$ $\mathcal{P}(3)$-valued tensor and of (b) the $10\times 10$ $\mathcal{P}(3)$-valued tensor from \cref{eq:2D-data-DTMRI}, both of which are linearized at the data barycentre. Both y-axes are scaled logarithmically.}
    \label{fig:p-extras-P3-Camino2D}
\end{figure}


\end{document}